\documentstyle[11pt, amstex]{amsart}

\makeatletter

\begin{document}

\bibliographystyle{amsalpha}

\newcommand{\e}{\epsilon}
\newcommand{\0}{{\bold 0}}
\newcommand{\w}{{\bold w}}
\newcommand{\y}{{\bold y}}
\newcommand{\balpha}{{\boldsymbol \alpha}}
\newcommand{\bt}{{\bold t}}
\newcommand{\z}{{\bold z}}
\newcommand{\x}{{\bold x}}
\newcommand{\N}{{\bold N}}
\newcommand{\Z}{{\bold Z}}
\newcommand{\F}{{\bold F}}
\newcommand{\R}{{\bold R}}
\newcommand{\Q}{{\bold Q}}
\newcommand{\C}{{\bold C}}
\newcommand{\BP}{{\bold P}}
\newcommand{\cO}{{\mathcal O}}
\newcommand{\cX}{{\mathcal X}}
\newcommand{\cH}{{\mathcal H}}
\newcommand{\cM}{{\mathcal M}}
\newcommand{\cD}{{\mathcal D}}
\newcommand{\cB}{{\mathcal B}}
\newcommand{\cC}{{\mathcal C}}
\newcommand{\cT}{{\mathcal T}}
\newcommand{\cI}{{\mathcal I}}
\newcommand{\cJ}{{\mathcal J}}
\newcommand{\cS}{{\mathcal S}}
\newcommand{\cF}{{\mathcal F}}
\newcommand{\cE}{{\mathcal E}}
\newcommand{\cP}{{\mathcal P}}
\newcommand{\cL}{{\mathcal L}}
\newcommand{\cY}{{\mathcal Y}}
\newcommand{\sB}{{\sf B}}
\newcommand{\sE}{{\sf E}}

\newcommand{\sA}{{\sf A}}
\newcommand{\ga}{{\sf a}}
\newcommand{\es}{{\sf s}}
\newcommand{\m}{{\bold m}}
\newcommand{\bS}{{\bold S}}

\newcommand{\ovf}{{\overline{f}}}

\newcommand{\ihra}{\stackrel{i}{\hookrightarrow}}
\newcommand\rank{\mathop{\rm rank}\nolimits}
\newcommand\im{\mathop{\rm Im}\nolimits}
\newcommand\coker{\mathop{\rm coker}\nolimits}
\newcommand\Li{\mathop{\rm Li}\nolimits}
\newcommand\NS{\mathop{\rm NS}\nolimits}
\newcommand\Hom{\mathop{\rm Hom}\nolimits}
\newcommand\Ext{\mathop{\rm Ext}\nolimits}
\newcommand\Pic{\mathop{\rm Pic}\nolimits}
\newcommand\Spec{\mathop{\rm Spec}\nolimits}
\newcommand\Hilb{\mathop{\rm Hilb}\nolimits}
\newcommand\Ker{\mathop{\rm Ker}\nolimits}
\newcommand{\length}{\mathop{\rm length}\nolimits}
\newcommand{\res}{\mathop{\sf res}\nolimits}

\newcommand\lra{\longrightarrow}
\newcommand\ra{\rightarrow}
\newcommand\la{\leftarrow}
\newcommand\JG{J_{\Gamma}}
\newcommand{\wvskp}{\vspace{1cm}}
\newcommand{\vskp}{\vspace{5mm}}
\newcommand{\nvskp}{\vspace{1mm}}
\newcommand{\nid}{\noindent}
\newcommand{\new}{\nvskp \nid}
\newtheorem{Assumption}{Assumption}[section]
\newtheorem{Theorem}{Theorem}[section]
\newtheorem{Lemma}{Lemma}[section]
\newtheorem{Remark}{Remark}[section]
\newtheorem{Corollary}{Corollary}[section]
\newtheorem{Conjecture}{Conjecture}[section]
\newtheorem{Proposition}{Proposition}[section]
\newtheorem{Example}{Example}[section]
\newtheorem{Definition}{Definition}[section]
\newtheorem{Question}{Question}[section]
\newtheorem{Fact}{Fact}[section]

\baselineskip=14pt

\title{Nodal curves and Riccati solutions of Painlev\'e equations}
\author{Masa-Hiko Saito and Hitomi Terajima}
\address{Department of Mathematics, Faculty of Science, 
Kobe University, Kobe, Rokko, 657-8501, Japan}
\email{mhsaito@math.kobe-u.ac.jp}
\email{terajima@math.kobe-u.ac.jp}
\thanks{Partly supported by Grant-in Aid
for Scientific Research  (B-12440008) and  (Houga-138740023)
the Ministry of Education, Science and Culture, Japan }
%\date{January, 23, 2002}

\begin{abstract}
In this paper, we study Riccati solutions of Painlev\'e equations 
from a view point of geometry of 
Okamoto-Painlev\'e pairs $(S,Y)$. 
After establishing  the correspondence between 
(rational) nodal curves on $S-Y$ and  Riccati solutions, 
we give the complete classification of  the configurations of 
nodal curves on $S-Y$ for each Okamoto--Painlev\'e pair $(S, Y)$. 
As an application of the classification, we prove 
the non-existence of 
Riccati solutions of Painlev\'e equations of 
types $P_{I}, P_{III}^{\tilde{D}_8}$ 
and $P_{III}^{\tilde{D}_7}$. We will also give  
a partial answer to the conjecture in \cite{STT} and \cite{T} that the 
dimension of the local cohomology 
$H^1_{Y_{red}}(S,\Theta_S(-\log Y_{red}))$ is one.   
\end{abstract}

\maketitle

\section{Introduction}

A pair   $(S, Y)$   of a projective 
smooth surface $S$ and an  effective anti-canonical divisor  
$Y$ on $S$ is called an Okamoto--Painlev\'e pair if it satisfies 
a suitable condition (see (\ref{eq:canonical}) in \S \ref{section:riccati}).  
In \cite{STT}, we established the theory of Okamoto--Painlev\'e 
pairs $(S, Y)$ and characterize 
the Painlev\'e equations by means of the special deformation 
of Okamoto--Painlev\'e pairs.  There exist 8 types of rational 
Okamoto--Painlev\'e pairs which correspond to the Painlev\'e equations. 
The types are classified by the types of the dual graphs of 
the configurations of $Y$, 
which are the affine Dynkin diagram of types 
$R = \tilde{D_k}, 4 \leq k \leq 8,  \tilde{E_l}, 6 \leq l \leq 8$.  
For each $R$, we obtain the global family of Okamoto--Painlev\'e pairs  
\begin{equation}\label{eq:globalf-intro} 
\begin{array}{ccl}
 {\cal S}_R &  \hookleftarrow & {\cal D}      \\
\pi \downarrow \hspace{0.3cm} &    \swarrow &  \varphi  \\
 \cM_R \times \cB_R  &  &  
\end{array} 
\end{equation}
where $\cB_R$ is an affine open subset of the $t$-affine line $\Spec \C[t]$. 
In \cite{STT}, the deformation with respect to the $t$-direction 
can be characterized by the local cohomology 
group  $H^1_{D}(S, \Theta_{S} ( - \log D))$ where $D = Y_{red}$.   
Furthermore we can show that
the vector field $\frac{\partial}{\partial t}$ has a unique lifting to 
a rational global vector field 
\begin{equation}\label{eq:gf-intro}
\tilde{v} \in H^0(\cS_R, \Theta_{\cS_R} ( - \log \cD) 
\otimes \cO_{\cS}(\cD))
\end{equation}
which  induces  the Painlev\'e differential equations on $\cS_R - \cD$.  
 
 In the theory of Painlev\'e equations, it is important to determine
 all classical solutions, like rational solutions and Riccati solutions.  
 (For the definition of 
 classical solutions of Painlev\'e equations, see \S 1 in \cite{U-W1}). 
In this direction, there are  a considerable number of works by many authors.  
Here, we  list up only a part of the references:
(e.g.,  \cite{DM}, \cite{Grm1}, \cite{Grm2}, 
\cite{Grm3}, \cite{Grm4}, \cite{Grm5}, \cite{Grm6}, 
\cite{Gr-Lu},  \cite{Gr-Ts}, \cite{Luk1}, 
\cite{Luk2}, \cite{Maz},  \cite{Mu1}, \cite{Mu2}, \cite{Ohy}, 
\cite{O3}, \cite{U1}, \cite{U2}, 
\cite{U-W1}, \cite{U-W2}, \cite{V}, \cite{W1}, 
\cite{W2}, \cite{Y}).   
 For example, in order to prove the irreducibility of the Painlev\'e equations, 
 one has to determine  the cases when the given Painlev\'e 
 equations admit the  Riccati solutions 
 (cf. \cite{U1}, \cite{U2}, \cite{NO}, \cite{U-W1}, \cite{U-W2}).

 \vspace{0.2cm}
 
 One of the main purpose of this paper is 
 to  characterize the Riccati solutions of Painlev\'e equations 
 by means of geometry of nodal curves on $S - Y_{red}$ 
 for the corresponding rational  Okamoto--Painlev\'e pairs $(S, Y)$. 
 Since our charcterization of Painlev\'e vector field 
 $\tilde{v}$ (\ref{eq:gf-intro}) in \cite{STT} is intrinsic, that is, coordinate free, 
 so is the characterization of Riccati solutions.  
 
  Moreover we shall give the complete classification 
theorem (Theorem \ref{thm:main-2}) of configurations of nodal curves on $S - Y_{red}$ 
for all rational Okamoto--Painlev\'e pairs $(S, Y)$ of  non-fibered type and of additive type. 
   As a corollary to Theorem (\ref{thm:main-2}), we can show that 
   Painlev\'e equations $P_{I}, P_{III}^{\tilde{D_7}}, P_{III}^{\tilde{D_8}}$ have no
   Riccati solutions for any parameters in the equations.

 The following is a rough outline of this paper.

 In \S \ref{section:riccati},  we characterize the Riccati 
 solutions of the Painlev\'e equations by means of 
 $(-2)$-curve (or nodal curve) $C$  on $S - Y_{red}$.   
 If for a given $(\balpha_0,t_0) 
 \in \cM_R \times \cB_R$, the fiber $S$ 
 of $\pi$ in (\ref{eq:globalf-intro}) over $(\balpha_0,t_0)$ 
 contains a nodal curve 
 $C \subset S - Y_{red}$, we can extend the nodal curve $C$ 
 in the $t$-direction and obtain 
 a family of nodal curves  $\cC \lra \{ \balpha_0 \} \times U$ 
 where $U$ is an (analytic or \'etale) open 
 neighborhood of $t_0$ in $\cB_R$.  Then the restriction $\tilde{v}_{|\cC}$ is 
 tangent to  $\cC$ which induces the Riccati equation on $\cC$.  
 It seems that this approach is  essentially equivalent to  
 Umemura's theory of invariant divisors for  
 the Painlev\'e equations (cf. e.g., \cite{U-W1}).  
 However we believe that  our approach gives a  clearer 
  geometric viewpoint of Riccati solutions of Painlev\'e equations.

In \S \ref{section:class}, 
we shall give the complete classification of configurations of 
nodal curves on $S - Y$ for all rational Okamoto--Painlev\'e 
pairs $(S, Y)$ of  non-fibered type and of additive type. 
The classification is based on the structure theorem of the lattice 
 induced by the intersection form on $H^2(S, \Z)$.  
 We can show that 
the sub-lattice generated by the nodal curves 
$C$ on $S - Y$ is a sub-lattice of $E_8^-$, the unique 
even unimodular negative-definite  lattice of rank $8$.  Then taking 
account into the sub-lattice generated by 
the irreducible components of $Y$, we can 
obtain the list of the possible configurations. 
For the existence of the possible configurations, we quote 
the Oguiso--Shioda's classification  theorem of singular fiber or 
Mordell--Weil group for rational elliptic surfaces.    Note that a rational 
elliptic surface with a fixed fiber is a  rational Okamoto--Painlev\'e 
pairs of fibered type in our terminology.  
Using  the Oguiso--Shioda's existence theorem 
and  the deformation theory of Okamoto--Painlev\'e pairs, we shall show 
the existence of all  possible configurations for some rational 
Okamoto--Painlev\'e pairs.  

In \S \ref{sec:non-exist}, as a corollary to the classification theorem, 
 we shall prove the non-existence of Riccati solutions 
of  the Painlev\'e equations of 
type $R = P_I, P_{III}^{\tilde{D_7}}, P_{III}^{\tilde{D_8}}$.
Though there are other proofs for this  result for $P_I$ and $P_{III}^{\tilde{D}_7}$  
 (e.g.,  \cite{U1}, \cite{U2} and \cite{Ohy}),  
our proof  clarify the point that the obstruction to  
the existence of Riccati solutions lies in the  topological conditions.  

In \S \ref{sec:example}, 
we give explicit examples of nodal curves and Riccati solutions of 
Painlev\'e equations associated to the nodal curves.  

In \S \ref{sec:conf}, we shall give an example of the 
confluence of the Riccati solutions for 
$R = \tilde{E_6}, (P_{IV})$ and also the confluence of nodal curves.  
Moreover 
we give a remark on rational solutions 
coming from the intersection of two different 
Riccati solutions.  

In Appendix \ref{sec:lcoh}, as a 
corollary to  Theorem \ref{thm:main-2}, 
we shall give a partial answer to the Conjecture \ref{conj:local} 
presented in \cite{STT} and \cite{T} about 
the dimension of the local cohomology group.

\section{\bf $(-2)$-curves (nodal curves) and Riccati solutions}
\label{section:riccati}

In this section, we shall review the theory of Okamoto--Painlev\'e pairs 
and their  relations to the Painlev\'e equations which were introduced 
in \cite{STT}.

\subsection{Okamoto--Painlev\'e pairs}

\begin{Definition}\label{def:op}
{\rm
Let $(S, Y)$ be a pair  of a  complex projective surface $S$  and an effective 
anti-canonical divisor $Y \in |-K_S|$ of $S$. Let $Y= \sum_{i=1}^r m_i Y_i$ 
be the irreducible decomposition of $Y$.    
 We call a pair $(S, Y)$ an {\em  Okamoto--Painlev\'e pair} if for all $i, 1 \leq i \leq r$, 
\begin{equation} \label{eq:canonical}
Y \cdot Y_i = \deg [Y]_{|Y_i} = 0.
\end{equation}
An Okamoto--Painlev\'e pair $(S, Y)$ is called {\em rational} if $S$ is a rational surface.
}
\end{Definition}
             
\begin{Remark} {\rm An  Okamoto--Painlev\'e pair $(S, Y)$ 
in Definition \ref{def:op} is called a generalized Okamoto--Painlev\'e pair 
in \cite{STT}.  However,  in this paper, we shall use this terminology.  
Note that in the original definition of an Okamoto--Painlev\'e pair $(S, Y)$ in \cite{STT} 
 we assume that $S - Y_{red}$ contains $\C^2$ as a Zariski 
 open set and $Y_{red}$ is a normal crossing divisor. (See also \cite{Sa-Ta}.)}
 \end{Remark}

\subsection{Okamoto--Painlev\'e pairs and Painlev\'e equations}

Let $(S,Y)$ be a rational Okamoto--Painlev\'e pair with the irreducible 
decomposition $Y= \sum_{i=1}^r m_i Y_i$ and set $D=Y_{red} = \sum_{i=1}^r  Y_i$. 
Denote by $M(Y)$ the sub-lattice of $\Pic(S) \simeq H^2(S, \Z)$ generated by the 
irreducible components $\{ Y_i \}_{i=1}^r$.  With the bilinear form on $M(Y)$ 
 which is $(-1)$ times the intersection pairing on $S$, $M(Y)$ becomes a root lattice 
of affine type (cf. [Section 1, \cite{STT}], \cite{Sakai}).  Let $R(Y)$ denote the 
type of the root lattice.  
One can classify rational Okamoto--Painlev\'e pairs $(S, Y)$ in terms of the type $R(Y)$. 
  (See [Section 1, \cite{STT}], \cite{Sakai}).

Not all types of rational Okamoto--Painlev\'e pairs correspond to the
 Painlev\'e equations.  The Table \ref{tab:type} is the list of the types 
of Okamoto--Painlev\'e pairs  which correspond  to  the Painlev\'e equations.  
We shall explain the meaning of the correspondence  in Theorem \ref{thm:global}.  
Note that classically, Painlev\'e equations were classified into 6 types, 
however now we should classify them into 8 types.  
Actually, the third Painlev\'e equations $P_{III}$ can be classified further  into 3 types 
$P^{\tilde{D}_6}_{III}, P^{\tilde{D}_7}_{III}$ and $ P^{\tilde{D}_8}_{III}$ 
corresponding to the types of $R = R(Y)$. 
The classical third  Painlev\'e  equations correspond to $P^{\tilde{D}_6}_{III}$, which form 
a two parameter family of equations. The equations $P^{\tilde{D}_7}_{III}$ and $P^{\tilde{D}_8}_{III}$ 
can be obtained by specializations of these parameters.  

\begin{table}[h]
\begin{center}
{\bf  Okamoto--Painlev\'e pairs and Painlev\'e equations}

\vspace{0.2cm}
\begin{tabular}{||c||c|c|c|c|c|c|c|c|c||} \hline
    &  & & &  & & & &  \\
$R= R(Y)$ & $\tilde{E_8}$ &  $\tilde{E_7}$& $\tilde{D_8}$   & $\tilde{D_7}$  & 
$\tilde{D_6}$ & $\tilde{E_6}$ & $\tilde{D_5}$ & $\tilde{D_4}$ \\
    &  & & & & & & &  \\ \hline
    &  & & & & & & &  \\
Painlev\'{e} equation & $P_{I}$ & $P_{II}$ &  
$P^{\tilde{D}_8}_{III}$   & $P^{\tilde{D}_7}_{III}$ & $P^{\tilde{D}_6}_{III}$ & $P_{IV}$
& $P_{V}$  &  $P_{VI}$ \\
    &   & & & & & & & \\ \hline
\end{tabular}
\vspace{0.3cm}
\caption{}
\label{tab:type}
\end{center}
\end{table}

Here we shall recall one more important definition (cf. [Section 1, \cite{STT}]).
\begin{Definition}{\rm 
A rational Okamoto--Painlev\'e pair $(S, Y)$ will 
be called {\em of fibered type} if there exists an elliptic 
fibration $f:S \lra \BP^1$ such that 
$f^{*}(\infty)  = Y$ as divisors.  We say that a rational Okamoto--Painlev\'e 
pair is {\em of non-fibered type} if $(S, Y)$ is not of fibered type.  
}
\end{Definition}

The following theorem (cf.\ [Proposition 5.1, Theorem 6.1,  \cite{STT}])  
explains how 
one can give correspondences between   rational Okamoto--Painlev\'e pairs 
and Painlev\'e equations in Table \ref{tab:type}.  

\begin{Theorem}[Proposition 5.1. \cite{STT}]\label{thm:global}
Let $R = R(Y)$ be one of types of the root systems 
in Table \ref{tab:type} 
$( i.e., R = \tilde{D}_{i}, 4 \leq i \leq 8$ 
or $\tilde{E}_{j}, 6 \leq j \leq 8)$ and let $r$ be the number of 
irreducible components 
of $D = Y_{red}$ and set $s = s(R) = 9 -r$.   
Then there exist   affine open subschemes $ \cM_{R}  \subset \C^s = 
\Spec \C[ \alpha_1, \cdots, \alpha_s]$,   
$\cB_R \subset \C = \Spec \C[t]$, and 
the following commutative diagram satisfying the 
conditions below:  
\begin{equation}\label{eq:globalf} 
\begin{array}{ccl}
 {\cal S} &  \hookleftarrow & {\cal D}      \\
\pi \downarrow \hspace{0.3cm} &    \swarrow &  \varphi  \\
 \cM_R \times \cB_R  &  &.  
\end{array} 
\end{equation}
\begin{enumerate}

\item $\cal S$ is a smooth quasi-projective manifold and $\cD$ is a 
divisor with normal crossing of $\cS$.  Moreover $\pi$ is a smooth and 
projective morphism  and $\varphi$ is a flat morphism such that the above 
diagram is a deformation of non-singular pairs of projective surfaces and 
normal crossing divisors in the sense of Kawamata \cite{Kaw}.

\item There is a rational relative 2-form 
\begin{equation}\label{eq:2form}
\omega_{\cal S} \in \Gamma({\cal S}, \Omega^2_{{\cal S}/\cM_R \times \cB_R}(*{\cal D}) ) 
\end{equation}
which has poles only along ${\cal D}$.  If we denote by ${\cal Y}$ the pole divisor 
of $\omega_{\cal S}$, then for each point 
$(\balpha, t) \in \cM_R \times \cB_R$,  $({\cal S}_{\balpha, t}, 
{\cal Y}_{\balpha, t})$ is a rational 
Okamoto--Painlev\'e pair of type $R = R(Y)$ and ${\cal Y}_{red} = {\cal D}$.

\item There is a unique global rational vector field 
\begin{equation}\label{eq:vf}
\tilde{v} \in \Gamma(\cS, \Theta_{\cS}( - \log \cD) \otimes \cO_{\cS}(\cD))
\end{equation}
on $\cS$ which is a lift of  $\frac{\partial}{\partial t}$, 
that is, $\pi_*(\tilde{v}) =   \frac{\partial}{\partial t}$.  
Moreover the restriction of 
$\tilde{v}$ to $\cS - \cD $ gives a regular 
algebraic vector field which 
corresponds to the  Painlev\'e equation of type $R$. 
   We call the systems of differential equations 
   determined by the vector field $\tilde v$  
   {\em the Painlev\'e system of type $R$}.  
  $($See $($\ref{eq:diffeq}$)$ below$)$.
\end{enumerate}
\end{Theorem}

We can state more about the family in (\ref{eq:globalf}) as follows.  

\begin{enumerate}

\item The family is semi-universal at each point 
$(\balpha, t) \in \cM_R \times \cB_R$, 
that is, the Kodaira--Spencer map 
\begin{equation}
\rho:T_{\balpha, t}(\cM_R \times \cB_R) 
\lra H^1({\cal S}_{\balpha, t}, 
\Theta_{{\cal S}_{\balpha, t}}(- \log {\cal D}_{\balpha, t}))
\end{equation}
is an isomorphism. For a point  $(\balpha, t) \in 
\cM_R \times \cB_R$ at which the corresponding 
Okamoto--Painlev\'e pair is of non-fibered type, 
one can obtain the following commutative diagram:
\small
\begin{equation}\label{eq:com-diag1}
\begin{array}{cccc} 
  &       &       &  0 \\
  &       &       &  \uparrow  \\
 0 \longrightarrow &  
 H^1_{{\cal D}_{\balpha, t}}({\cal S}_{\balpha, t},  
 \Theta_{{\cal S}_{\balpha, t}}(- \log {\cal D}_{\balpha, t})) & \rightarrow & 
  H^1({\cal S}_{\balpha, t}, \Theta_{{\cal S}_{\balpha, t}}(- \log {\cal D}_{\balpha, t}))  \\ 
        &  \uparrow    &    &  ||     \\ 
 0 \longrightarrow  &  H^0( {\cal D}_{\balpha, t},  \Theta_{\cS_{\balpha, t}}( - \log \cD) 
     \otimes N_{{\cal D}_{\balpha, t}}) \simeq \C  
  \cdot \rho( \frac{\partial }{\partial t})  &   
  \hookrightarrow  &   
  H^1({\cal S}_{\balpha, t}, \Theta_{{\cal S}_{\balpha, t}}(- \log {\cal D}_{\balpha, t}))    \\ 
     &  \wr \uparrow    &    &  \quad \wr \uparrow  \rho       \\ 
    0  \longrightarrow &   T_{\balpha, t}(\cB_R) \simeq \C \cdot
    \frac{\partial }{\partial t} \ &  \hookrightarrow     &         T_{\balpha, t}(\cM_R \times \cB_R) \\
      &   \uparrow    &     & \uparrow     \\
      &   0   &     & 0     
\end{array}
\end{equation}

\item  Let $M_R$ and $B_R$ denote the affine coordinate rings of 
$\cM_R$ and $\cB_R$ respectively so that $ \cM_R= \Spec M_R$ and $\cB_R = \Spec B_R $.  
(Note that $M_R$ and $B_R$ are obtained by some localizations of 
$ \C[\alpha_1, \cdots, \alpha_s] $ and $ \C[t]$ respectively).   

 There exists an affine open covering  
 $\{ \tilde{U}_i \}_{i=1}^{l+k}$ of ${\cal S}$ 
such that  for each $i$
\begin{equation}\label{eq:covering1}
  \tilde{U}_i \   \simeq \  
     \Spec \left( (M_R \otimes B_R)  
     [x_i, y_i, \frac{1}{f_i(x_i, y_i, \balpha, t)}] \right) \subset  
      \Spec \C[ \balpha, t, x_i, y_i] \simeq \C^{s+3} \simeq \C^{12 - r}.
\end{equation}
Here $f_i(x_i, y_i, \balpha, t)$ is a polynomial in 
 $(M_R \otimes B_R) [x_i, y_i ] $.   Moreover, 
 we may assume that $\cS -{\cal D}$ can be covered by 
 $\{ \tilde{U_i} \}_{i=1}^l$,  and for each $i$, the 
 restriction of the 
 rational $2$-form $\omega_{\cS}$ can be written as 
\begin{equation}
\omega_{{\cal S}|\tilde{U}_i} = 
\frac{dx_i \wedge dy_i}{f_i(x_i, y_i, \balpha, t)^{m_i}}.
\end{equation}

\item  By using the local  coordinates of $\cS - \cD$, 
the global rational  vector field 
$\tilde{v}$ on $\cS$  obtained in (\ref{eq:vf}) 
can be written on  each open set $\tilde{U_i}$ for 
 $ 1 \leq i \leq l $ (corresponding to the open  coverings of $\cS - \cD$)  as 
\begin{equation}\label{eq:global-v}
\tilde{v}_{| \tilde{U}_i}  \  =  \  
\frac{\partial}{\partial t} \  - \  \theta_i \  = \  
\frac{\partial}{\partial t} \  - \  
\eta_i \frac{\partial}{\partial x_i } - \zeta_i \frac{\partial}{\partial y_i}
\end{equation}
where $\theta_i = \eta_i \frac{\partial}{\partial x_i }
 + \zeta_i \frac{\partial}{\partial y_i}$
  is a regular algebraic vector field on $\tilde{U_i}$.  

This explicit expression of $\tilde{v}_{|\tilde{U_i}}$ 
gives  a system of differential equations 
\begin{equation}\label{eq:diffeq}
\renewcommand{\arraystretch}{2.2}
\left\{
\begin{array}{ccl}
\displaystyle{\frac{d x_i }{d t}} &  = & -  \eta_i (x_i, y_i, \balpha, t) \\
\displaystyle{\frac{d y_i }{d t}} & = & -\zeta_i (x_i, y_i,  \balpha, t) 
\end{array}
\right.
\end{equation}
which is equivalent to the Painlev\'e equation of type $R$.

\end{enumerate}

\begin{Remark}{\rm  One can show that the deformation corresponding to 
$\rho(\frac{\partial}{\partial t})$ preserves the relative rational 2-form
$\omega_{\cS}$ in (\ref{eq:2form}).  
This fact explains the reason why the systems of differential 
equations in (\ref{eq:diffeq}) can be written in Hamiltonian systems. 
For more details, see [\S 6, \cite{STT}]. } 
\end{Remark}

\subsection{Riccati equations}\label{ss:riccati}

Let $U \subset \C $ be an open complex domain (in analytic topology) 
with a local analytic coordinate $t$ and
 $ a(t), b(t), c(t) $  holomorphic functions defined in $U$.

Consider a 
Riccati equation 
\begin{equation} \label{eq:riccati}
 x' = a(t) x^2 + b(t) x + c(t).  
\end{equation}
By the change of unknown 
\begin{equation}\label{eq:riccati-sol}
x = - \frac{1}{a(t)}\frac{d}{dt} \log (u) = - \frac{1}{a(t)} \frac{u'}{u}, 
\end{equation}
the equation (\ref{eq:riccati}) is transformed into 
the linear equation
\begin{equation}\label{eq:riccati-lin}
u'' - [ \frac{a'(t)}{a(t)} + b(t) ]  u' + a(t) c(t) u = 0.
\end{equation}
Therefore the movable singularities  of the solution 
$x(t) =  - \frac{1}{a(t)} \frac{u'}{u}$ of 
(\ref{eq:riccati}) are only poles.  
This condition is called 
{\em the Painlev\'e property} for 
an algebraic ordinary differential equation.  (Cf. [3.1, Ch.\ 3, \cite{IKSY}]). 

\begin{Remark}\label{rm:riccati}
{\rm Riccati equations above are defined  
in the space $\BP^1 \times U $ with the  
coordinates $(x, t)$.  The equation (\ref{eq:riccati}) 
is  equivalent to a rational global vector field on $\BP^1 \times U $ as
\begin{equation}\label{eq:riccati-vf}
\tilde{v}= \frac{\partial}{\partial t} + 
[a(t) x^2 + b(t) x + c(t)] \frac{\partial}{\partial x}.  
\end{equation}
By the coordinate change 
$ u = \frac{1}{x} $, $\tilde{v}$ can be transformed into the form 
$$
\tilde{v} = \frac{\partial}{\partial t} - 
[a(t) + b(t) u + c(t) u^2] \frac{\partial}{\partial u}.  
$$
This shows that the vector field $\tilde{v}$ is holomorphic 
 even at  $x = \infty$, 
 hence $\tilde{v}$ is a global holomorphic vector field 
 on $\BP^1 \times U$.  (Conversely, one can show that any  
 holomorphic vector field on $\BP^1 
\times U$ which is a lift of $\frac{\partial}{\partial t}$ 
can be written as in (\ref{eq:riccati-vf})).   
Therefore the space $\BP^1 \times U $ can be considered as the space of 
 initial conditions for the Riccati equation above.  }
\end{Remark}

\subsection{Nodal curves on Okamoto--Painlev\'e pairs and Riccati equations}
Let $(S, Y)$ be a rational Okamoto--Painlev\'e pair of type $R=R(Y)$ 
corresponding to   Painlev\'e equations of type $R$.  Then,  as we see in 
Theorem \ref{thm:global},  one  can 
construct a global rational vector field $\tilde{v}$ 
on the semi-universal deformation 
family of $(S, Y)$ which gives  the Painlev\'e equation of type $R$.

In what follows, we will show that Painlev\'e equations  
can be reduced to the 
Riccati equations if and only if the corresponding rational Okamoto--Painlev\'e 
pair $(S, Y)$ contains $\BP^1$ on $ S - Y_{red}$. 
Roughly speaking,  we have the following correspondences.

\vspace{0.2cm}
\begin{equation}\label{eq:scheme}
\begin{array}{ccc}
\fbox{Painlev\'e  equations} & \Leftrightarrow & 
\fbox{Special deformations of Okamoto--
Painlev\'e pairs $(S, Y)$} \\
\cup  &    & \cup \\
\fbox{ Riccati equations}&  \Leftrightarrow  &      
\fbox{ Nodal curves $ C \simeq \BP^1   \subset S - Y_{red}$ }  
\end{array}
\end{equation}

\vspace{0.2cm}

In order to explain this scheme more explicitly, 
 let us consider the Hamiltonian systems of the 
Painlev\'e equation of type $\tilde{E_6}$ ($ = P_{IV}$) with 
two auxiliary parameters $\kappa_0, \kappa_{\infty}$.   

\begin{equation} \label{eq:pk0-1}
\left\{
\begin{array}{ccl}
\displaystyle{\frac{d x_0 }{d t}} &  = & 4x_0 y_0 -x_0^2 -2 t x_0 -2 \kappa_0 \\
   &    &    \\
\displaystyle{\frac{d y_0 }{d t}} & = & -2y_0^2 +2(x_0+t)y_0 -\kappa_\infty
\end{array}
\right..
\end{equation}
When $\kappa_0 = 0$, if we set  $x_0 \equiv 0$, the first equation of 
the  system (\ref{eq:pk0-1}) is  automatically satisfied, and 
the second equation can be reduced to the equation  
\begin{equation}
\displaystyle{\frac{d y_0 }{d t}}  =  -2y_0^2 +2ty_0 -\kappa_\infty, 
\end{equation}
which is nothing but a Riccati equation.  
One can easily check that $\{ x_0 =0 \}$ defines a smooth $\BP^1$ on $S - Y_{red}$. (See \S 4).   

Note that if $C \subset S - Y_{red}$ is a smooth irreducible 
rational curve in $S - Y_{red}$, we see that $K_S \cdot C = -Y \cdot C = 0$, 
hence, by the adjunction formula, we have    
$$
C^2 = K_S \cdot C + C^2 =  -2.
$$
Therefore a smooth irreducible rational curve $C \subset S - Y_{red}$ is always 
{\em a $(-2)$-curve or a nodal curve}.

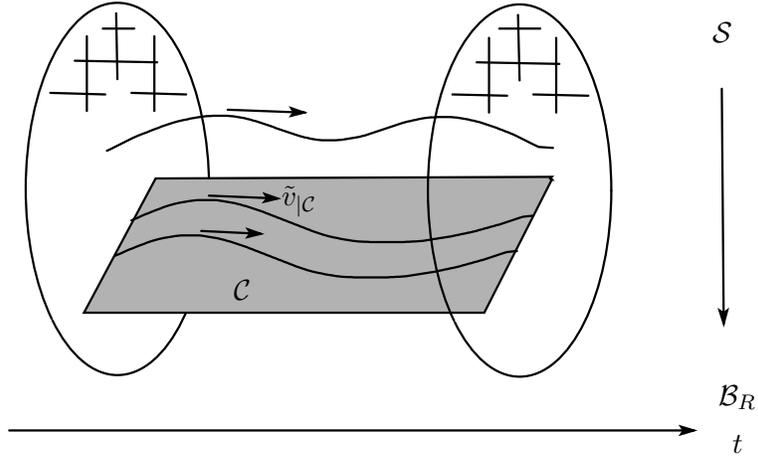
\begin{figure}[h]
\begin{center}
%WinTpicVersion2.15
\unitlength 0.1in
\begin{picture}(37.84,28.10)(13.40,-36.30)
% VECTOR 1 0 3 0
% 2 1340 3460 4924 3460
% 
\special{pn 13}%
\special{pa 1340 3060}%
\special{pa 4924 3060}%
\special{fp}%
\special{sh 1}%
\special{pa 4924 3060}%
\special{pa 4857 3040}%
\special{pa 4871 3060}%
\special{pa 4857 3080}%
\special{pa 4924 3060}%
\special{fp}%
% ELLIPSE 1 0 3 0
% 4 1908 2196 1428 3172 1428 3172 1364 3276
% 
\special{pn 13}%
\special{ar 1908 1796 480 976  2.3682869 6.2831853}%
\special{ar 1908 1796 480 976  0.0000000 2.3561945}%
% POLYGON 1 0 1 0
% 10 2108 2140 1732 2844 3828 2844 3828 2844 4180 2148 4172 2148 4180 2148 4172 2140 4180 2132 2108 2140
% 
\special{pn 13}%
\special{sh 0.300}%
\special{pa 2108 1740}%
\special{pa 1732 2444}%
\special{pa 3828 2444}%
\special{pa 3828 2444}%
\special{pa 4180 1748}%
\special{pa 4172 1748}%
\special{pa 4180 1748}%
\special{pa 4172 1740}%
\special{pa 4180 1732}%
\special{pa 2108 1740}%
\special{fp}%
% LINE 1 0 3 0
% 2 1684 1524 2116 1532
% 
\special{pn 13}%
\special{pa 1684 1124}%
\special{pa 2116 1132}%
\special{fp}%
% LINE 1 0 3 0
% 2 1900 1276 1908 1620
% 
\special{pn 13}%
\special{pa 1900 876}%
\special{pa 1908 1220}%
\special{fp}%
% LINE 1 0 3 0
% 2 1748 1396 1748 1788
% 
\special{pn 13}%
\special{pa 1748 996}%
\special{pa 1748 1388}%
\special{fp}%
% LINE 1 0 3 0
% 2 1556 1692 1844 1700
% 
\special{pn 13}%
\special{pa 1556 1292}%
\special{pa 1844 1300}%
\special{fp}%
% LINE 1 0 3 0
% 2 2100 1396 2100 1788
% 
\special{pn 13}%
\special{pa 2100 996}%
\special{pa 2100 1388}%
\special{fp}%
% LINE 1 0 3 0
% 2 1980 1692 2268 1700
% 
\special{pn 13}%
\special{pa 1980 1292}%
\special{pa 2268 1300}%
\special{fp}%
% ELLIPSE 1 0 3 0
% 4 4012 2204 3532 3180 3532 3180 3468 3284
% 
\special{pn 13}%
\special{ar 4012 1804 480 976  2.3682869 6.2831853}%
\special{ar 4012 1804 480 976  0.0000000 2.3561945}%
% LINE 1 0 3 0
% 2 3788 1532 4220 1540
% 
\special{pn 13}%
\special{pa 3788 1132}%
\special{pa 4220 1140}%
\special{fp}%
% LINE 1 0 3 0
% 2 4004 1284 4012 1628
% 
\special{pn 13}%
\special{pa 4004 884}%
\special{pa 4012 1228}%
\special{fp}%
% LINE 1 0 3 0
% 2 3852 1404 3852 1796
% 
\special{pn 13}%
\special{pa 3852 1004}%
\special{pa 3852 1396}%
\special{fp}%
% LINE 1 0 3 0
% 2 3660 1700 3948 1708
% 
\special{pn 13}%
\special{pa 3660 1300}%
\special{pa 3948 1308}%
\special{fp}%
% LINE 1 0 3 0
% 2 4204 1404 4204 1796
% 
\special{pn 13}%
\special{pa 4204 1004}%
\special{pa 4204 1396}%
\special{fp}%
% LINE 1 0 3 0
% 2 4084 1700 4372 1708
% 
\special{pn 13}%
\special{pa 4084 1300}%
\special{pa 4372 1308}%
\special{fp}%
% STR 2 0 3 0
% 3 5020 1348 5020 1428 2 0
% $\cS$
\put(50.2000,-10.2800){\makebox(0,0)[lb]{$\cS$}}%
% SPLINE 1 0 3 0
% 6 1980 2360 2324 2240 3052 2456 4060 2344 4084 2336 4084 2336
% 
\special{pn 13}%
\special{pa 1980 1960}%
\special{pa 2010 1946}%
\special{pa 2040 1931}%
\special{pa 2071 1918}%
\special{pa 2101 1904}%
\special{pa 2131 1892}%
\special{pa 2161 1880}%
\special{pa 2191 1869}%
\special{pa 2222 1860}%
\special{pa 2252 1852}%
\special{pa 2282 1846}%
\special{pa 2312 1841}%
\special{pa 2343 1839}%
\special{pa 2373 1838}%
\special{pa 2403 1839}%
\special{pa 2433 1842}%
\special{pa 2464 1847}%
\special{pa 2494 1853}%
\special{pa 2524 1860}%
\special{pa 2554 1869}%
\special{pa 2585 1878}%
\special{pa 2615 1889}%
\special{pa 2646 1900}%
\special{pa 2676 1912}%
\special{pa 2707 1924}%
\special{pa 2737 1937}%
\special{pa 2768 1950}%
\special{pa 2799 1963}%
\special{pa 2830 1976}%
\special{pa 2861 1989}%
\special{pa 2892 2001}%
\special{pa 2923 2013}%
\special{pa 2954 2025}%
\special{pa 2985 2036}%
\special{pa 3017 2046}%
\special{pa 3048 2055}%
\special{pa 3080 2063}%
\special{pa 3112 2070}%
\special{pa 3144 2076}%
\special{pa 3175 2081}%
\special{pa 3207 2085}%
\special{pa 3240 2088}%
\special{pa 3272 2090}%
\special{pa 3304 2091}%
\special{pa 3336 2092}%
\special{pa 3368 2091}%
\special{pa 3401 2090}%
\special{pa 3433 2088}%
\special{pa 3465 2085}%
\special{pa 3497 2082}%
\special{pa 3530 2078}%
\special{pa 3562 2073}%
\special{pa 3594 2068}%
\special{pa 3626 2062}%
\special{pa 3658 2056}%
\special{pa 3690 2049}%
\special{pa 3722 2042}%
\special{pa 3754 2035}%
\special{pa 3786 2027}%
\special{pa 3818 2018}%
\special{pa 3849 2009}%
\special{pa 3880 2000}%
\special{pa 3912 1991}%
\special{pa 3943 1982}%
\special{pa 3974 1972}%
\special{pa 4005 1962}%
\special{pa 4035 1952}%
\special{pa 4066 1942}%
\special{pa 4084 1936}%
\special{sp}%
% SPLINE 1 0 3 0
% 6 1895 2545 2239 2425 2967 2641 3975 2529 3999 2521 3999 2521
% 
\special{pn 13}%
\special{pa 1895 2145}%
\special{pa 1925 2131}%
\special{pa 1955 2116}%
\special{pa 1986 2103}%
\special{pa 2016 2089}%
\special{pa 2046 2077}%
\special{pa 2076 2065}%
\special{pa 2106 2054}%
\special{pa 2137 2045}%
\special{pa 2167 2037}%
\special{pa 2197 2031}%
\special{pa 2227 2026}%
\special{pa 2258 2024}%
\special{pa 2288 2023}%
\special{pa 2318 2024}%
\special{pa 2348 2027}%
\special{pa 2379 2032}%
\special{pa 2409 2038}%
\special{pa 2439 2045}%
\special{pa 2469 2054}%
\special{pa 2500 2063}%
\special{pa 2530 2074}%
\special{pa 2561 2085}%
\special{pa 2591 2097}%
\special{pa 2622 2109}%
\special{pa 2652 2122}%
\special{pa 2683 2135}%
\special{pa 2714 2148}%
\special{pa 2745 2161}%
\special{pa 2776 2174}%
\special{pa 2807 2186}%
\special{pa 2838 2198}%
\special{pa 2869 2210}%
\special{pa 2900 2221}%
\special{pa 2932 2231}%
\special{pa 2963 2240}%
\special{pa 2995 2248}%
\special{pa 3027 2255}%
\special{pa 3059 2261}%
\special{pa 3090 2266}%
\special{pa 3122 2270}%
\special{pa 3155 2273}%
\special{pa 3187 2275}%
\special{pa 3219 2276}%
\special{pa 3251 2277}%
\special{pa 3283 2276}%
\special{pa 3316 2275}%
\special{pa 3348 2273}%
\special{pa 3380 2270}%
\special{pa 3412 2267}%
\special{pa 3445 2263}%
\special{pa 3477 2258}%
\special{pa 3509 2253}%
\special{pa 3541 2247}%
\special{pa 3573 2241}%
\special{pa 3605 2234}%
\special{pa 3637 2227}%
\special{pa 3669 2220}%
\special{pa 3701 2212}%
\special{pa 3733 2203}%
\special{pa 3764 2194}%
\special{pa 3795 2185}%
\special{pa 3827 2176}%
\special{pa 3858 2167}%
\special{pa 3889 2157}%
\special{pa 3920 2147}%
\special{pa 3950 2137}%
\special{pa 3981 2127}%
\special{pa 3999 2121}%
\special{sp}%
% VECTOR 1 0 3 0
% 2 2385 2230 2745 2246
% 
\special{pn 13}%
\special{pa 2385 1830}%
\special{pa 2745 1846}%
\special{fp}%
\special{sh 1}%
\special{pa 2745 1846}%
\special{pa 2679 1823}%
\special{pa 2692 1844}%
\special{pa 2678 1863}%
\special{pa 2745 1846}%
\special{fp}%
% VECTOR 1 0 3 0
% 2 2345 2415 2657 2431
% 
\special{pn 13}%
\special{pa 2345 2015}%
\special{pa 2657 2031}%
\special{fp}%
\special{sh 1}%
\special{pa 2657 2031}%
\special{pa 2591 2008}%
\special{pa 2604 2028}%
\special{pa 2589 2048}%
\special{pa 2657 2031}%
\special{fp}%
% STR 2 0 3 0
% 3 5052 3244 5052 3324 2 0
% $\cB_R$
\put(50.5200,-29.2400){\makebox(0,0)[lb]{$\cB_R$}}%
% STR 2 0 3 0
% 3 5124 3500 5124 3580 2 0
% $t$
\put(51.2400,-31.8000){\makebox(0,0)[lb]{$t$}}%
% STR 2 0 3 0
% 3 2515 2695 2515 2775 2 0
% ${\cal C}$
\put(25.1500,-23.7500){\makebox(0,0)[lb]{${\cal C}$}}%
% SPLINE 1 0 3 0
% 6 1852 1996 2564 1804 3020 1940 3612 1804 4188 1980 4188 1980
% 
\special{pn 13}%
\special{pa 1852 1596}%
\special{pa 1883 1581}%
\special{pa 1914 1565}%
\special{pa 1945 1550}%
\special{pa 1976 1535}%
\special{pa 2007 1520}%
\special{pa 2038 1506}%
\special{pa 2069 1493}%
\special{pa 2100 1479}%
\special{pa 2131 1467}%
\special{pa 2162 1455}%
\special{pa 2193 1444}%
\special{pa 2224 1434}%
\special{pa 2255 1425}%
\special{pa 2286 1417}%
\special{pa 2317 1410}%
\special{pa 2347 1404}%
\special{pa 2378 1399}%
\special{pa 2409 1396}%
\special{pa 2440 1394}%
\special{pa 2471 1394}%
\special{pa 2501 1396}%
\special{pa 2532 1399}%
\special{pa 2563 1404}%
\special{pa 2593 1410}%
\special{pa 2624 1419}%
\special{pa 2654 1428}%
\special{pa 2685 1439}%
\special{pa 2716 1450}%
\special{pa 2746 1462}%
\special{pa 2777 1474}%
\special{pa 2807 1485}%
\special{pa 2838 1497}%
\special{pa 2869 1507}%
\special{pa 2899 1517}%
\special{pa 2930 1525}%
\special{pa 2961 1532}%
\special{pa 2992 1537}%
\special{pa 3023 1540}%
\special{pa 3054 1541}%
\special{pa 3085 1539}%
\special{pa 3116 1536}%
\special{pa 3147 1530}%
\special{pa 3178 1524}%
\special{pa 3209 1516}%
\special{pa 3241 1507}%
\special{pa 3272 1497}%
\special{pa 3303 1487}%
\special{pa 3334 1476}%
\special{pa 3366 1466}%
\special{pa 3397 1455}%
\special{pa 3428 1445}%
\special{pa 3460 1435}%
\special{pa 3491 1427}%
\special{pa 3522 1419}%
\special{pa 3553 1412}%
\special{pa 3584 1407}%
\special{pa 3615 1404}%
\special{pa 3646 1402}%
\special{pa 3677 1403}%
\special{pa 3708 1405}%
\special{pa 3739 1408}%
\special{pa 3769 1413}%
\special{pa 3800 1420}%
\special{pa 3831 1427}%
\special{pa 3861 1436}%
\special{pa 3892 1446}%
\special{pa 3923 1457}%
\special{pa 3953 1469}%
\special{pa 3984 1482}%
\special{pa 4014 1495}%
\special{pa 4044 1509}%
\special{pa 4075 1524}%
\special{pa 4105 1539}%
\special{pa 4136 1554}%
\special{pa 4166 1569}%
\special{pa 4188 1580}%
\special{sp}%
% VECTOR 1 0 3 0
% 4 2484 1780 2884 1796 2884 1796 2884 1796
% 
\special{pn 13}%
\special{pa 2484 1380}%
\special{pa 2884 1396}%
\special{fp}%
\special{sh 1}%
\special{pa 2884 1396}%
\special{pa 2818 1373}%
\special{pa 2831 1394}%
\special{pa 2817 1413}%
\special{pa 2884 1396}%
\special{fp}%
\special{pa 2884 1396}%
\special{pa 2884 1396}%
\special{fp}%
% STR 2 0 3 0
% 3 2765 2200 2765 2280 2 0
% $\tilde{v}_{|{\cal C}}$
\put(27.6500,-18.8000){\makebox(0,0)[lb]{$\tilde{v}_{|\cC}$}}%
% STR 2 0 3 0
% 3 2270 3930 2270 4030 2 0
% 
\put(22.7000,-36.3000){\makebox(0,0)[lb]{}}%
% VECTOR 1 0 3 0
% 2 5070 1670 5080 2900
% 
\special{pn 13}%
\special{pa 5070 1270}%
\special{pa 5080 2500}%
\special{fp}%
\special{sh 1}%
\special{pa 5080 2500}%
\special{pa 5099 2433}%
\special{pa 5080 2447}%
\special{pa 5059 2433}%
\special{pa 5080 2500}%
\special{fp}%
% LINE 1 0 3 0
% 2 3905 1355 4120 1355
% 
\special{pn 13}%
\special{pa 3905 955}%
\special{pa 4120 955}%
\special{fp}%
% LINE 1 0 3 0
% 2 1830 1355 1995 1355
% 
\special{pn 13}%
\special{pa 1830 955}%
\special{pa 1995 955}%
\special{fp}%
\end{picture}%
\end{center}
\caption{Nodal curves and Riccati equations for $\tilde{E_6}$ ($P_{IV}$)}
\label{fig:riccati} 
\end{figure}

The following proposition gives a characterization of Riccati 
equations obtained from the Painlev\'e equations in 
terms of rational nodal curves on   
Okamoto--Painlev\'e pair  $(S, Y)$. 
 (See Figure \ref{fig:riccati}).

\begin{Proposition}\label{prop:riccati}
Under the same notation as in Theorem \ref{thm:global}, let us consider the 
family $\pi:\cS \longrightarrow \cM_R \times \cB_R$ of 
the Okamoto--Painlev\'e pairs of type $R$  in (\ref{eq:globalf}).  
\begin{enumerate}
\item Assume that for a point 
 $t_0' = (\balpha_0, t_0) \in \cM_R \times \cB_R$, 
 there exists a smooth rational curve 
 $C \subset \cS_{(\balpha_0, t_0)} - \cD_{(\balpha_0, t_0)}$. 
  Then there exists an (analytic or \'etale)  open neighborhood  
$U $ of $t_0 $ of  $\cB_R$  satisfying the following conditions.   
\begin{enumerate}
\item There exist a flat family of rational curves 
$ \varphi: \cC \lra \{ \balpha_0 \} \times U $ 
and an inclusion $ \iota: \cC   \hookrightarrow  
\cS-\cD|_{\{\balpha_0 \} \times 
U}$ 
such that the following diagram is commutative:
\begin{equation}\label{eq:nodal-fam}
\begin{array}{ccccl}
C &\hookrightarrow& \cC &  \stackrel{\iota}{\hookrightarrow} &
 \cS-\cD|_{\{\balpha_0 \} \times U }      \\
\downarrow & & \varphi \downarrow \quad  &    \swarrow &  \hspace{-0.3cm} \pi  \\
(\balpha_0, t_0)  & \in &\{\balpha_0 \} \times 
U &  &
\end{array} 
\end{equation}

\item The restriction of the 
vector field 
$\tilde{v} \in 
\Gamma(\cS, \Theta_{\cS}( - \log \cD) \otimes \cO_{\cS}(\cD))$
 in (\ref{eq:vf}) to $\cC$ is tangent to $\cC$, that is, 
\begin{equation}\label{eq:riccati-vf1}
\tilde{v}|_{\cC} \in H^0(\cC, \Theta_\cC).
\end{equation}
Moreover  $\tilde{v}_{|\cC} $ defines a Riccati equation.   
\end{enumerate}

\item Conversely, assume that the restriction of 
Painlev\'e equation $\tilde{v}_{|\cS'}$ to 
the family $\pi': \cS':= \cS_{|\{ \balpha_0\} \times \cB_R} 
\lra \{ \balpha_0\} \times \cB_R $ can be reduced 
to a Riccati equation on 
an open neighborhood $ \{ \balpha_0 \} \times U$  of 
 a point $(\balpha_0, t_0)  \in \{ \balpha_0 \} \times \cB_R$.  
 Then there exist a family of rational nodal curves
$\cC \lra \{ \balpha_0 \} \times U$ on $\pi':\cS' - \cD' 
\lra \{ \balpha_0 \} \times U$.   
 
\end{enumerate}
\end{Proposition}

{\it Proof.}
Let us set $ \cB_{\balpha_0} = \{ \balpha_0 \} \times \cB_R
\hookrightarrow \cM_R  \times \cB_R $,   
$  t_0' = (\balpha_0, t_0)$. Restricting  the family
 $ \cS \lra \cM_R \times \cB_R $ to $ \cB_{\balpha_0}$, we obtain 
 a smooth projective family of surfaces:
 $$
\pi':\cS' :=  \cS_{|\cB_{\balpha_0}} \lra \cB_{\balpha_0}. 
 $$
 Moreover, we set $S_{t_0'} ={\pi'}^{-1}(t_0')$.  Fix a relatively ample 
 line bundle $H$ for  $\pi':\cS' \lra \cB_{\balpha_0}$. 
Consider the connected component $ T $ of the Hilbert scheme 
$\Hilb (\cS'/\cB_{\balpha_0})$ which contains a point $[C]$ and let 
$\cC \lra T$ denote the corresponding universal family.  
(Since $\pi'$ is projective and smooth, the universal family 
 $\tau:\cC \lra T$ exists (cf.\ [Theorem 1.4, Ch. I, \cite{Kol}]).) 
 
 Moreover we have a natural morphism 
 $\phi: T \rightarrow  \cB_{\balpha_0} $  
and a natural  inclusion 
 $\iota: \cC \hookrightarrow T \times_{\cB_{\balpha_0}} \cS'$, so that  
 $\tau:\cC \lra T$ can be factorized  into $\tau = p_1 \circ \iota$ where 
 $p_1$ denotes the first projection.

 Let $(Q, m_Q)$ be the local ring of $T$ at $[C]$.  Then from 
 [Theorem 2.10, Ch. I, \cite{Kol}], one can  see the following: 
  \begin{enumerate}
\item  The $\cO_{ \cB_{\balpha_0}, t_0'}$ -algebra $Q$ can be written as 
 the quotient of  a local $\cO_{ \cB_{\balpha_0},t_0'}$ -algebra $P$, 
 where
 $$
 \Spec P \lra \cB_{\balpha_0}
 $$ 
 is smooth of relative dimension $ d = \dim H^0(C, N_{C/S_{t_0'}})$.

\item  The kernel $K = \ker [P \rightarrow Q]$ is generated by 
 $\dim Obs(C)$ elements where $Obs(C)$ denotes the space of 
 obstructions.

\end{enumerate}

Since $C \subset S_{t_0'}$ is a $(-2)$-curve, we see that $N_{C/S_{t_0'}} 
\simeq \cO_{C}(-2)$, and hence  we have 
$H^0(C, N_{C/S_{t_0'}}) = H^0(\BP^1, \cO_{\BP^1}(-2)) = \{0\}$.  Therefore 
$\Spec P \lra \cB_{\balpha_0}$ is smooth of relative dimension $0$. 
Now we claim that:
\begin{equation}\label{eq:claim-ob}
 \fbox{{\bf Claim:}  $Obs(C) = \{ 0 \}$.}  
\end{equation}

Assuming the claim, we see that
$$
P \simeq Q \simeq  \cO_{ \cB_{\balpha_0}, t_0'}, 
$$
hence this implies that $T$ is a smooth variety of dimension $1$ at 
the point $[C]$ and 
the morphism $\phi:T \lra \cB_{\balpha_0}$ 
is  also an isomorphism near $[C]$ (\'etale or analytic) locally.  Hence 
we obtain an open neighborhood $U'$ of $[C]$ in $T$ on which the 
morphism $\phi$ induces the isomorphism 
$ \phi_{|U'}: U' \stackrel{\simeq}{\lra}\phi(U') \subset \cB_{\balpha_0}$.   
 It is clear that $U' = {\balpha_0} \times U$
for some open neighborhood of $t_0$ in $\cB_{R}$  and the restriction of the 
family $ \cC \lra T$ to $U'$  gives a family of rational curves 
$\cC \lra \{\balpha_0\} \times U$ 
which is a deformation of the rational curve $C$ in  $S_{t_0'}$.  

Now we show  the claim (\ref{eq:claim-ob}). 

From [Proposition 2.14, Ch.\ 1., \cite{Kol}], 
one see that the space of the obstructions $Obs(C)$ lies 
in $H^1(C, N_{C/S_{t_0'}}) $.  
Consider the natural homomorphisms of cohomology groups
$$
H^1(S_{t_0'}, \Theta_{S_{t_0'}}) \stackrel{\nu}{\lra} H^1( C,  \Theta_{S_{t_0'}|C} ) 
\stackrel{\mu}{\lra} H^1( C, N_{C/S_{t_0'}}).   
$$
Combining the Kodaira-Spencer homomorphism $\rho: T_{\cB_{\balpha_0}, t_0'} 
\lra H^1(S_{t_0'}, \Theta_{S_{t_0'}})$, it is easy to  see that 
\begin{equation}\label{eq:obs}
Obs(C) =   \mu \circ \nu \circ \rho (T_{\cB_{\balpha_0}, t_0'} ). 
\end{equation}

For simplicity, we set $S = S_{t_0'}, Y = \cY_{t_0'}, D = Y_{red} = \sum_{i=1}^r   Y_i  $.  

Since $ C \subset S - D$, we see that 
$$ 
\Theta_{S}(- \log (D +C))_{|D} 
\simeq \Theta_S ( -\log D)_{|D}.
$$
Therefore we have the following exact sequence 
$$
0 \lra \Theta_S ( - \log (D+C)) \lra 
\Theta_S (- \log (D+C))(D) \lra \Theta_S(-\log D) \otimes N_{D/S} \lra 0.  
$$
For an Okamoto--Painlev\'e pair $(S, Y)$ of non-fibered type, 
we have $H^0(S, \Theta_S (- \log (D+C))(D)) = \{ 0 \}$ 
(cf. [Proposition 2.1, \cite{STT}]). Hence, this  gives an injective homomorphism 
\begin{equation}\label{eq:normal}
0 \ra H^0(D, \Theta_S(-\log D) \otimes N_{D/S}) 
\ra H^1(S, \Theta_S ( - \log (D+C)) \ra H^1 (S,  \Theta_{S}( - \log(D+C))(D)).   
\end{equation}
We also have the following commutative 
 diagram of  sheaves (cf. [Lemma 2.1, \cite{STT}]):
$$
\begin{array}{ccccccc}
                &  0  & & 0&    &  0   &   \\   
            &  \downarrow    &  & \downarrow&    &  \downarrow   &   \\   
0 \lra & \Theta_S ( - \log (D+C)) &  \lra  &  \Theta_S (- \log D)&  
\lra & N_{C/S}&  \lra 0  \\
           &  \downarrow    &  & \downarrow&   &  \downarrow   &   \\   
0 \lra & \Theta_S ( - \log C ) &  \lra &\Theta_S&  \lra & N_{C/S}&  \lra 0 \\
           &  \downarrow    &  & \downarrow&   &  \downarrow   &   \\   
0 \lra & \oplus_{i=1}^r N_{Y_i/S} & \lra &\oplus_{i=1}^r N_{Y_i/S}& \lra& 0 &  \lra 0 \\ 
           &  \downarrow    &  & \downarrow&   &  \downarrow   &   \\   
                      &  0    &  &0&   &  0   &    .   
\end{array}
$$
Since $ N_{Y_i/S} = \cO_{Y_i}(-2)$ and $N_{C/S} = \cO_{C}(-2)$, we have  
the inclusions
$$
H^1( S, \Theta_S( - \log (D+C))) \hookrightarrow H^1( S, \Theta_S( - \log C))
\hookrightarrow H^1(S, \Theta_S).  
$$
Combining this and (\ref{eq:normal}), we see that 
\begin{equation}\label{eq:kernel}
 H^0(D, \Theta_S(-\log D) \otimes N_{D/S}) \hookrightarrow 
H^1( S, \Theta_S( - \log (D+C))) \subset \ker [ H^1(S, \Theta_S) 
\stackrel{\mu \circ  \nu}{\lra} H^1(C, N_{C/S})].   
\end{equation}
From (\ref{eq:com-diag1}), we have  
$$
\rho(T_{\cB_{\balpha_0}, t_0'} ) \simeq H^0(D, \Theta_S(-\log D) \otimes N_{D/S}), 
$$
and hence 
$$
\mu \circ \nu \circ \rho(T_{\cB_{\balpha_0}, t_0'} ) = \{ 0 \}.   
$$
Together with (\ref{eq:obs}), this shows the claim (\ref{eq:claim-ob}). 

Next, let us consider the family 
\begin{equation}\label{eq:dia}
\begin{array}{ccc}
  \cC &  \hookrightarrow &  \cS'_{|U}  \\
    &  \searrow & \quad \downarrow \pi \\
    &       &   U  .
    \end{array}
\end{equation} 
Since $\cD \cap \cC = \emptyset$, 
 we have  $\Theta_{\cS'|\cC} = 
 \Theta_{\cS'} ( -\log \cD) \otimes \cO_{\cS'}(\cD)_{|\cC}, $
and hence we obtain the following exact sequence: 
 $$
0 \lra \Theta_{\cC} \lra 
\Theta_{\cS'}( -\log \cD) \otimes \cO_{\cS'}(\cD)_{|\cC}
\lra N_{\cC/\cS'} \lra 0. 
$$
Since $N_{\cC/\cS'|\cS'_t} = N_{\cC_t/\cS'_t} = \cO_{\cC_t}(-2)$, 
we can show that $\pi_*(N_{\cC/\cS'}) = \{ 0 \}$.  Then we have 
$\Gamma(\cC, N_{\cC/\cS'}) = \{ 0 \}$.  This implies that 
$$
H^0(\cC, \Theta_{\cC}) \simeq 
H^0(\cS', \Theta_{\cS'}( -\log \cD) \otimes \cO_{\cS'}(\cD)_{|\cC}).
$$
Hence $\tilde{v}_{|\cC} \in H^0(\cC, \Theta_{\cC})$. 

Moreover, we may assume that $\cC \lra U$ is a trivial 
   $\BP^1$-bundle, that is,  $\cC \simeq \BP^1 \times U$ 
analytically.   
Since $\tilde{v}_{|\cC}$ defines a holomorphic vector field on 
$\BP^1 \times U$,  
it is easy to see that  $\tilde{v}_{|\cC}$ is equivalent 
to a Riccati equation (cf. Remark \ref{rm:riccati}).  

\nvskp
The second assertion is now obvious, because  the space of initial 
conditions of a Riccati equation must be a family of $\BP^1$.  
(Cf. Remark \ref{rm:riccati}).
\qed

\begin{Remark}[Global deformation of a $(-2)$-curve $C$]{\em 
Let us consider 
the connected component $T$ of the Hilbert scheme $\Hilb(\cS'/\cB_{\balpha_0})$ which contains a point $[C]$ and the corresponding universal family $
\tau:\cC \lra T$ in the proof of Proposition \ref{prop:riccati}.  
The argument in the proof shows that $\dim T = 1$ and 
the natural morphism 
$$
\phi: T \lra \cB_{\balpha_0}
$$
is projective, and hence surjective. We see that $\phi$ is a finite morphism of degree $d \geq 1$.  
Assume that $\phi$ is an isomorphism, i.e., $d =1$.   Then we have the global 
family of rational curves $\cC \subset \cS'$ over the affine curve 
$\cB_{\balpha_0}$:
\begin{equation}
\begin{array}{ccccc}
C & \subset & \cC &  \hookrightarrow  &  \cS' \\
\downarrow&  &\downarrow   &  & \downarrow  \\
 t_\0' & \in   &  T &  = &  \cB_{\balpha_0}
   \end{array}.
 \end{equation}
 Then the vector field $\tilde{v}_{|\cC}$ becomes an algebraic regular vector 
 field on $\cC$ and 
 defines a Riccati equation over the affine algebraic curve  
 $\cB_{\balpha_0}$.  
In this case, we call the differential equation defined by   
$\tilde{v}_{|\cC}$ the {\em Riccati equation} associated to the rational curve 
 $C (\subset S - D)$.  
 
  We do not know whether the case with $ d > 1$ occurs  or does
   not occur. However, 
if $\phi:T \lra \cB_{\balpha_0}$ is of degree $d > 1$, we see that 
$\phi^{-1}( \phi([C]))$ consists of $d$ rational curves of $\cS_{t_0'}$ 
$C_1 := C, C_2, \cdots, C_d$ which are in the flat family of 
rational curves in $\cS'$  parameterized by a connected variety $T$.  

In \S 2, we see that there exists an 
Okamoto--Painlev\'e pair $(S, Y)$ 
which contains more than one  rational curves $C_i \subset S - D$,  $i \geq 2$.

}
\end{Remark}

\begin{Definition}\label{def:riccati}
{\rm   Under the same notation and assumptions in Proposition \ref{prop:riccati},  
we call the differential equations determined 
by the vector field $\tilde v|_\cC$ in (\ref{eq:riccati-vf})
 {\em Riccati equation associated with the rational 
curve  $C \subset S - Y_{red}$}. Moreover  we call a solution of the 
Riccati equation $\tilde v|_\cC$  {\em a  Riccati solution of the 
 Painlev\'e system (associated with $C \subset S- Y_{red}$}).
(Note that all solutions of $\tilde v|_\cC$ remain in the family of rational 
curves in (\ref{eq:nodal-fam})).  }
\end{Definition}

\section{\bf Classification of $(-2)$-rational 
curves (nodal rational curves) on $S-D$}\label{section:class}

Let $(S,Y)$ be a rational Okamoto--Painlev\'e pair  of non-fibered 
type which corresponds to a Painlev\'e equation (cf.  Table \ref{tab:type}).
 
In this section, we will  classify all configurations of  
$(-2)$-curves on $S-D$ for a rational Okamoto--Painlev\'e pair $(S, Y)$ of 
non-fibered type.   
The classification of the configurations are based on 
the similar classification for rational Okamoto--Painlev\'e pairs  $(S, Y)$ of 
fibered type with the elliptic fibration $f:S \lra \BP^1$  and some deformation arguments.

\subsection{Notations and the Result}

Let $S$ be a projective smooth surface over $\C$. 
We denote by $Div(S)$ the free abelian group  generated by all irreducible 
curves on $S$.  Let $ \sim_a$ and $\sim$  denote the algebraic equivalence  
and the linear equivalence of divisors respectively.  We define 
the N\'eron--Severi group and the Picard group of $S$ by 
\begin{eqnarray}
    \NS(S) & = & Div(S)/ \sim_a,   \\
    \Pic(S) & = & Div(S)/\sim.   
\end{eqnarray} 

In what follows, we assume that $S$ is a rational surface.  Then
we have the natural isomorphisms
\begin{equation} \label{eq:picard}
\Pic(S) \simeq \NS(S) \simeq H^2(S, \Z), 
\end{equation}
and these groups are free $\Z$-modules of rank $b_2(S)$.   
For any divisor $C$, we also denote by the same letter 
$C$ the class of the divisor 
in $\NS(S) \simeq H^2(S, \Z)$.  Moreover $C = D$ means that 
the two divisors are linear equivalent to each other.  
We can consider the lattice structure on these free $\Z$-modules 
by the intersection form $< \  , \  >$ on $\NS(S)$ or equivalently by 
the cup product on $H^2(S, \Z)$.  
Let $E_8$ be the unique even unimodular positive-definite lattice of rank $8$. 
For a lattice $L = (L, < \  , \  >)$, we denote by 
$L^- = (L,(-1) \times  < \  , \  >)$, the opposite lattice of $L$.  
Note that the opposite lattice $E_8^-$ of $E_8$ is negative-definite.

Let $(S, Y)$ be a rational Okamoto--Painlev\'e pair and let 
\begin{equation}\label{eq:decomp}
Y = \sum_{i=1}^r m_i Y_i 
\end{equation}
be the irreducible decomposition of $Y$.   
Since $S$ is a rational surface with $b_2(S) = \rank H^2(S, \Z) = 10$, 
by the Hodge index theorem, the   bilinear form 
$ <\ , \ >$ on $H^2(S, \Z)$ can be written as the diagonal 
matrix $ (1, \underbrace{-1, \cdots, -1}_{9}) $.  
 The sub-lattice $M(Y)$ generated by $ \{ Y_i \}_{i=1}^r$ in $H^2(S, \Z)$  
is a root lattice of an affine type, say $R = R(Y)$.  
Since $S$ is not relatively minimal, $S$ 
contains a $(-1)$-rational curve $O$ on $S$.  
Then by the adjunction formula, one has $Y \cdot O = -K_S \cdot O = 1$.  
Hence, there exists a $i_0, 1 \leq i_0 \leq r $ such that  
$m_{i_0} =1$ and  $Y_{i_0} \cdot O =1$.  By renumbering  $i$, 
we may assume that $i_0 = 1$. 
Define  the sub-lattice by 
\begin{equation}\label{eq:classical}
M'(Y) =  \langle  Y_2, \cdots, Y_r \rangle_{\Z} \subset M(Y), 
\end{equation}
which is a root lattice of classical type $R'$. 
 For example, if $R = \tilde{D}_4$, then $R' = D_4$.   
Let  $M(S-Y_{red})$ be the sub-lattice $H^2(S, \Z)$ 
generated by all $(-2)$-curves $C$ on $S - Y$.  
Note that we have the orthogonal sum 
\begin{equation}\label{eq:orth}
 M'(Y) \oplus M(S - Y_{red}) \subset H^2(S, \Z).  
\end{equation}

\begin{Lemma} \label{lem:root}
 Assume that $(S, Y)$ is of non-fibered type. 
 Then $M'(Y) \oplus M(S-Y_{red}) $ is a root sub-lattice of $E_8^-$. 

\end{Lemma}

\noindent
{\it Proof.}  \quad  
  The sub-lattice $\langle Y,O \rangle_\Z$
generated by $Y$ and $O$ has the intersection matrix  $ \left( \begin{array}{cc} 0 & 1 \\ 1 & -1  
\end{array}  \right) $.  

 Then  the orthogonal complement $\langle Y,O \rangle^\perp $ 
  in  $H^2(S,\Z)$ is an even, negative-definite unimodular 
  lattice  of rank $8$, which is 
isomorphic to the root lattice $E_8^{-}$.  (Since $K_S = -Y$, the 
adjunction formula implies that $\langle Y,O \rangle^\perp $ 
is even). 
Since $Y \cdot O = 1$, we see that 
the orthogonal complement $\langle Y \rangle^{\perp}$ is given by 
$$
\langle Y \rangle^{\perp} \simeq  
\langle Y, O \rangle^{\perp} \oplus \Z Y \simeq \tilde{E}_8^{-}. 
$$

Since $M(S-Y_{red})$ is generated by $(-2)$-curves on $S - Y_{red}$, 
we see that $M(S-Y_{red}) \subset  \langle Y \rangle^{\perp}$.  
Moreover by definition of Okamoto--Painlev\'e pair 
(cf.\ (\ref{eq:canonical})),  $M'(Y) \subset  \langle Y \rangle^{\perp}$.  (In fact, 
we have $ M'(Y) \subset \langle Y, O \rangle^{\perp}$).  
Set 
\begin{equation}
N(Y) := M'(Y) \oplus M(S - Y_{red}).
\end{equation}
Then $N(Y) \subset \langle Y \rangle^{\perp}$.  

Let us consider the natural projection map
$$
\pi:\langle Y \rangle^{\perp} \simeq \langle Y, O \rangle^{\perp} \oplus \Z 
Y \lra  \langle Y, O \rangle^{\perp}.  
$$
We claim that: 
\begin{equation}\label{eq:claim-2}
\mbox{\bf Claim}: \quad \pi_{|N(Y)} \quad \mbox{is injective}.
\end{equation}
 If the claim is true, 
we see that $N(Y) \simeq \pi(N(Y)) \subset \langle Y, O 
\rangle^{\perp} \simeq E_8^{-}$. This implies that $N(Y)$ 
is a negative-definite lattice generated by $(-2)$-elements.  
Hence one can see that 
  $N(Y)$ is a root lattice which  is a direct sum of root
   lattices of type $A_i$, $D_j$, $E_k$.  (This also implies that $M'(Y)$ and $M(S-Y_{red})$ are 
  direct sums of root lattices of type $A_i$, $D_j$, $E_k$.) 
To show the claim (\ref{eq:claim-2}),  it suffices to show that 
$ \Ker \pi_{|N(Y)} = \Ker \pi \cap N(Y) = \{ 0 \}$.
Since $\Ker \pi = \Z [Y]$ with $Y^2 = 0$ and 
 $M'(Y)$ is negative-definite, we have  
$$
\Ker \pi \cap N(Y)  = \Z [Y] \cap N(Y) = \Z [Y] \cap 
(M'(Y) \oplus M(S - Y_{red})) = \Z[Y] \cap  M(S - Y_{red}).   
$$  
Hence we have to show that $\Z[Y] \cap  M(S - Y_{red}) =\{ 0 \}$.  
Take $\gamma \in \Ker \pi_{|M(S-Y_{red})}$ and assume that 
$\gamma \not= 0$.   Since $\Ker \pi = \Z \cdot Y$, we can write 
$ \gamma $  as $\gamma = b \cdot Y$ with $b \not= 0$. We may assume that 
$ b > 0$.  On the other hand, since $\gamma \in M(S - Y_{red})$, 
we can write $\gamma$ as 
$$
\gamma = C - D 
$$
with  
$$
C = \sum_{i=1}^l a_i C_i, \quad  D = \sum_{j=1}^t b_j D_j 
$$
where $C_i$ ($1 \leq i \leq l$) and $D_j$  ($ 1 \leq j \leq t$) 
are different $(-2)$-curves in $S - Y_{red}$ and $a_i \geq 0, b_j \geq  0$.  
Assume that $D = 0$. Then we see that $b Y$ and  $ C$ are linear equivalent 
to each other. Since $bY$ and $C$ are different effective divisors, 
we see that $\dim H^0(S, \cO_S(b Y)) \geq 2$.  
This contradicts to the fact that $(S, Y)$ is of 
non-fibered type (cf. Proposition 1.3, \cite{STT}). 
Therefore we may assume that both of 
$C$ and $D$ are non-zero effective divisors.  
Recall that the lattice $\langle Y \rangle^{\perp}$  
is negative semi-definite. Hence one has 
$$
0 \geq C^2 = (D + bY)^2 = D^2 = D \cdot C \geq 0.
$$
(Here we used the fact that $D \cdot Y = C \cdot Y = 0$).  
This implies that 
$$
C^2 = D^2 = C \cdot D = 0.
$$
An element $G  \in \langle Y \rangle^{\perp}$ with $G^2 = 0$ must be 
proportional to $Y$, that is, $ G = c Y$.  Therefore we see that 
$ C = b'Y $ with $ b' > 0$,  which again contradicts to the fact 
that $(S, Y)$ is of non-fibered type.   We have proved that 
$\Ker \pi_{|M(S- Y_{red})} = \{ 0 \}$ and hence $\Ker \pi_{|N(Y)} = \{ 0 \}$
as in (\ref{eq:claim-2}).  
\qed

\vspace{0.5cm}
By Lemma \ref{lem:root}, there are only finitely many $(-2)$ curves 
$\{ C_i \}_{i=1}^l$ on $S - Y_{red}$.  The dual graph of configurations of 
$(-2)$-curves on $S$ can be classified by the Dynkin diagram of ADE types. 
The following theorem is the main theorem in this section.

\begin{Theorem} \label{thm:main-2}
Let $(S,Y)$ be a  rational Okamoto--Painlev\'e pair of 
non-fibered type which corresponds to a Painlev\'e equation 
(cf. Table \ref{tab:type}).  
The type of the root lattice $M(S - Y_{red})$,  or equivalently, 
the dual graph of the configuration of $(-2)$-curves 
on $S-Y$ are classified in  Table \ref{tab:config-2nonfibered}.
\small
\begin{table}[h]
\renewcommand{\arraystretch}{1.2}
\begin{center}
\begin{tabular}{|c|c @{\ \vrule width0.8pt \quad} l|}
\hline
Painlev\'e & $ R(Y)$  & the type of the dual graph of configuration of   \\
equations &  &  $
(-2)$-curves on $S-Y$ \\
\noalign{\hrule height 0.8pt}
$P_{VI}$&$\tilde D_4$ & $ D_4 ,\quad ( A_1, A_1, A_1, A_1),\quad 
 A_3 ,\quad ( A_1, A_1, A_1) ,\quad  A_2 ,\quad ( A_1, A_1) ,\quad  A_1$ \\
\hline
$P_{V}$ & $\tilde D_5$ & $ A_3 ,\quad  A_2 ,\quad( A_1, A_1) ,\quad  A_1$   \\
\hline
$P_{III}^{\tilde{D_6}}$ & $\tilde D_6$ & $( A_1, A_1) ,\quad  A_1 $ \\
\hline
$P_{III}^{\tilde{D_7}}$ & $\tilde D_7$ & none  \\
\hline
$P_{III}^{\tilde{D_6}}$ & $\tilde D_8$ & none \\
\hline
$P_{IV}$ & $\tilde E_6$ & $ A_2 ,\quad  A_1$   \\
\hline
$P_{II}$ & $\tilde E_7$ & $ A_1$  \\
\hline
$ P_I $ &$ \tilde E_8$ & none  \\
\hline
\end{tabular}
\caption{Configuration of $(-2)$-curves  on $S-Y$ for
 a rational Okamoto--Painlev\'e pair $(S, Y)$ of  non-fibered type}
\label{tab:config-2nonfibered}
\end{center}
\end{table}
\end{Theorem}

%\newpage

\subsection{The case of fibered type}\label{subsec:-2fibered}

Oguiso and Shioda \cite{O-S} give  the complete structure 
theorem of the Mordell-Weil group of rational 
elliptic surfaces $f:S \lra \BP^1$ with a section.  
Let $(S,Y)$ be a rational Okamoto--Painlev\'e pair of fibered type,  i.e. 
there exists an elliptic fibration $f:S \to \BP^1$ such that $f^*(\infty)=Y$. 
Since $K_S = f^{*}( - \infty ) = -Y $, by the adjunction formula, 
it is easy to check that an irreducible curve 
 $C$ is a $(-2)$-curve if and only if it is one of the irreducible 
 components of the reducible singular fibers. 
Hence, to give the complete structure of $(-2)$-curves 
on $S-Y$, we quote the structure of the reducible singular 
fibers which is a part of the structure 
theorem of the Mordell-Weil group of $f:S \lra \BP^1$. 

\nvskp
We will introduce some notations. 
Let  $(S, Y)$ be  a rational Okamoto--Painlev\'e pairs of 
fibered type with an elliptic fibration $f:S \lra \BP^1$ 
such that $f^{*}(\infty) = Y$. (Here,  we do not 
assume that the type of $Y$ is in Table \ref{tab:type}).  
We also assume that there exists
 a section $O \subset S$ and we denote by $F$ the class 
 of a general fiber of $f$ so that $Y$ and $F$ are linearly 
equivalent to each other, or equivalently,  
have the same class in $H^2(S, \Z)$.  For a lattice  $L$, 
let us denote by $L^-$ the opposite lattice of $L$, i.e.,
$$
L^-= \mbox{the module }L\mbox{ with the pairing }(-1) \times < \ ,\ >.
$$

Let $F_v := f^{-1}(v)$ denote the fiber over the closed point $v \in \BP^1$, and set
$$
\begin{array}{c}
Sing(f):= \{ v \in \BP^1 | F_v = f^{-1}(v) \mbox{ is singular } \}, \\
{\cal R} =Red(f):= \{ v \in \BP^1 | F_v = f^{-1}(v) \mbox{ is reducible } \}.
\end{array}
$$
For each $v \in {\cal R}$, let
$$
F_v = f^{-1}(v)= \Theta_{v,0}+\sum_{i=1}^{m_v-1} \mu_{v,i} \Theta_{v,i} 
\quad (\mu_{v,i} \geq 1,\mu_{v,0} = 1)$$
be the irreducible decomposition of $F_v$ where $\Theta_{v,0}$ 
is the unique component of $F_v$ 
meeting the zero section $O$ and $m_v$ is 
the number of irreducible components. We set 
\begin{equation}
T_v:= \langle \Theta_{v,i} |1 \le i \le m_v-1 \rangle_\Z \subset \NS(S),
\end{equation}
and
\begin{equation}
T:=\bigoplus_{v\in {\cal R}} T_v.
\end{equation}
Note that the notation $T$ is used for another lattice in \cite{Shi}.

By the classification of singular fibers (cf. \cite{Kod}), 
(and using the intersection matrix $(\Theta_{v,i} 
\cdot \Theta_{v,j})_{i \le i,j \le m_v-1}$), we 
have the following
\begin{Lemma}[Lemma 7.2 \cite{Shi}]\label{lem:typeF_v}
The opposite lattice $T_v^-$ is a root lattice of rank $m_v-1$, 
determined by the type of the singular fiber $F_v$ as follows:

\begin{center}
\begin{tabular}{|c@{\quad \vrule width0.8pt \quad}c|c|c|c|c|c|c|} \hline
Type of $F_v$ & $I_m$ & $I_m^*$ & $II^*$ & $III^*$ & $IV^*$ & $IV$ & $III$ \\
\hline 
$T_v^-$ & $A_{m-1}$ & $D_{m+4}$ & $E_8$ & $E_7$ & $E_6$ & $A_2$ & $A_1$ \\
\hline 
\end{tabular}
\end{center}

\end{Lemma}

Furthermore, we have (cf. (7.2) \cite{Shi})
$$
\langle O, F, \Theta_{v,i} \ (0 \le i \le m_v-1,v \in {\cal R}) 
\rangle_\Z =\langle O,F \rangle_\Z \oplus T \subset \NS(S) 
\quad (\mbox{orthogonal direct sum})
$$
where $F$ is the class of  a fiber of $f$. As we see in the previous 
subsection, we see that $\langle O, F \rangle^{\perp} \simeq E_8^{-}$. 

Hence we have an embedding
\begin{equation}
 T^-=\bigoplus_{v\in {\cal R}} T_v^- \hookrightarrow E_8. 
\end{equation}

Now we recall Dynkin's results on the classification 
of root lattices contained in $E_8$, which is equivalent 
to the classification of regular semisimple subalgebras 
of the exceptional Lie algebra of type $E_8$. 
\begin{Theorem}[Ch. II, Table 11 \cite{D}]\label{thm:dynkin}
Let $L$ be a root lattice of $\rank \  s$ which is embedded as 
a sub-lattice of $E_8$, other than $\{0\}$ and $E_8$. 
Then $L$ is isomorphic to one in Table \ref{tab:subE8}.
\small
\begin{table}[h] 
\renewcommand{\arraystretch}{1.2}
\begin{center}
\begin{tabular}{|c @{\ \vrule width0.8pt \quad} l|}
\hline
$\ s\ $ & $L$ \\
\noalign{\hrule height 0.8pt}
$8$ & $A_8,\ D_8,\ A_7\oplus A_1,\ A_5\oplus A_2\oplus A_1,
\ A_4^{\oplus2},\ A_2^{\oplus4},\ E_6\oplus A_2,\ E_7 \oplus A_1,
\ D_6\oplus A_1^{\oplus2},\ D_5\oplus A_3,$\\
& $D_4^{\oplus2},\ D_4\oplus A_1^{\oplus4},
\ A_3^{\oplus2}\oplus A_1^{\oplus2},\ A_1^{\oplus8}$\\
\hline
$7$ &$ A_6\oplus A_1,\ A_4\oplus A_2\oplus A_1,
\ A_5\oplus A_2,\ A_2^{\oplus3}\oplus A_1,\ E_6\oplus A_1,
\ E_7,\ D_7,\ D_5\oplus A_1^{\oplus2},$\\
& $D_4\oplus A_1^{\oplus3},\ A_3^{\oplus2}\oplus A_1,
\ A_1^{\oplus7},\ D_6\oplus A_1,\ D_5\oplus A_2,\ 
A_3\oplus A_2\oplus A_1^{\oplus2},\ D_4\oplus A_3,\ A_3\oplus A_1^{\oplus4},$\\
& $A_4\oplus A_3,\ A_5\oplus A_1^{\oplus2},\ A_7$\\
\hline
$6$ & $A_2^{\oplus3},\ E_6,\ D_6,\ D_4\oplus A_1^{\oplus2},
\ A_3^{\oplus2},\ D_5\oplus A_1,\ A_3\oplus A_1^{\oplus3},
\ D_4\oplus A_2,\ A_1^{\oplus6},\ A_2\oplus A_1^{\oplus4},$\\
& $A_4\oplus A_1^{\oplus2},\ A_6,\ A_3\oplus A_2\oplus A_1,
\ A_5\oplus A_1,\ A_4\oplus A_2,\ A_2^{\oplus2}\oplus A_1^{\oplus2}$\\
\hline
$5$ & $D_5,\ A_3\oplus A_1^{\oplus2},\ A_3\oplus A_2,\ A_5\ ,
A_1^{\oplus5},\ A_4\oplus A_1,\ D_4\oplus A_1,\ A_2\oplus A_1^{\oplus3},
\ A_2^{\oplus2}\oplus A_1$ \\
\hline
$4$ & $D_4,\ ,A_1^{\oplus4},\ A_2\oplus A_1^{\oplus2},\ 
A_2^{\oplus2},\ A_3\oplus A_1,\ A_4$\\
\hline
$3$ & $A_3,\ A_2\oplus A_1,\ A_1^{\oplus3}$\\
\hline
$2$ & $A_2,\ A_1^{\oplus2}$\\
\hline
$1$ & $A_1$\\
\hline
\end{tabular}
\caption{Root sub-lattice of $E_8$}
\label{tab:subE8}
\end{center}
\end{table}

\end{Theorem}
\normalsize
From Theorem \ref{thm:dynkin}, one can classify 
the root sub-lattice of $E_8$, hence $T$ must be one of the root lattices 
in the  Table \ref{tab:subE8}. 

However, as for the existence, we quote the following 

\begin{Theorem}[cf. Remark 2.7 \cite{O-S}]
For every type given in Table \ref{tab:subE8} except for the type 
$$
D_4\oplus A_1^{\oplus4},\ A_1^{\oplus8} \mbox{ and } A_1^{\oplus7},
$$
there exists a rational elliptic surface whose $T^-$ is of given type.
\end{Theorem}

\begin{Remark}[cf. Remark 3.4. \cite{O-S}]{\rm 
The sum of the local Euler number of the reducible singular 
fibers cannot exceed $12$, the Euler number of a rational elliptic surface.
Therefore, the types $D_4\oplus A_1^{\oplus4},
\ A_1^{\oplus8} \mbox{ and } A_1^{\oplus7}$ do not appear.

}
\end{Remark}

In the case of a  rational Okamoto--Painlev\'e pair $(S,Y)$ of fibered type 
in Table \ref{tab:type}, 
the type of root lattice $T_\infty$ is determined by the type of $Y$. 
Thus, we obtain the classification theorem as follows.

\begin{Proposition}\label{prop:fibered}
Let $(S,Y)$ be a  rational Okamoto--Painlev\'e pair 
of fibered type in Table 
\ref{tab:type}. The type of root lattice 
$\bigoplus_{v\in {\cal R}-\infty} T_v^-$ 
are classified by Table \ref{tab:subE8fibered}.

\begin{table}[h]
\renewcommand{\arraystretch}{1.2}
\begin{center}
\begin{tabular}{|l @{\ \vrule width0.8pt \quad} l|}
\hline
Type of $Y$  & $\bigoplus_{v\in {\cal R}-\infty} T_v^-$ \\
\noalign{\hrule height 0.8pt}
$\tilde D_4=I^*_0$ & $D_4,\quad A_3 ,\quad A_1^{\oplus3} ,
\quad  A_2 ,\quad A_1^{\oplus2} ,\quad A_1$ \\
\hline
$\tilde D_5=I^*_1$ & $A_3 ,\quad A_2 ,\quad A_1^{\oplus2} ,\quad A_1$   \\
\hline
$\tilde D_6=I^*_2$ & $A_1^{\oplus2} ,\quad A_1 $ \\
\hline
$\tilde D_7=I^*_3$ & none  \\
\hline
$\tilde D_8=I^*_4$ & none \\
\hline
$\tilde E_6=IV^*$ & $ A_2 ,\quad  A_1$   \\
\hline
$\tilde E_7=III^*$ & $ A_1$  \\
\hline
$\tilde E_8=II^*$ & none  \\
\hline
\end{tabular}
\caption{The lists of root lattice 
$\bigoplus_{v\in {\cal R}-\infty} T_v^-$ (fibered type)}
\label{tab:subE8fibered}
\end{center}
\end{table}

\end{Proposition}

\begin{Remark}{\rm
By Lemma \ref{lem:typeF_v} and Proposition 
\ref{prop:fibered}, the structure of configuration of 
$(-2)$-curves (i.e. type of $F_v$'s) is `almost' 
determined. For $T_v=A_1$, the type of $F_v$ cannot be distinguished between 
$I_2$ and $III$. Similarly, for $T_v=A_2$, 
the type of $F_v$ cannot be distinguished between 
$I_2$ and $IV$. 
(For other types, we can determine the type of $F_v$.)   
}
\end{Remark}

\subsection{Proof of Theorem 
\ref{thm:main-2}}\label{sssec:-2nonfib}
\ 
Now we  prove Theorem \ref{thm:main-2}.  
Let $(S, Y)$ be a rational Okamoto--Painlev\'e pair of non-fibered type 
with a given type $R=R(Y)$ of $Y$ in the Table \ref{tab:type}.  
Let $M'(Y)$ and $M(S - Y_{red})$ be the sub-lattices defined in (\ref{eq:orth}).
  By Lemma \ref{lem:root},  
 the orthogonal sum  
$M'(Y)^{-} \oplus M(S- Y_{red})^{-} $ is a root sub-lattice of $E_8$.  
Then since the type $R'(Y)$ of $M'(Y)^{-}$ is $D_{k}, 4 \leq k \leq 8$ or 
$E_6, E_7, E_8$, by the Classification Theorem \ref{thm:dynkin}, we can 
obtain the list of possible types for  $M(S- Y_{red})^{-} $ as in  Table 
\ref{tab:config-2nonfibered}.   

Therefore, it suffices to show for each type $R''$ of root lattices listed in 
Table \ref{tab:config-2nonfibered},  there exsits  a rational 
Okamoto--Painlev\'e pair $(S, Y)$ of non-fibered types with 
the root sub-lattice $M(S- Y_{red})$ of type $R''$.

First, let $(S, Y)$ be a rational Okamoto--Painlev\'e 
pair $(S, Y)$ of fibered type
with a given type of $Y$ in the Table \ref{tab:type} and let 
$f:S \lra \BP^1$ be the elliptic fibration with $f^{*}(\infty) = Y$. 

From Proposition \ref{prop:fibered}, we can  determine the possible 
configuration of $(-2)$-curves on $S - Y_{red} = S - f^{-1}(\infty)$ 

by the classification of the other reducible singular fibers.    
(Note that Proposition \ref{prop:fibered} says the existence of 
such a fibration.)
Let $ Y = \sum_{i=1}^r m_i Y_i$ be the irreducible decomposition of $Y$. 
Set $D  = Y_{red} =  \sum_{i=1}^r Y_i$, and take all $(-2)$ curves $
\{ C_1, \cdots, C_l \}$ on $S - Y_{red}$.   Note that 
each $C_i$ is an irreducible component of  reducible singular fibers of $f$.

Now we will use the following deformation 
argument. 

\begin{Lemma}\label{lem:defm nonfibered}
Let $(S,Y)$ be a rational Okamoto--Painlev\'e pair of fibered type 
with the irreducible decomposition $Y = \sum_{i=1}^{r} m_i Y_i$  such that 
 $D=Y_{red}$ is a normal crossing divisor, 
 and let $C=\sum_{j=1}^s C_j$ be a normal crossing divisor of $S$ 
satisfying the following conditions:
\begin{enumerate}
\item\label{as:c1} $C \subset S-D$,
\item $C_j \simeq \BP^1$,
\item The classes of curves $\{Y_i , C_j \ |\ 
 1 \le i \le r, 1 \le j \le s  \}$ are  
 linearly independent in  $H^2(S,\C) \simeq \Pic (S)\otimes_{\Z}\C $.
\end{enumerate}

Then there exists  a rational Okamoto--Painlev\'e pair $(S', Y')$ such that 
\begin{enumerate}
\item $(S', Y')$ is of non-fibered type,  
\item the type of $Y'$ is same as the type of $Y$, 
\item $S'- Y'_{red}$ contains $(-2)$ curves $\{ C'_j \}_{j=1}^s$ with the 
same configurations as $\{ C_j \}_{j=1}^s$, and  
\item $S'$ is a deformation of $S$. 
\end{enumerate}

\end{Lemma}

{\it Proof.}
Let $F$ be an arbitrary fiber at $\BP^1 - \{ \infty \} - Sing(f)$, 
which is an elliptic curve and $F \subset S-(D+C)$.
Let us consider the exact sequence of sheaves
$$
0 \to \Theta_S(-\log (D+C+F)) \to \Theta_S(-\log (D+C)) \to N_F \to 0, 
$$
which yields the exact sequence 
\small
\begin{equation}\label{eq:coh}
H^1(\Theta_S(-\log (D+C+F))) \to H^1(\Theta_S(-\log (D+C))) 
\stackrel{\phi}{\to} H^1(N_F) \to H^2(\Theta_S(-\log (D+C+F))).
\end{equation}
\normalsize
Since $F$ and $Y$ are linearly equivalent, 
we get $N_F = [F]_{|F} =[Y]_{| F} = \cO_F$, and hence  
$H^1(N_F)=H^1( \cO_F)=\C$.  In Lemma \ref{lem:h2tlDCE}, 
we will show 
\begin{equation}\label{eq:vanish}
H^2(\Theta_S(-\log (D+C+F))= \{ 0 \}. 
\end{equation}
From (\ref{eq:coh}) together with (\ref{eq:vanish}), 
we see that there exists an element $\theta 
\in H^1(\Theta_S(-\log (D+C))$ such that $\phi(\theta) \neq 0$. 
 Such an element 
$\theta$ induces an infinitesimal deformation of the pair $(S, D+C)$ which 
does not preserve the elliptic curve $F$.  
Since we see that $H^2(S, \Theta_S( -\log(D+C))) = \{0 \}$, 
such an  infinitesimal deformation $\theta$ induces  
 a one parameter deformation 
 $$
 \begin{array}{ccc}
 \cS & \hookleftarrow & {\cD + {\cal C}} \\
 \varphi \downarrow \quad & \swarrow &  \\
 \Delta &   &   \\
 \end{array}
$$ 
of $(S, D +C)$
where $\Delta = \{ z \in \C \ | \ |z| < \epsilon \}$ is a small 
neighborhood of the origin.  
Note that we also have the  relative divisor $\cY_i$ for $\varphi$ 
which gives the deformation of $Y_i$.  
Hence we  have the relative divisor 
$ \cY = \sum_{i=1}^{r} m_i \cY_i$.  
For $z \in \Delta$, denote by $\cS_z, \cY_{i, z}, \cD_z, \cC_z, 
\cY_z$ the corresponding fibers of $\cS, \cY_i, 
\cD, \cC$ and $\cY$ over $z$ respectively.   
It is obvious that for every $z \in \Delta$ 
each $\cY_{i, z}$ is a $(-2)$-curve on $\cS_z$ and $\cY_z$  
satisfies the numerical condition 
(\ref{eq:canonical}) that  $ \cY_z \cdot \cY_{i, z} = 0$ for all $i$.

%%%%
Consider the divisor $ K_{\cS} + \cY $   
on $\cS$ and set $\cL =\cO_{\cS}(K_{\cS} + \cY)$.  
We know the following two facts:
\begin{enumerate}
\item $\cL_{|\cS_0} \sim \cO_{\cS_0}$.
\item  Since $\cS_z$ is a projective smooth rational surface for every 
$z \in \Delta$, we see that $H^i(\cS_z, \cO_{\cS_z}) = 0$ 
for $ i \geq 1$ and  every $z \in \Delta$. In particular, 
$R^i \pi_* \cO_{\cS} = 0$ for $ i \geq 1$.  
\end{enumerate} 
Then by the upper-semicontinuity theorem, we see that 
$\dim H^i(\cS_{z}, \cL_{|\cS_z}) = 0$ for every $i \geq 1$.  
Noting  that $\pi_{*} \cL \simeq \cO_{\Delta}$, we see that 
there is a non-trivial homomorphism  
$ s:  \pi^{*}(\pi_{*}\cL ) = \cO_{\cS} \ra \cL$. 
Applying the same argument for the dual sheaf $\cL^{\vee}$, we also have 
a non-trivial homomorphism $s': \cO_{\cS} \ra \cL^{\vee}$.  
Then we conclude that 
$\cO_{\cS}(K_{\cS} + \cY) = \cL \simeq \cO_{\cS}$.  
Therefore 
we see that $K_{\cS} = - \cY$ and hence 
$K_{\cS_{z}} \sim - \cY_z$ for every $z \in \Delta$.   
This implies that $(\cS_z, \cY_z)$ is a rational 
 Okamoto--Painlev\'e pair for $z \in \Delta$.  
Next we claim that if $z \in \Delta- \{ 0 \}$, 
 then  $\dim H^0(\cS_z, \cY_z) = 1$ which also 
implies that $(\cS_z, \cY_z)$ is of non-fibered type.  
  If $\dim H^0(\cS_z, \cY_z) \geq  2 $, 
  we can  show that there exists an elliptic fibration $f_z:\cS_z \lra 
\BP^1$ with $f_z^{*}(\infty) = \cY_z$ which 
is a deformation of the original elliptic 
fibration $f:\cS_0 \lra \BP^1$. 
 Since the general fiber 
$F$ of $f$ does not extend over 
 $ z \in \Delta - \{ 0 \} $, 
 this deduces the contradiction. 
 Note that the type of $(\cS_z, \cY_z) $ 
 is same as the type of $(\cS_0, \cY_0) = (S, Y)$ and 
 $\cS_z - (\cY_z)_{red}$ contains $(-2)$-curves $\cC_z$ 
 whose configuration is 
 same as the configuration of $\cC_0 = C$. \qed.

\vspace{0.5cm}

Now we shall prove the claim (\ref{eq:vanish}).  

\begin{Lemma}\label{lem:h2tlDCE}
Under the same assumption of Lemma \ref{lem:defm nonfibered}, we have
$$
H^2(S, \Theta_S(-\log (D+C+F))= \{ 0 \}, 
$$
where $F$ is a smooth fiber of the elliptic fibration $f:S \lra \BP^1$.  
\end{Lemma}
{\it Proof.} By the Serre duality, it suffices to show that 
\begin{equation}\label{eq:vanishing}
H^0(S, \Omega^1_S( \log (D + C + F)) \otimes K_S ) 
\simeq H^0(S, \Omega^1_S(\log(D+C+F) ) (-F) ) = \{ 0 \}.  
\end{equation}
(Note that $K_S \sim -F $). 
Set $\tilde{D} = \prod_{i=1}^r Y_i, \tilde{C} = \prod_{j=1}^s C_j $.  
Then we have the following commutative diagram of sheaves:
\footnotesize
\begin{equation}\label{eq:diagram}
\begin{array}{ccccccc}
   &          &         &    &     &  0 &         \\ 
   & &         &     &     &  \downarrow &         \\ 
       &       0   &         &   0  & \lra     &  \cO_{F}(-F) &  
       \stackrel{\delta}{\lra}       \\  
       &\downarrow &         & \downarrow    &     &  \downarrow &         \\ 
 0 \lra &  \Omega^1_S(-F)  & \lra & \Omega^1_S  &  \lra 
& (\Omega_S^1)_{|F} & \lra 0   \\ 
     &\downarrow &         & \downarrow    &     &  \quad \downarrow \mu &         \\ 
0 \lra &  \Omega^1_S(\log(D + C+ F))(-F) & \lra & 
\Omega^1_S(\log(D + C+ F))  &  \lra 
& \Omega^1_S(\log( F))_{|F} & \lra  0 \\ 
      &\quad \quad \downarrow P.R. &        
       & \quad \quad \downarrow P.R.   &     &  \downarrow &         \\  
 \stackrel{\delta}{\lra} &
 \oplus_{i=1}^r \cO_{Y_i} \oplus_{j=1}^s \cO_{C_j} \oplus \cO_F(-F) 
  & \lra & \oplus_{i=1}^r \cO_{Y_i} \oplus_{j=1}^s \cO_{C_j} \oplus \cO_F &  \lra & 
 \cO_F & \lra  0  \\
      &\downarrow &         & \downarrow    &     &  \downarrow &         \\  
       &       0   &         &   0  &     &  0 &    .    \\  
          &          &         &    &     &   &         \\ 
 \end{array}
\end{equation}
\normalsize
Here the map $ P.R.:\Omega^1_S(\log(D + C+ F)) \lra 
\oplus_{i=1}^r \cO_{Y_i} \oplus_{j=1}^s \cO_{C_j}$ 
is  the Poincar\'e residue map and the image of 
$\mu:(\Omega^1_{S})_{|F} \lra \Omega^1_S(\log F)_{|F}$ 
 coincides with 
$ \Omega^1_F$ so that the following sequences are exact.
\begin{equation}
\begin{array}{ccccccc}
0 \lra & \cO_F(-F) &  \lra & 
(\Omega^1_S)_{|F} & \lra & \Omega^1_F & \lra 0,  \\
\end{array}
\end{equation}
\begin{equation}
\begin{array}{ccccccc}
 0 \lra & \Omega^1_F &  \lra & 
 (\Omega^1_S(\log F))_{|F} & \lra & \cO_F & \lra 0.\\
\end{array}
\end{equation}
Noting  that $N^{\vee}_F \simeq \cO_F(-F)
 \simeq \cO_F$ and $H^0(\Omega_S^1)= 0$, 
from the first and second rows of (\ref{eq:diagram}),  
we obtain the exact sequence of cohomology
\small
\begin{equation}\label{eq:injective}
\begin{array}{cccc}
       &   0    &   &\\
       & \downarrow &  &\\
 & H^0(\cO_F(-F)) \simeq \C &  &\\
       & \quad \downarrow H_1 &  \searrow & \\
0 \lra & H^0((\Omega_S^1)_{|F}) & \ra& H^1(\Omega^1_S(-F)) 
\end{array}
\end{equation}
\normalsize
From the first column of (\ref{eq:diagram}), 
$H^0(\Omega^1_S(\log(D + C+ F))(-F))$ 
is isomorphic to the kernel of  Gysin map 
\small
\begin{equation}\label{eq:Gysin}
\oplus_{i=1}^r H^0(\cO_{Y_i}) \oplus_{j=1}^s 
H^0(\cO_{C_j}) \oplus H^0(\cO_F(-F)) \stackrel{G_1}{\lra} H^1(\Omega^1_S(-F)),  
\end{equation}
\normalsize
We will show that the Gysin map $G_1$ is injective, which implies the 
assertion (\ref{eq:vanishing}).  

By (\ref{eq:diagram}) and (\ref{eq:injective}),  
we can decompose the map $G_1$ as follows:
\begin{equation}\label{eq:cd2}
\begin{array}{ccc}
  0 &    &           0         \\
 \downarrow   &    &     \downarrow    \\
H^0(\cO_F(-F)) \simeq  \C    & \stackrel{H_1}{\lra} &  
  H^0((\Omega_S^1)_{|F})  \\   
   \quad  \downarrow \tau  &    &  \downarrow      \\
\left. \begin{array}{c}
H^0(\cO_F(-F)) \simeq  \C  \\
\oplus  \\                                
\oplus_{i=1}^r H^0(\cO_{Y_i}) \oplus_{j=1}^s H^0(\cO_{C_j}) 
\end{array} \right\}
 & \stackrel{G_1}{\lra} & 
H^1(\Omega^1_S(-F)) \\
\quad \downarrow \mu_1  &         & \downarrow \nu    \\ 
 \oplus_{i=1}^r H^0(\cO_{Y_i}) \oplus_{j=1}^s H^0(\cO_{C_j}) &   
  \stackrel{G_2}{\lra} &  H^1(\Omega^1_S).  
\end{array}
\end{equation}
Here $\mu_1$ is just the projection and $G_2$ is the natural Gysin map. 
Since $H_1$ is injective (cf. (\ref{eq:injective})), 
a diagram chasing shows that $G_1$ is injective if $G_2$ is injective.   
The image of $1_{Y_i}$ and $1_{C_j}$ by $G_2$ are the class of 
the divisors of $Y_i$ and $C_j$ in $H^1(\Omega^1_S) \simeq H^1(S, \C)$. 
Since $\{  Y_i,  C_j,  1 \leq i \leq r, \quad 1 \leq j \leq s \}$ are
linearly independent in 
$H^1(S, \C) \simeq H^1(\Omega^1_S)$ by assumption 
of Lemma \ref{lem:defm nonfibered}, 
$G_2$ is injective, hence we have proved 
the assertion. \qed.  

\vspace{0.5cm}
Now together with Lemma \ref{lem:defm nonfibered} and 
Proposition \ref{prop:fibered} the following 
lemma shows the existence part of 
Theorem \ref{thm:main-2} and hence  
completes the proof of Theorem \ref{thm:main-2}.  
(See Example \ref{rem:exist}).  

\begin{Lemma}\label{lem:indep}
Let $R$ be a type of affine root lattice in 
Table \ref{tab:type}, that is $R = \tilde{E}_k,
 (k =8,7, 6)$ or $R = \tilde{D}_l, (l =8, 7, 6, 5,4)$. 
Let $(S, Y)$ be a rational  Okamoto--Painlev\'e pair 
of fibered-type and   let $f:S \lra \BP^1$ be 
 the elliptic fibration  with 
 $ f^{*}(\infty) = Y = \sum_{i=1}^r m_i Y_i $. 
  Let $ \{ C_j \}_{j=1}^s $ be a set of different 
irreducible $(-2)$ curves on $S - Y_{red}$ 
such that no linear combination of 
$ \{ C_j \}_{j=1}^s $ has 
the same class of general fiber (= the class of $Y$).  Then 
$ \{  \ Y_i,   C_{j}  \ | \ 1 \leq i \leq r , 1 \leq j \leq s \}$
are linearly independent in $ H^1(S, \Q ) $.  
\end{Lemma}

{\it Proof.}    
From the condition of the set $ \{ C_j \}_{j=1}^{s} $, 
we see that the sub-lattice $ \langle C_j \rangle_{ j=1}^s \subset 
H^2(S, \Z) $ generated by 
$ \{ C_j \}_{j=1}^s $ is negative-definite. Then  we have 
an orthogonal decomposition 
$$
\langle C_j \rangle_{j=1}^s  \oplus \langle Y_i \rangle_{i=1}^r \subset
H^2(S, \Z)  
$$ 
which shows the assertion.  

\qed

\begin{Example}\label{rem:exist}{
\rm
From Proposition \ref{prop:fibered}, we have an Okamoto-Painlev\'e 
pair $(S, Y)$ of fibered type with another singular fiber $F_1$ where the pair
$(Y, F_1)$ has the type $(\tilde{D}_4, \tilde{D}_4)$. 
  Take a proper subset $ \{ C_j \}_{j=1}^s$ of all of 
irreducible components of $F_1$.  Then the type of $M_1$ coincides with 
the proper subgraph of the 
Dynkin diagram $\tilde{D}_4$ of $F_1$, 
that is, one of the types;  
$D_4$, $( A_1, A_1, A_1, A_1)$ $A_3$, 
$(A_1, A_1, A_1)$, $A_2$, $(A_1, A_1)$ 
and  $A_1$.  It is easy to see that the set of classes 
$ \{ Y_i, C_j, 1 \leq i \leq 5, 1 \leq j \leq s \} $ 
are linearly independent in $H^2(S, \Q)$.  
Therefore from Lemma \ref{lem:defm nonfibered} and Lemma 
\ref{lem:indep}, we see that there exists 
a rational Okamoto--Painlev\'e pair 
$(S', Y')$  of non-fibered type, such that:
\begin{enumerate}
\item  the type of $Y'$ is $\tilde{D}_4$, 
\item there exist $(-2)$-curves $\{ C'_j \}_{j=1}^s$ on $S' - Y'_{red}$ with 
the same Dynkin type of  $ \{ C_j \}_{j=1}^s$.
\end{enumerate}
Therefore, we can obtain the assertion of 
Theorem \ref{thm:main-2} for $\tilde{D}_4$. 
 We can treat the other cases similarly.  
}
\end{Example}

\section{\bf Non-existence of Riccati solutions for $P_I$, 
$P_{III}^{ \tilde{D}_8}$, $P_{III}^{\tilde{D}_7}$ } \label{sec:non-exist}

As a corollary to Theorem \ref{thm:main-2}, we obtain the following

\begin{Corollary}\label{cor:non-existence} Let $(S, Y)$ be a rational  
Okamoto--Painlev\'e pair  of non-fibered type,  with the  
type $R = R(Y) = \tilde{E}_8, \tilde{D}_8$ or $\tilde{D}_7$.  
 Then $S - Y_{red}$ does not contain  
a rational nodal curve $C$.  Therefore  all 
the Painlev\'e equations of types 
$P_I$, $P_{III}^{\tilde{D}_8}$, $P_{III}^{\tilde{D}_7}$
 do not admit   
Riccati solutions.  
\end{Corollary} 

{\it Proof.}  The first assertion  directly follows from 
Theorem \ref{thm:main-2} and the last assertion 
follows from the first and Proposition 
\ref{prop:riccati}.  
\qed.

\begin{Remark} {\rm 
\begin{enumerate}
\item 
Umemura proved that the Painlev\'e  equation of type 
$P_I$ has no classical solution 
 and  hence in particular no Riccati solution 
 (cf. \cite{U1}, \cite{U2}).  

\item  Ohyama \cite{Ohy} showed that all the Painlev\'e 
equations of type $P_{III}^{\tilde{D_7}}$ 
has no Riccati solutions by proving that 
they have no invariant divisor with respect to the 
vector field (\ref{eq:vf}).  

\item 
It is worth while remarking that the obstruction to 
the existence of nodal curves in $S - Y_{red}$ is a 
topological one and hence 
so is  the obstruction to the  existence of Riccati solutions.   
In fact, the sub-lattice $M(S - Y_{red})$ is classified only by the 
intersection theory of the surface $S$ and 
the structure of the sub-lattice 
does not depend on the complex structure of $S$.
\end{enumerate} } 
\end{Remark}

%%%
For other types $R$, by the similar argument 
in the proof of Lemma (\ref{lem:h1tlDC}),  
we can show the following porposition.  
This proposition shows that 
for a general parameter ${\balpha} \in \cM_{R}$,  
the corresponding Painlev\'e equations 
do not admit any Riccati solution. 

\begin{Proposition} Let $(S, Y)$ be a 
rational Okamoto--Painlev\'e pair of non-fibered type and of type 
$R$ which corresponds to a Painlev\'e equation and assume that  $S - Y_{red}$
contains a nodal curve $ C$.  
Then there exists a one parameter deformation of 
Okamoto--Painlev\'e pairs of non-fibered type and of the given type $R$, 
$\cY \hookrightarrow \cS \lra \Delta = \{ z \in \C | |z| < \epsilon \}$ 
of $(S, Y)$ such that $\cS_z - \cY_z$ does not contains any nodal curve for 
$z \in \Delta - \{ 0 \}$.  Hence for $z \in \Delta - \{ 0 \}$ the 
Painlev\'e equation corresponding to $(\cS_z, \cY_z)$ does 
not admit any Riccati solutions   
\end{Proposition}

\section{\bf Examples of $(-2)$-curves on $S-D$}\label{sec:example}

In this section, 
we will give  examples of $(-2)$-curves $C$ on $S - D$ for 
some rational Okamoto--Painlev\'e pairs $(S, Y)$ 
and Riccati equations associated to $C$.  

Here we will use the explicit description of  families 
of Okamoto--Painlev\'e pairs 
\begin{equation}\label{eq:globalf2} 
\begin{array}{ccl}
 {\cal S} &  \hookleftarrow & {\cal D}      \\
\pi \downarrow \hspace{0.3cm} &    \swarrow &  \varphi  \\
 \cM_R \times \cB_R  &  &
\end{array} 
\end{equation}
%$\pi:\cS  \lra \cM_R \times \cB_R$  
in \cite{Sa-Te}.  
As we explained in Section \ref{section:riccati}, we have isomorphisms 
$\cM_R = \Spec M_R $ and $\cB_R = \Spec B_R$ such that  
$\Spec M_R$ and $\Spec B_R$ are affine open subschemes of   
$ \Spec \C[\alpha_1, \cdots, \alpha_s] \simeq \C^s $ and 
$\Spec \C[t]$  respectively.  
Moreover 
 $\cS$ can be covered by affine open sets 
$\{ \tilde{U}_i \}_{i=1}^{l+k}$ such that  for each $i$ 
\begin{equation}\label{eq:covering2}
  \tilde{U}_i \   \simeq \  
     \Spec \left( (M_R \otimes B_R)  [x_i, y_i, \frac{1}{f_i(x_i, y_i, \balpha, t)}] \right) 
     \subset   \Spec \C[ \balpha, t, x_i, y_i] \simeq \C^{s+3} \simeq \C^{12 - r}, 
\end{equation}
where $f_i(x_i, y_i, \alpha, t)$ is an element of $(M_R \otimes B_R)  [x_i, y_i] $.  
( Note that in most cases $f_i(x_i, y_i, \alpha, t) \equiv 1$ 
and we may assume that $\cS - \cD$ is covered by $\{ \tilde{U}_i \}_{i=1}^{l}$.) 
%%%

For a given point $\balpha = (\alpha_1, \cdots, \alpha_s) 
\in \cM_R $, we denote the restriction of the family 
$\pi:\cS \lra \cM_R \times \cB_R$ to $ \{ \balpha \}  \times \cB_R$ by 
$\cS_{\balpha} \lra \{ \balpha \} \times \cB_R$.  Moreover we set 
\begin{equation}\label{eq:covering}
{U_i}_{\balpha}:=\tilde{U}_i \cap \cS_{\balpha} \subset  \Spec B_R  
[x_i,y_i], \quad \quad  {U_i}_{(\balpha, t)} = \tilde{U}_i \cap \cS_{\balpha, t}
\subset \Spec \C[x_i, y_i]
\end{equation}\
where $\cS_{\balpha, t}= \pi^{-1}((\balpha, t))$.  

Next let us consider the smooth variety 
obtained by patching affine planes $W_i = \Spec \C[x_i, y_i] \simeq \C^2$ 
( $i =1, 2$) by the coordinate transformation 
\begin{equation}\label{eq:coordinate}
x_1=\frac1{x_2}, \quad y_1={x_2}^2 y_2.
\end{equation}
It is easy to see that the equations $ \{ y_1 = y_2 = 0 \} $  
define a $(-2)$-curve $C$ in $W$.

\begin{Example}[$\tilde E_7$--type ($P_{II}$)]{\rm
In the case of $R = \tilde{E}_7$ ($P_{II}$), the family is constructed as follows 
(cf. \cite{MMT}, \cite{Sa-Te}, \cite{SU}).
Let us set 
$$
\cM_R = \Spec \C[\alpha] \simeq \C, \quad \cB_R = \Spec \C[t] \simeq \C.
$$
Here we only give the affine covering of the family 
$\pi:\cS - \cD \lra \cM_R \times \cB_R$. 
Take three affine schemes $i =0, 1, 2$
\begin{equation}
\tilde{U}_{i} = \Spec \C[\alpha, t, x_i, y_i] \simeq \C^4, 
\end{equation}
and patch these affine schemes by 
the coordinate transformations:
\begin{equation}\label{eq:coor-e7}
\renewcommand{\arraystretch}{1.3}
\begin{array}{lll}
x_0 & = \displaystyle{ \frac{1}{x_1}} & = \displaystyle{\frac{1}{x_2}} ,\\
y_0 & = x_1 ((- \alpha-\frac{1}{2}) - x_1y_1) & 
=2 x_2^{-2} + t + (\alpha-\frac{1}{2})x_2 - y_2 x_2^2 . \\ 
\end{array}
\end{equation}
On $\tilde{U}_0$, the Painlev\'e vector field $\tilde{v}$ 
 in (\ref{eq:vf}) is explicitly given by 
\begin{equation}\label{eq:vf-2}
\tilde{v} = \frac{\partial}{\partial t } + 
\left[ y_0  - x_0^2 - \frac{t}{2} \right] \frac{\partial}{\partial x_0} 
+  \left[ 2x_0 y_0+ \alpha + \frac{1}{2} \right] \frac{\partial}{\partial y_0}.
\end{equation}
which is equivalent to the equation:
\begin{equation}\label{eq:p-2}
\left \{
\begin{array}{ccc}
\displaystyle{\frac{d x_0}{ d t}}&  = &  \displaystyle{ y_0  - x_0^2 - \frac{t}{2} } \\
    &  &  \\
\displaystyle{\frac{d y_0}{d t}} & = &  \displaystyle{2x_0 y_0+ \alpha + \frac{1}{2}}
\end{array} \right. 
\end{equation}

Then for $\alpha = - \frac{1}{2}$, on $U_{0, -\frac{1}{2}} \cup U_{1, -\frac{1}{2}}$, we obtain 
a family of $(-2)$-curves ${\cal C}_{-\frac{1}{2}}  \lra \{- \frac{1}{2} \} \times \cB_{\tilde{E}_7}$ 
defined by 
\begin{equation}\label{eq:-2-curve}
{\cal C}_{-\frac{1}{2}} = \{ y_0 = y_1 = 0 \} \subset U_{0, -\frac{1}{2}} \cup U_{1, -\frac{1}{2}} 
\subset \cS_{- \frac{1}{2}} - \cD_{-\frac{1}{2}}.
\end{equation}
Moreover,   on  the family ${\cal C}_{-\frac{1}{2}}  \lra \{ - \frac{1}{2} \} \times \cB_{\tilde{E}_7}$,  
the equation (\ref{eq:p-2}) can be reduced to  
\begin{equation}\label{eq:p2-riccati}
\displaystyle{\frac{d x_0}{ d t}}  =   \displaystyle{ - x_0^2 - \frac{t}{2} }. 
\end{equation}
It is known that B\"acklund transformations give  isomorphisms between  $\cS_{\alpha}$  and 
$\cS_{\alpha \pm 1}$.  Hence for $ \alpha \in - \frac{1}{2} + \Z $ , the family $\cS_{\alpha} - \cD_{\alpha}$ 
also contains a  family of (-2)-curves (cf. \cite{SU}, \cite{U-W1}).  Moreover,  Noumi and Okamoto \cite{NO} proved the following Theorem (cf. [Theorem 2, \cite{NO}] and remark after it).  
(See also [Theorem 2,1, \cite{U-W1}]).  

\begin{Theorem}\label{thm:u-w} $($ $[$Theorem 2, \cite{NO}$]$$)$.  
Let us denote by $P_{II}(\alpha)$ the equation in (\ref{eq:p-2}).  Then 
\begin{enumerate}
\item For every integer $\alpha \in \Z$, there exists a unique rational solution of the system $P_{II}(\alpha)$. 
\item For every $ \alpha \in \frac{1}{2} + \Z$, there exists a unique one parameter family of classical 
solutions of $P_{II}(\alpha)$, of which each solution is rationally written by a solution of the 
Riccati equation (\ref{eq:p2-riccati}).  
\item Let $(x_0, y_0)$ be a solution of $P_{II}(\alpha)$ different from those mentioned above.  Then neither 
$x_0$ nor $y_0$ is classical, hence  a solution of a Riccati equation. 
\end{enumerate}
\end{Theorem}
Note that for  $\alpha = 0$, $P_{II}(0)$ in (\ref{eq:p-2}) has a rational solution 
$(x_0, y_0) = (0, \frac{t}{2})$.  Theorem \ref{thm:u-w} says that this rational solution 
is the unique rational solution for $P_{II}(0)$. 

}
\end{Example}

\begin{Example}[$\tilde D_4$ ($P_{VI}$)]{\rm
Next let us show examples of $(-2)$-curves for $R= \tilde{D}_4$ (cf. \cite{Sa-Te}).  
The parameter space of the semiuniversal family 
 $\cS - \cD \lra \cM_{\tilde{D}_4} \times  \cB_{\tilde{D}_4}$ are given by 
$$
\cM_R = \Spec \C[\kappa_0, \kappa_1, \kappa_\infty, \kappa_t] 
\simeq \C^4, \quad \cB_R = \Spec \C[t,1/t,1/(t-1)] \simeq \C-\{ 0,1 \}.  
$$
(Here we use the parameters $\kappa_i$, 
$i = 0, 1, \infty, t$ for $\cM_R$ as in \cite{MMT} and \cite{Sa-Te}.)  
Take affine schemes $i =0,1, 2, 3,4,5$
\begin{equation}
\tilde{U}_{i} = \Spec \C[x_i, y_i, \kappa_0, \kappa_1, \kappa_\infty, 
\kappa_t, t,1/t,1/(t-1)] \simeq \C^2 \times \cM_R \times \cB_R.
\end{equation}
and patch them by the coordinate transformations:
\begin{equation}
\renewcommand{\arraystretch}{1.3}
\begin{array}{ll}
x_0 = \displaystyle{ y_1(\kappa_0-x_1 y_1) }, & y_0 = \displaystyle{ \frac{1}{y_1} } ,\\
x_1 = \displaystyle{ y_0(\kappa_0-x_0 y_0) }, & y_1 = \displaystyle{ \frac{1}{y_0} } ,\\
x_0 = \displaystyle{ 1+y_2(\kappa_1-x_2 y_2) } ,& y_0 = \displaystyle{ \frac{1}{y_2} } ,\\
x_2 = \displaystyle{ y_0(\kappa_1+ y_0-x_0 y_0) }, & y_2 = \displaystyle{ \frac{1}{y_0} } ,\\
x_0 = \displaystyle{ t+y_3(\kappa_t-x_3 y_3) }, & y_0 = \displaystyle{ \frac{1}{y_3} } ,\\
x_3 = \displaystyle{ y_0(\kappa_t+ t y_0 -x_0 y_0) }, & y_3 = \displaystyle{ \frac{1}{y_0} } ,\\
x_0 = \displaystyle{ \frac{1}{x_4} } ,& y_0 = \displaystyle{ x_4(\frac{\kappa_0+\kappa_1+\kappa_t -1 +\kappa_\infty}{2} -x_4 y_4) } ,\\
x_4 = \displaystyle{ \frac{1}{x_0} } ,& y_4 = \displaystyle{ x_0(\frac{\kappa_0+\kappa_1+\kappa_t -1 +\kappa_\infty}{2} -x_0 y_0) } ,\\
x_4 = \displaystyle{ y_5(\kappa_\infty -x_5 y_5) } ,& y_4 = \displaystyle{ \frac{1}{y_5} } ,\\
x_5 = \displaystyle{ y_4(\kappa_\infty -x_4 y_4) } ,& y_5 = \displaystyle{ \frac{1}{y_4} }
\end{array}
\end{equation}

On $\tilde{U_0}$, the Painlev\'e vector field $\tilde{v}$ in (\ref{eq:vf}) is given by
\begin{equation}\label{eq:vf-d4}
\tilde{v} = \frac{\partial}{\partial t} + A(x, y, t)  \frac{\partial}{\partial x_0}  + B(x, y, t) \frac{\partial}{\partial y_0},  
\end{equation}
where 
$$
\begin{array}{ll}
A(x, y, t) := & \displaystyle{ \frac{x_0(x_0-1)(x_0-t)}{t(t-1)}  \left[2y_0   -(\frac{\kappa_0}{x_0}+\frac{\kappa_1}{(x_0-1)} +\frac{(\kappa_t -1)}{(x_0-t)}) \right]},  \\
B(x, y, t) := & \displaystyle{ - \frac{1}{t(t-1)}} \left[ (3 x_0^2 -2(t+1)x_0 +t)y_0^2  \right.  \\ 
  & \quad \left. -(2(\kappa_0+\kappa_1+\kappa_t -1)x_0 -(\kappa_0+\kappa_1)t-\kappa_0-\kappa_t +1)y_0  +\frac{(\kappa_0+\kappa_1+\kappa_t -1)^2 -\kappa_\infty^2}{4} \right].
\end{array}
$$
This  is equivalent to the equation:
\begin{equation}\label{eq:p-6}
\left \{
\begin{array}{ccc}
\displaystyle{\frac{d x_0}{ d t}}&  = &  A(x, y, t) \\
    &  &  \\
\displaystyle{\frac{d y_0}{d t}} & = &  B(x, y, t). 
\end{array} \right.
\end{equation}

\normalsize
Let us set the hyperplanes of the parameter space $\cM_{\tilde D_4} \times \cB_{\tilde D_4}$ as follows:
\begin{equation}
\begin{array}{l}
H_0=\{ \kappa_0=0  \}, \quad H_1=\{ \kappa_1=0 \}, \quad H_t=\{ \kappa_t=0 \}, \\
 H_\epsilon=\{ \kappa_0+\kappa_1+\kappa_t+\kappa_\infty-1=0 \}, \quad H_\infty =\{ \kappa_\infty=0 \}.
\end{array}
\end{equation}
Note that each hyperplane $H_i$ is a direct product of $H'_i \subset \cM_{\tilde{D}_4}$ and 
$\cB_{\tilde{D}_4}$, i.e., $H_i = H'_i \times \cB_{\tilde{D}_4}$.  
Remark also that each hyperplane is one of the reflection hyperplanes of the 
affine Weyl group $W(\tilde{D}_4)$ generated by B\"acklund transformations (cf. \cite{NTY}). 
 
We consider the deformation 
$$
\begin{array}{ccc}
\pi ^*(H_0)& \subset & \cS-\cD \\
 \pi \downarrow \hspace{0.3cm} & & \pi \downarrow \hspace{0.3cm} \\
H_0 & \subset & \cM_{\tilde D_4} \times \cB_{\tilde D_4}
\end{array} 
$$
which is given by restricting the parameter space $ \cM_{\tilde D_4} \times \cB_{\tilde D_4}$ to $H_0$.
For subfamily  $(\cS-\cD)_{(0,\kappa_1,\kappa_t,\kappa_\infty)}$ over  $H_0$, 
the coordinate transformation between ${U_0}_{(0,\kappa_1,\kappa_t,\kappa_\infty)}$ and ${U_1}_{(0,\kappa_1,\kappa_t,\kappa_\infty)}$ is given by
$$
x_0=-x_1{y_1}^2, \quad y_0=\frac1{y_1}.
$$
Therefore 
$$
{\cal C}_{0, (0,\kappa_1,\kappa_t,\kappa_\infty)}:=\{ x_0=x_1=0 \}
$$
determines a family of   $(-2)$-curves 
\begin{equation}
\begin{array}{ccc}
{\cal C}_{0, (0,\kappa_1,\kappa_t,\kappa_\infty)} & \hookrightarrow & 
(\cS - \cD)_{(0,\kappa_1,\kappa_t,\kappa_\infty)} \\
\downarrow & \swarrow & \\
H_{0}= H'_0 \times \cB_R. &  &  
\end{array} 
\end{equation} 

In the same way, we obtain families of $(-2)$-curves over each hyperplane $H_i$ as follows:
\begin{equation}
\renewcommand{\arraystretch}{1.2}
\begin{array}{lcll}
H_0 &:& {\cal C}_{0, (0,\kappa_1,\kappa_t,\kappa_\infty)}:=\{ x_0=x_1=0 \} &\subset (\cS-\cD)_{(0,\kappa_1,\kappa_t,\kappa_\infty)}\\
H_1 &:& {\cal C}_{1, (\kappa_0,0,\kappa_t,\kappa_\infty)}:=\{ x_0=1, \ x_2=0 \} &\subset (\cS-\cD)_{(\kappa_0,0,\kappa_t,\kappa_\infty)}\\
H_t &:& {\cal C}_{t, (\kappa_0,\kappa_1,0,\kappa_\infty)}:=\{ x_0=t, \ x_3=0 \} &\subset (\cS-\cD)_{(\kappa_0,\kappa_1,0,\kappa_\infty)}\\
H_\epsilon &:& {\cal C}_{\epsilon, (\kappa_0,\kappa_1,\kappa_t,1-(\kappa_0+\kappa_1+\kappa_t))}:=\{ y_0=y_4=0 \} &\subset (\cS-\cD)_{(\kappa_0,\kappa_1,\kappa_t,1-(\kappa_0+\kappa_1+\kappa_t))} \\
H_\infty &:& {\cal C}_{\infty, (\kappa_0,\kappa_1,\kappa_t,0)}:=\{ x_4=x_5=0 \} &\subset (\cS-\cD)_{(\kappa_0,\kappa_1,\kappa_t,0)}
\end{array}
\end{equation}

\vskp

By restricting the (extended) Hamiltonian system to each $\cC_j$, we obtain the following Riccati equation.
\begin{itemize}
\item On ${\cal C}_{0, (0,\kappa_1,\kappa_t,\kappa_\infty)} \cap {U_0}_{(0,\kappa_1,\kappa_t,\kappa_\infty)}$:
$$
x_0 \equiv 0  ,  \quad
\frac{dy_0}{dt}=-\frac{1}{t(t-1)}( ty_0^2 +(\kappa_1 t+\kappa_t -1)y_0+\frac{(\kappa_1+\kappa_t -1)^2 -\kappa_\infty^2}{4} ) .
$$

\item On ${\cal C}_{1, (\kappa_0,0,\kappa_t,\kappa_\infty)} \cap {U_0}_{(\kappa_0,0,\kappa_t,\kappa_\infty)}$:
$$
x_0 \equiv 1 ,  \quad
\frac{dy_0}{dt}=-\frac{1}{t(t-1)}( (1-t)y_0^2 -((\kappa_0+\kappa_t -1) -\kappa_0 t)y_0+\frac{(\kappa_0+\kappa_t -1)^2 -\kappa_\infty^2}{4} ) .
$$

\item On ${\cal C}_{t, (\kappa_0,\kappa_1,0,\kappa_\infty,t)} \cap {U_0}_{(\kappa_0,\kappa_1,0,\kappa_\infty,t)}$:
$$
x_0  \equiv t ,  \quad
\frac{dy_0}{dt}=-\frac{1}{t(t-1)}( t(t-1)y_0^2 -((\kappa_0+\kappa_1 -2)t -\kappa_0 +1)y_0+\frac{(\kappa_0+\kappa_1 -1)^2 -\kappa_\infty^2}{4} ) .
$$

\item On ${\cal C}_{\epsilon, (\kappa_0,\kappa_1,\kappa_t,1-(\kappa_0+\kappa_1+\kappa_t),t)} \cap {U_0}_{(\kappa_0,\kappa_1,\kappa_t,1-(\kappa_0+\kappa_1+\kappa_t),t)}$:
$$
\frac{dx_0}{dt}=-\frac{1}{t(t-1)}(\kappa_0(x_0-1)(x_0-t)+\kappa_1 x_0(x_0-t) +(\kappa_t -1)x_0(x_0-1)) ,  \quad
y_0 \equiv 0.
$$

\item On ${\cal C}_{\infty, (\kappa_0,\kappa_1,\kappa_t,0,t)} \cap {U_4}_{(\kappa_0,\kappa_1,\kappa_t,0,t)}$:
$$
x_4 \equiv 0 ,  \quad
\frac{dy_4}{dt}=-\frac{1}{t(t-1)}(  y_4^2 +( (\kappa_t-1)t+\kappa_1)y_4      +\frac{(\kappa_1+\kappa_t -1)^2 -\kappa_0^2}{4} t )  .
$$

\end{itemize}

\vskp
Next, choose four hyperplanes from the five hyperplanes and  consider the fibers over 
the intersection of them. 
For each fiber $(\cS-\cD)_{(0,0,0,1,t)}$ of $(0,0,0,1,t) \in H_0 \cap H_1 \cap H_t \cap H_\epsilon$, we can see that ${\cal C}_{0, (0,0,0,1,t)},{\cal C}_{1, (0,0,0,1,t)},$ and ${\cal C}_{t, (0,0,0,1,t)}$ do not intersect each other but they intersect with ${\cal C}_{\epsilon, (0,0,0,1,t)}$ respectively. Hence the type of the configuration of these curves is $D_4$. By checking the other cases, we obtain the following.
\begin{center}
\renewcommand{\arraystretch}{1.2}
\begin{tabular}{|c|c|c|}
\hline
fiber & $(-2)$-curves & configuration \\
\noalign{\hrule height 0.8pt}
$(\cS-\cD)_{(0,0,0,1)}$ & $\{ \cC_0, \cC_1 , \cC_t, \cC_\epsilon \}$ & $D_4$ \\
\hline
$(\cS-\cD)_{(0,0,1,0)}$ & $\{ \cC_0, \cC_1 , \cC_\epsilon, \cC_\infty \}$ & $D_4$ \\
\hline
$(\cS-\cD)_{(0,1,0,0)}$ & $\{ \cC_0, \cC_t, \cC_\epsilon, \cC_\infty \}$ & $D_4$ \\
\hline
$(\cS-\cD)_{(1,0,0,0)}$ & $\{ \cC_1, \cC_t, \cC_\epsilon, \cC_\infty \}$ & $D_4$ \\
\hline
$(\cS-\cD)_{(0,0,0,0)}$ & $\{ \cC_0, \cC_1 , \cC_t, \cC_\infty \}$ & $A_1,A_1,A_1,A_1$ \\
\hline
\end{tabular}
\end{center}

\begin{figure}[h]
\begin{center}
%WinTpicVersion2.15
\unitlength 0.1in
\begin{picture}(46.90,35.14)(0.10,-39.43)
% LINE 0 0 3 0
% 2 220 1048 1991 1048
% 
\special{pn 20}%
\special{pa 220 648}%
\special{pa 1991 648}%
\special{fp}%
% LINE 0 0 3 0
% 2 507 915 178 1727
% 
\special{pn 20}%
\special{pa 507 515}%
\special{pa 178 1327}%
\special{fp}%
% LINE 0 0 3 0
% 2 920 915 920 1741
% 
\special{pn 20}%
\special{pa 920 515}%
\special{pa 920 1341}%
\special{fp}%
% LINE 0 2 3 0
% 2 1753 922 1984 1741
% 
\special{pn 20}%
\special{pa 1753 522}%
\special{pa 1984 1341}%
\special{dt 0.054}%
\special{pa 1984 1341}%
\special{pa 1984 1340}%
\special{dt 0.054}%
% STR 2 0 3 0
% 3 2019 936 2019 1006 2 0
% $\cC_{\epsilon}$
\put(20.1900,-6.0600){\makebox(0,0)[lb]{$\cC_{\epsilon}$}}%
% STR 2 0 3 0
% 3 80 1846 80 1916 2 0
% $\cC_{0}$
\put(0.8000,-15.1600){\makebox(0,0)[lb]{$\cC_{0}$}}%
% STR 2 0 3 0
% 3 724 1853 724 1923 2 0
% $\cC_{1}$
\put(7.2400,-15.2300){\makebox(0,0)[lb]{$\cC_{1}$}}%
% STR 2 0 3 0
% 3 1382 1853 1382 1923 2 0
% $\cC_{t}$
\put(13.8200,-15.2300){\makebox(0,0)[lb]{$\cC_{t}$}}%
% STR 2 0 3 0
% 3 2054 1853 2054 1923 2 0
% $\cC_{\infty}$
\put(20.5400,-15.2300){\makebox(0,0)[lb]{$\cC_{\infty}$}}%
% LINE 0 0 3 0
% 2 1340 908 1340 1734
% 
\special{pn 20}%
\special{pa 1340 508}%
\special{pa 1340 1334}%
\special{fp}%
% LINE 0 0 3 0
% 2 2845 1041 4616 1041
% 
\special{pn 20}%
\special{pa 2845 641}%
\special{pa 4616 641}%
\special{fp}%
% LINE 0 0 3 0
% 2 3132 908 2803 1720
% 
\special{pn 20}%
\special{pa 3132 508}%
\special{pa 2803 1320}%
\special{fp}%
% LINE 0 0 3 0
% 2 3545 908 3545 1734
% 
\special{pn 20}%
\special{pa 3545 508}%
\special{pa 3545 1334}%
\special{fp}%
% LINE 0 0 3 0
% 2 4378 915 4609 1734
% 
\special{pn 20}%
\special{pa 4378 515}%
\special{pa 4609 1334}%
\special{fp}%
% STR 2 0 3 0
% 3 4644 929 4644 999 2 0
% $\cC_{\epsilon}$
\put(46.4400,-5.9900){\makebox(0,0)[lb]{$\cC_{\epsilon}$}}%
% STR 2 0 3 0
% 3 2705 1839 2705 1909 2 0
% $\cC_{0}$
\put(27.0500,-15.0900){\makebox(0,0)[lb]{$\cC_{0}$}}%
% STR 2 0 3 0
% 3 3349 1846 3349 1916 2 0
% $\cC_{1}$
\put(33.4900,-15.1600){\makebox(0,0)[lb]{$\cC_{1}$}}%
% STR 2 0 3 0
% 3 4007 1846 4007 1916 2 0
% $\cC_{t}$
\put(40.0700,-15.1600){\makebox(0,0)[lb]{$\cC_{t}$}}%
% STR 2 0 3 0
% 3 4679 1846 4679 1916 2 0
% $\cC_{\infty}$
\put(46.7900,-15.1600){\makebox(0,0)[lb]{$\cC_{\infty}$}}%
% LINE 0 2 3 0
% 2 3965 901 3965 1727
% 
\special{pn 20}%
\special{pa 3965 501}%
\special{pa 3965 1327}%
\special{dt 0.054}%
\special{pa 3965 1327}%
\special{pa 3965 1326}%
\special{dt 0.054}%
% LINE 0 0 3 0
% 2 150 2364 1921 2364
% 
\special{pn 20}%
\special{pa 150 1964}%
\special{pa 1921 1964}%
\special{fp}%
% LINE 0 0 3 0
% 2 437 2231 108 3043
% 
\special{pn 20}%
\special{pa 437 1831}%
\special{pa 108 2643}%
\special{fp}%
% LINE 0 2 3 0
% 2 850 2231 850 3057
% 
\special{pn 20}%
\special{pa 850 1831}%
\special{pa 850 2657}%
\special{dt 0.054}%
\special{pa 850 2657}%
\special{pa 850 2656}%
\special{dt 0.054}%
% LINE 0 0 3 0
% 2 1683 2238 1914 3057
% 
\special{pn 20}%
\special{pa 1683 1838}%
\special{pa 1914 2657}%
\special{fp}%
% STR 2 0 3 0
% 3 1949 2252 1949 2322 2 0
% $\cC_{\epsilon}$
\put(19.4900,-19.2200){\makebox(0,0)[lb]{$\cC_{\epsilon}$}}%
% STR 2 0 3 0
% 3 10 3162 10 3232 2 0
% $\cC_{0}$
\put(0.1000,-28.3200){\makebox(0,0)[lb]{$\cC_{0}$}}%
% STR 2 0 3 0
% 3 654 3169 654 3239 2 0
% $\cC_{1}$
\put(6.5400,-28.3900){\makebox(0,0)[lb]{$\cC_{1}$}}%
% STR 2 0 3 0
% 3 1312 3169 1312 3239 2 0
% $\cC_{t}$
\put(13.1200,-28.3900){\makebox(0,0)[lb]{$\cC_{t}$}}%
% STR 2 0 3 0
% 3 1984 3169 1984 3239 2 0
% $\cC_{\infty}$
\put(19.8400,-28.3900){\makebox(0,0)[lb]{$\cC_{\infty}$}}%
% LINE 0 0 3 0
% 2 1270 2224 1270 3050
% 
\special{pn 20}%
\special{pa 1270 1824}%
\special{pa 1270 2650}%
\special{fp}%
% LINE 0 0 3 0
% 2 2866 2371 4637 2371
% 
\special{pn 20}%
\special{pa 2866 1971}%
\special{pa 4637 1971}%
\special{fp}%
% LINE 0 2 3 0
% 2 3153 2238 2824 3050
% 
\special{pn 20}%
\special{pa 3153 1838}%
\special{pa 2824 2650}%
\special{dt 0.054}%
\special{pa 2824 2650}%
\special{pa 2824 2649}%
\special{dt 0.054}%
% LINE 0 0 3 0
% 2 3566 2238 3566 3064
% 
\special{pn 20}%
\special{pa 3566 1838}%
\special{pa 3566 2664}%
\special{fp}%
% LINE 0 0 3 0
% 2 4399 2245 4630 3064
% 
\special{pn 20}%
\special{pa 4399 1845}%
\special{pa 4630 2664}%
\special{fp}%
% STR 2 0 3 0
% 3 4665 2259 4665 2329 2 0
% $\cC_{\epsilon}$
\put(46.6500,-19.2900){\makebox(0,0)[lb]{$\cC_{\epsilon}$}}%
% STR 2 0 3 0
% 3 2726 3169 2726 3239 2 0
% $\cC_{0}$
\put(27.2600,-28.3900){\makebox(0,0)[lb]{$\cC_{0}$}}%
% STR 2 0 3 0
% 3 3370 3176 3370 3246 2 0
% $\cC_{1}$
\put(33.7000,-28.4600){\makebox(0,0)[lb]{$\cC_{1}$}}%
% STR 2 0 3 0
% 3 4028 3176 4028 3246 2 0
% $\cC_{t}$
\put(40.2800,-28.4600){\makebox(0,0)[lb]{$\cC_{t}$}}%
% STR 2 0 3 0
% 3 4700 3176 4700 3246 2 0
% $\cC_{\infty}$
\put(47.0000,-28.4600){\makebox(0,0)[lb]{$\cC_{\infty}$}}%
% LINE 0 0 3 0
% 2 3986 2231 3986 3057
% 
\special{pn 20}%
\special{pa 3986 1831}%
\special{pa 3986 2657}%
\special{fp}%
% LINE 0 2 3 0
% 2 1445 3638 3216 3638
% 
\special{pn 20}%
\special{pa 1445 3238}%
\special{pa 3216 3238}%
\special{dt 0.054}%
\special{pa 3216 3238}%
\special{pa 3215 3238}%
\special{dt 0.054}%
% LINE 0 0 3 0
% 2 1732 3505 1403 4317
% 
\special{pn 20}%
\special{pa 1732 3105}%
\special{pa 1403 3917}%
\special{fp}%
% LINE 0 0 3 0
% 2 2145 3505 2145 4331
% 
\special{pn 20}%
\special{pa 2145 3105}%
\special{pa 2145 3931}%
\special{fp}%
% LINE 0 0 3 0
% 2 2978 3512 3209 4331
% 
\special{pn 20}%
\special{pa 2978 3112}%
\special{pa 3209 3931}%
\special{fp}%
% STR 2 0 3 0
% 3 3244 3526 3244 3596 2 0
% $\cC_{\epsilon}$
\put(32.4400,-31.9600){\makebox(0,0)[lb]{$\cC_{\epsilon}$}}%
% STR 2 0 3 0
% 3 1305 4436 1305 4506 2 0
% $\cC_{0}$
\put(13.0500,-41.0600){\makebox(0,0)[lb]{$\cC_{0}$}}%
% STR 2 0 3 0
% 3 1949 4443 1949 4513 2 0
% $\cC_{1}$
\put(19.4900,-41.1300){\makebox(0,0)[lb]{$\cC_{1}$}}%
% STR 2 0 3 0
% 3 2607 4443 2607 4513 2 0
% $\cC_{t}$
\put(26.0700,-41.1300){\makebox(0,0)[lb]{$\cC_{t}$}}%
% STR 2 0 3 0
% 3 3279 4443 3279 4513 2 0
% $\cC_{\infty}$
\put(32.7900,-41.1300){\makebox(0,0)[lb]{$\cC_{\infty}$}}%
% LINE 0 0 3 0
% 2 2565 3498 2565 4324
% 
\special{pn 20}%
\special{pa 2565 3098}%
\special{pa 2565 3924}%
\special{fp}%
\end{picture}%
\end{center}
\caption{Maximal configurations for $R = \tilde{D_4}$.}
\label{fig:d-4max}
\end{figure}
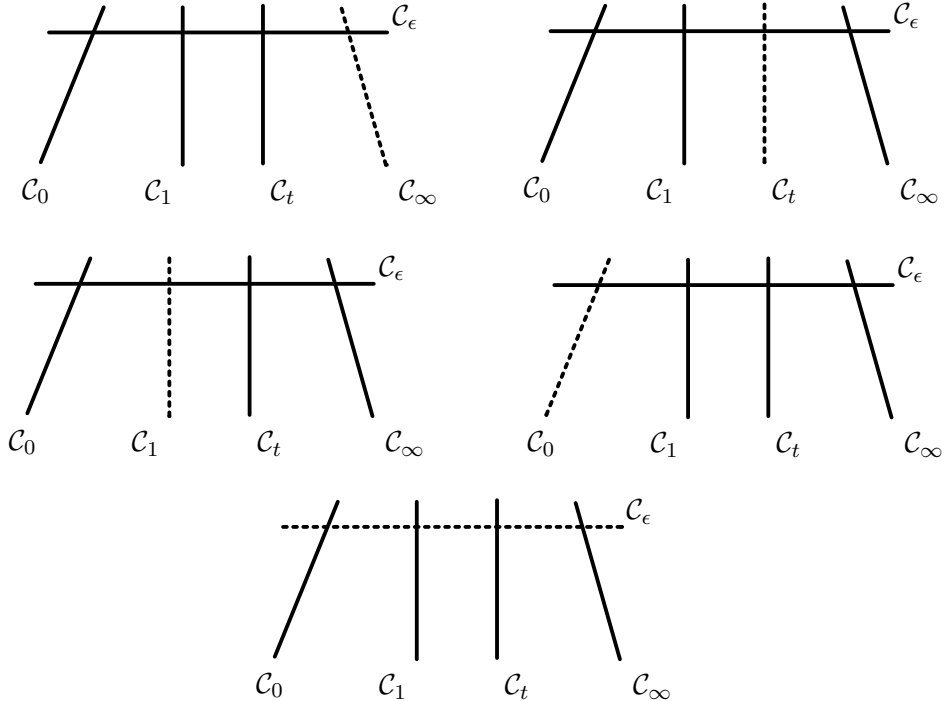

}
\end{Example}

Below, we only give the tables for $\tilde{D}_k, k =5, 6$.  
The case $ \tilde{E}_6$ will be treated in Section \ref{sec:conf}.  
For  parameters and the coordinate transformations, see \cite{Sa-Te}. 

\begin{Example}[$\tilde D_5$ ($P_{V}$)]{\rm
$$
\cM_{\tilde{D}_5} = \Spec \C[ \kappa_0, \kappa_t, \kappa_{\infty}] \simeq \C^3,  
 \quad \cB_{\tilde{D}_5} =\Spec \C[t, t^{-1}] \simeq \C^{\times}. 
$$  
$$
\renewcommand{\arraystretch}{1.2}
\begin{array}{lcll}
H_0=\{ \kappa_0=0  \} &:& {\cal C}_{0, (0,\kappa_t,\kappa_\infty)}:=\{ x_0=x_1=0 \} &\subset (\cS-\cD)_{(0,\kappa_t,\kappa_\infty)}\\
H_\epsilon=\{ \kappa_0+\kappa_t+\kappa_\infty=0  \} &:& {\cal C}_{\epsilon, (\kappa_0,\kappa_t,-(\kappa_0+\kappa_t))}:=\{ y_0=y_3=0 \} &\subset (\cS-\cD)_{(\kappa_0,\kappa_t,-(\kappa_0+\kappa_t))} \\
H_\infty=\{ \kappa_\infty=0  \} &:& {\cal C}_{\infty, (\kappa_0,\kappa_t,0)}:=\{ x_3=x_4=0 \} &\subset (\cS-\cD)_{(\kappa_0,\kappa_t,0)}
\end{array}
$$

\begin{itemize}
\item On $ {\cal C}_{0, (0,\kappa_t,\kappa_\infty)}\cap {U_0}_{(0,\kappa_t,\kappa_\infty)}$
$$
x_0 \equiv 0 ,  \quad
\frac{dy_0}{dt}=-\frac{1}{t}(y_0^2+(\kappa_t - t)y_0+\frac{\kappa_t^2-\kappa_\infty^2}{4})      .
$$

\item On ${\cal C}_{\epsilon, (\kappa_0,\kappa_t,-(\kappa_0+\kappa_t))} \cap {U_0}_{(\kappa_0,\kappa_t,-(\kappa_0+\kappa_t))}$
$$
\frac{dx_0}{dt}=-\frac{1}{t}(\kappa_0 (x_0-1)^2 +\kappa_t x_0(x_0-1)+ t x_0) ,  \quad
y_0 \equiv 0 .
$$

\item On ${\cal C}_{\infty, (\kappa_0,\kappa_t,0)} \cap {U_3}_{(\kappa_0,\kappa_t,0)}$
$$
x_3 \equiv 0  ,  \quad
\frac{dy_3}{dt}=-\frac{1}{t}(y_3^2 +( \kappa_t +t)y_3 +\frac{\kappa_t^2-\kappa_0^2}{4}  ).
$$

\end{itemize}

\vskp

\begin{center}
\renewcommand{\arraystretch}{1.2}
\begin{tabular}{|c|c|c|}
\hline
fiber & $(-2)$-curves & configuration \\
\noalign{\hrule height 0.8pt}
$(\cS-\cD)_{(0,0,0)}$ & $\{ \cC_0, \cC\epsilon ,\cC_\infty \}$ & $A_3$ \\
\hline
\end{tabular}
\end{center}
}
\end{Example}

\vskp

\begin{Example}[$\tilde D_6$ ($P_{III}$)]{\rm

$$
\cM_R = \Spec \C[\kappa_0,\kappa_\infty] \simeq \C^2, \quad \cB_R = \Spec \C[t,t^{-1}] \simeq  \C^\times. 
$$

$$
\renewcommand{\arraystretch}{1.2}
\begin{array}{lcll}
H_1=\{ \kappa_0+\kappa_\infty=0  \} &:& {\cal C}_{1, (\kappa_0,-\kappa_0)}:=\{ y_0=y_2=0 \} & \subset (\cS-\cD)_{(\kappa_0,-\kappa_0)} \\
H_2=\{ \kappa_0-\kappa_\infty=0  \} &:& {\cal C}_{2, (\kappa_0,\kappa_0)}:=\{ y_0=t,y_3=0 \} &\subset (\cS-\cD)_{(\kappa_0,\kappa_0)}\\
H_3=\{ \kappa_0-\kappa_\infty+2=0  \} &:& {\cal C}_{3, (\kappa_0,\kappa_0+2)}:=\{ y_1=0,y_2=t \} &\subset (\cS-\cD)_{(\kappa_0,\kappa_0+2)}\\
H_4=\{ \kappa_0+\kappa_\infty+2=0  \} &:& {\cal C}_{4, (\kappa_0,-\kappa_0-2)}:=\{ y_1=t,y_3=t \} &\subset (\cS-\cD)_{(\kappa_0,-\kappa_0-2)}
\end{array}
$$

\vskp

\begin{itemize}
\item On ${\cal C}_{1, (\kappa_0,-\kappa_0)} \cap {U_0}_{(\kappa_0,-\kappa_0)}$
$$
\frac{dx_0}{dt}=\frac{1}{t}(-2 t x_0^2-(2 \kappa_0+1) x_0+2 t) ,  \quad
y_0 \equiv 0  .
$$

\item On ${\cal C}_{2, (\kappa_0,\kappa_0)} \cap {U_0}_{(\kappa_0,\kappa_0)}$
$$
\frac{dx_0}{dt}=\frac{1}{t}(2 t x_0^2-(2 \kappa_0+1) x_0+2 t) ,  \quad
y_0 \equiv t.
$$

\item On ${\cal C}_{3, (\kappa_0,\kappa_0+2)} \cap {U_1}_{(\kappa_0,\kappa_0+2)}$
$$
\frac{dx_1}{dt}=\frac{1}{t}( -2 t x_1^2 +(2\kappa_0+3)x_1 -2 t) ,  \quad
y_1 \equiv 0 .
$$

\item On ${\cal C}_{4, (\kappa_0,-\kappa_0-2)} \cap {U_1}_{(\kappa_0,-\kappa_0-2)}$
$$
\frac{dx_1}{dt}=\frac{1}{t}(2 t x_1^2 +(2\kappa_0+3)x_1 -2 t) ,  \quad
y_1 \equiv 1.
$$

\end{itemize}

\vskp

\begin{center}
\renewcommand{\arraystretch}{1.2}
\begin{tabular}{|c|c|c|}
\hline
fiber & $(-2)$-curves & configuration \\
\noalign{\hrule height 0.8pt}
$(\cS-\cD)_{(0,0)}$ & $\{ \cC_1, \cC_2 \}$ & $A_1,A_1$ \\
\hline
$(\cS-\cD)_{(-1,1)}$ & $\{ \cC_1, \cC_3 \}$ & $A_1,A_1$ \\
\hline
$(\cS-\cD)_{(-1,-1)}$ & $\{ \cC_2, \cC_4 \}$ & $A_1,A_1$ \\
\hline
$(\cS-\cD)_{(-2,0)}$ & $\{ \cC_3, \cC_4 \}$ & $A_1,A_1$ \\
\hline
\end{tabular}
\end{center}

}
\end{Example}

\vskp

\section{\bf Confluences of Nodal Curves and Riccati Equations}\label{sec:conf}

In this section,  we will discuss   the confluence of nodal curves and Riccati 
equations for Painlev\'e equations.  
We will deal with only the case $R = \tilde{E_6}$ ($P_{IV}$), however one can easily extend the 
result to  other cases like $\tilde{D_5}$ and $\tilde{D_4}$.

\subsection{The confluence of nodal curves}  
\begin{Example}[$\tilde E_6$ ($P_{IV}$)]{\rm
$$
\cM_R = \Spec \C[\kappa_0, \kappa_\infty ] \simeq \C^2, \quad \cB_R = \Spec \C[t] \simeq \C.  
$$
An open covering of $\cS-\cD$ is given by 
$$
\cS - \cD = \bigcup_{i=0}^3 \tilde{U}_{i}
$$
where for $i = 0, 1, 2, 3$ 
$$
\tilde{U}_{i} = \Spec \C[ x_i, y_i, \kappa_0, \kappa_\infty, t ] \simeq \C^5.  
$$
Moreover the coordinate transformations are given by
$$
\renewcommand{\arraystretch}{1.3}
\begin{array}{ll}
x_0 = \displaystyle{ y_1(\kappa_0-x_1 y_1) }, & y_0 = \displaystyle{ \frac{1}{y_1} }, \\
x_1 = \displaystyle{ y_0(\kappa_0-x_0 y_0) }, & y_1 = \displaystyle{ \frac{1}{y_0} }, \\
x_0 = \displaystyle{ \frac{1}{x_2} }, & y_0 = \displaystyle{ x_2(\kappa_\infty-x_2 y_2) }, \\
x_2 = \displaystyle{ \frac{1}{x_0} }, & y_2= \displaystyle{ x_0(\kappa_\infty-x_0 y_0) }, \\

x_2 = \displaystyle{ x_3 }, & y_2= \displaystyle{ -\frac{1/2}{x_3^3}-\frac{t}{x_3^2}+\frac{2\kappa_\infty -\kappa_0+1}{x_3}+y_3 }, \\
x_3 = \displaystyle{ x_2 }, & y_3= \displaystyle{ \frac{1/2}{x_2^3}+\frac{t}{x_2^2}-\frac{2\kappa_\infty -\kappa_0+1}{x_2}+y_2 }.
\end{array}
$$
Finally, on the affine open set $\tilde{U}_0$, 
the Painlev\'e system of type $\tilde{E_6}$ which is equivalent to $P_{IV}$ is given  as follows. 
\begin{equation}\label{eq:p4}
\left\{
\begin{array}{ccl}
\displaystyle{\frac{d x_0 }{d t}} &  = & 4x_0 y_0 -x_0^2 -2 t x_0-2 \kappa_0 \\
    &  & \\
\displaystyle{\frac{d y_0 }{d t}} & = & -2y_0^2 +2(x_0+t)y_0 -\kappa_\infty
\end{array}
\right..
\end{equation}

We have two  hyperplanes $H_0$ and $H_{\infty}$ on $\cM_R \times \cB_R$ and 
families of $(-2)$-curves $\cC_{0}$ and $\cC_{\infty}$ over $H_0$ and $H_{\infty}$ as 
follows. 
\begin{equation}\label{eq:e6-h}
\renewcommand{\arraystretch}{1.2}
\begin{array}{lcll}
H_0=\{ \kappa_0=0  \} &:& \cC_{0, (0,\kappa_\infty)}:=\{ x_0=x_1=0 \} &\subset (\cS-\cD)_{(0,\kappa_\infty)}\\
H_\infty=\{ \kappa_\infty=0  \} &:& \cC_{\infty, (\kappa_0,0)}:=\{ y_0=y_2=0 \} &\subset (\cS-\cD)_{(\kappa_0,0)}
\end{array}
\end{equation}

Then now it is easy to see that the Painlev\'e system (\ref{eq:p4}) can be reduced to 
the following Riccati equations on $\cC_{0}$ and $\cC_{\infty}$ respectively.

\begin{itemize}
\item On ${\cal C}_{0, (0,\kappa_\infty)} \cap {U_0}_{(0,\kappa_\infty)}$
\begin{equation}\label{eq:p4-0}
x_0 \equiv 0 ,  \quad
\frac{dy_0}{dt}=-2y_0^2 +2 t y_0 -\kappa_\infty   .
\end{equation}

\item On ${\cal C}_{\infty, (\kappa_0,0)} \cap {U_0}_{(\kappa_0,0)}$
\begin{equation}\label{eq:p4-inf}
\frac{dx_0}{dt}=-x_0^2 -2 t x_0-2 \kappa_0 ,  \quad
y_0 \equiv 0      .
\end{equation}
\end{itemize}

\begin{center}
\renewcommand{\arraystretch}{1.2}
\begin{tabular}{|c|c|c|}
\hline
fiber & $(-2)$-curves & configuration \\
\noalign{\hrule height 0.8pt}
$(\cS-\cD)_{(0,0)}$ & $\{ \cC_0, \cC_\infty \}$ & $A_2$ \\
\hline
\end{tabular}
\end{center}

}
\end{Example}

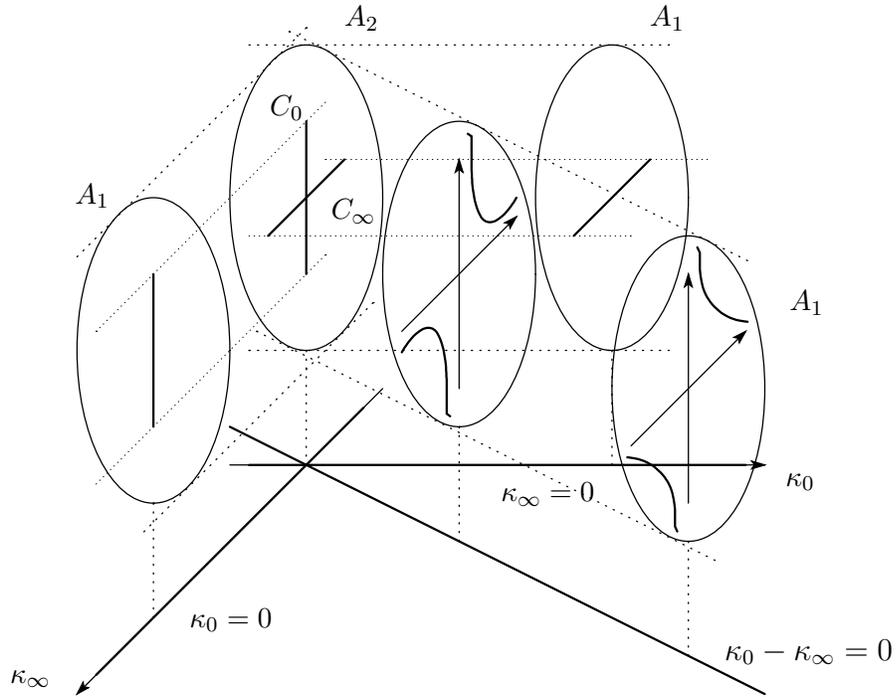
\begin{figure}
\begin{center}%WinTpicVersion2.15
\unitlength 0.1in
\begin{picture}(40.80,35.90)(2.50,-38.00)
% VECTOR 2 0 3 0
% 4 1400 3000 4200 3000 2200 2600 600 4200
% 
\special{pn 8}%
\special{pa 1400 2600}%
\special{pa 4200 2600}%
\special{fp}%
\special{sh 1}%
\special{pa 4200 2600}%
\special{pa 4133 2580}%
\special{pa 4147 2600}%
\special{pa 4133 2620}%
\special{pa 4200 2600}%
\special{fp}%
\special{pa 2200 2200}%
\special{pa 600 3800}%
\special{fp}%
\special{sh 1}%
\special{pa 600 3800}%
\special{pa 661 3767}%
\special{pa 638 3762}%
\special{pa 633 3739}%
\special{pa 600 3800}%
\special{fp}%
% LINE 1 0 3 0
% 2 1500 3000 4100 3000
% 
\special{pn 13}%
\special{pa 1500 2600}%
\special{pa 4100 2600}%
\special{fp}%
% LINE 1 0 3 0
% 2 2100 2700 700 4100
% 
\special{pn 13}%
\special{pa 2100 2300}%
\special{pa 700 3700}%
\special{fp}%
% LINE 2 2 3 0
% 2 1800 3000 1800 2400
% 
\special{pn 8}%
\special{pa 1800 2600}%
\special{pa 1800 2000}%
\special{dt 0.045}%
\special{pa 1800 2000}%
\special{pa 1800 2001}%
\special{dt 0.045}%
% ELLIPSE 2 0 3 1
% 4 1800 1600 2200 2400 2200 2400 2200 2400
% 
\special{pn 8}%
\special{ar 1800 1200 400 800  0.0000000 6.2831853}%
% LINE 1 0 3 2
% 4 1800 1200 1800 2000 1600 1800 2000 1400
% 
\special{pn 13}%
\special{pa 1800 800}%
\special{pa 1800 1600}%
\special{fp}%
\special{pa 1600 1400}%
\special{pa 2000 1000}%
\special{fp}%
% LINE 2 2 3 0
% 2 1000 3800 1000 3200
% 
\special{pn 8}%
\special{pa 1000 3400}%
\special{pa 1000 2800}%
\special{dt 0.045}%
\special{pa 1000 2800}%
\special{pa 1000 2801}%
\special{dt 0.045}%
% ELLIPSE 2 0 3 0
% 4 1000 2400 1400 3200 1400 3200 1400 3200
% 
\special{pn 8}%
\special{ar 1000 2000 400 800  0.0000000 6.2831853}%
% LINE 1 0 3 0
% 2 1000 2000 1000 2800
% 
\special{pn 13}%
\special{pa 1000 1600}%
\special{pa 1000 2400}%
\special{fp}%
% LINE 2 2 3 0
% 2 3400 3000 3400 2400
% 
\special{pn 8}%
\special{pa 3400 2600}%
\special{pa 3400 2000}%
\special{dt 0.045}%
\special{pa 3400 2000}%
\special{pa 3400 2001}%
\special{dt 0.045}%
% ELLIPSE 2 0 3 0
% 4 3400 1600 3800 2400 3800 2400 3800 2400
% 
\special{pn 8}%
\special{ar 3400 1200 400 800  0.0000000 6.2831853}%
% LINE 1 0 3 0
% 2 3200 1800 3600 1400
% 
\special{pn 13}%
\special{pa 3200 1400}%
\special{pa 3600 1000}%
\special{fp}%
% LINE 3 2 3 0
% 2 1500 1800 3500 1800
% 
\special{pn 4}%
\special{pa 1500 1400}%
\special{pa 3500 1400}%
\special{dt 0.027}%
\special{pa 3500 1400}%
\special{pa 3499 1400}%
\special{dt 0.027}%
% LINE 3 2 3 0
% 2 700 3100 1900 1900
% 
\special{pn 4}%
\special{pa 700 2700}%
\special{pa 1900 1500}%
\special{dt 0.027}%
\special{pa 1900 1500}%
\special{pa 1900 1500}%
\special{dt 0.027}%
% LINE 3 2 3 0
% 2 1900 1100 700 2300
% 
\special{pn 4}%
\special{pa 1900 700}%
\special{pa 700 1900}%
\special{dt 0.027}%
\special{pa 700 1900}%
\special{pa 700 1900}%
\special{dt 0.027}%
% LINE 3 2 3 0
% 2 3900 1400 1900 1400
% 
\special{pn 4}%
\special{pa 3900 1000}%
\special{pa 1900 1000}%
\special{dt 0.027}%
\special{pa 1900 1000}%
\special{pa 1901 1000}%
\special{dt 0.027}%
% STR 2 0 3 0
% 3 1610 1070 1610 1170 2 0
% c0
\put(16.1000,-7.7000){\makebox(0,0)[lb]{$C_0$}}%
% STR 2 0 3 0
% 3 1930 1610 1930 1710 2 0
% ci
\put(19.3000,-13.1000){\makebox(0,0)[lb]{$C_\infty$}}%
% STR 2 0 3 0
% 3 250 4040 250 4140 2 0
% ki
\put(2.5000,-37.4000){\makebox(0,0)[lb]{$\kappa_\infty$}}%
% STR 2 0 3 0
% 3 4310 3010 4310 3110 2 0
% k0
\put(43.1000,-27.1000){\makebox(0,0)[lb]{$\kappa_0$}}%
% STR 2 0 3 0
% 3 2820 3090 2820 3190 2 0
% ki=0
\put(28.2000,-27.9000){\makebox(0,0)[lb]{$\kappa_\infty=0$}}%
% STR 2 0 3 0
% 3 1190 3750 1190 3850 2 0
% k0=0
\put(11.9000,-34.5000){\makebox(0,0)[lb]{$\kappa_0=0$}}%
% STR 2 0 3 0
% 3 2000 680 2000 780 2 0
% a2
\put(20.0000,-3.0000){\makebox(0,0)[lb]{$A_2$}}%
% STR 2 0 3 0
% 3 590 1530 590 1630 2 0
% a1
\put(5.9000,-12.3000){\makebox(0,0)[lb]{$A_1$}}%
% STR 2 0 3 0
% 3 3600 690 3600 790 2 0
% a1
\put(36.0000,-3.0000){\makebox(0,0)[lb]{$A_1$}}%
% LINE 2 2 3 0
% 4 1500 800 3700 800 3700 2400 1500 2400
% 
\special{pn 8}%
\special{pa 1500 400}%
\special{pa 3700 400}%
\special{dt 0.045}%
\special{pa 3700 400}%
\special{pa 3699 400}%
\special{dt 0.045}%
\special{pa 3700 2000}%
\special{pa 1500 2000}%
\special{dt 0.045}%
\special{pa 1500 2000}%
\special{pa 1501 2000}%
\special{dt 0.045}%
% LINE 2 2 3 0
% 2 950 3350 2150 2150
% 
\special{pn 8}%
\special{pa 950 2950}%
\special{pa 2150 1750}%
\special{dt 0.045}%
\special{pa 2150 1750}%
\special{pa 2150 1750}%
\special{dt 0.045}%
% LINE 2 2 3 0
% 2 600 1910 1800 710
% 
\special{pn 8}%
\special{pa 600 1510}%
\special{pa 1800 310}%
\special{dt 0.045}%
\special{pa 1800 310}%
\special{pa 1800 310}%
\special{dt 0.045}%
% LINE 1 0 3 0
% 2 1400 2800 4200 4200
% 
\special{pn 13}%
\special{pa 1400 2400}%
\special{pa 4200 3800}%
\special{fp}%
% LINE 2 2 3 0
% 2 3800 4000 3800 3400
% 
\special{pn 8}%
\special{pa 3800 3600}%
\special{pa 3800 3000}%
\special{dt 0.045}%
\special{pa 3800 3000}%
\special{pa 3800 3001}%
\special{dt 0.045}%
% ELLIPSE 2 0 3 0
% 4 3800 2600 4200 3400 4200 3400 4200 3400
% 
\special{pn 8}%
\special{ar 3800 2200 400 800  0.0000000 6.2831853}%
% VECTOR 2 0 3 0
% 4 3500 2900 4100 2300 3800 3200 3800 2000
% 
\special{pn 8}%
\special{pa 3500 2500}%
\special{pa 4100 1900}%
\special{fp}%
\special{sh 1}%
\special{pa 4100 1900}%
\special{pa 4039 1933}%
\special{pa 4062 1938}%
\special{pa 4067 1961}%
\special{pa 4100 1900}%
\special{fp}%
\special{pa 3800 2800}%
\special{pa 3800 1600}%
\special{fp}%
\special{sh 1}%
\special{pa 3800 1600}%
\special{pa 3780 1667}%
\special{pa 3800 1653}%
\special{pa 3820 1667}%
\special{pa 3800 1600}%
\special{fp}%
% SPLINE 1 0 3 0
% 42 4110 2250 4097 2248 4084 2245 4072 2243 4059 2240 4048 2236 4036 2232 4024 2229 4013 2223 4002 2218 3990 2211 3980 2205 3970 2199 3960 2191 3950 2184 3942 2175 3933 2167 3923 2158 3916 2148 3909 2139 3902 2128 3895 2118 3889 2107 3883 2095 3877 2085 3872 2072 3867 2060 3864 2047 3859 2034 3855 2021 3852 2007 3851 1993 3848 1980 3846 1964 3845 1951 3843 1937 3843 1921 3844 1907 3844 1891 3845 1876 3847 1861 3847 1861
% 
\special{pn 13}%
\special{pa 4110 1850}%
\special{pa 4079 1844}%
\special{pa 4048 1836}%
\special{pa 4018 1826}%
\special{pa 3989 1811}%
\special{pa 3963 1793}%
\special{pa 3939 1772}%
\special{pa 3917 1749}%
\special{pa 3898 1723}%
\special{pa 3883 1695}%
\special{pa 3869 1666}%
\special{pa 3860 1636}%
\special{pa 3852 1605}%
\special{pa 3847 1573}%
\special{pa 3844 1541}%
\special{pa 3844 1509}%
\special{pa 3845 1477}%
\special{pa 3847 1461}%
\special{sp}%
% SPLINE 1 0 3 0
% 42 3480 2960 3493 2962 3506 2965 3518 2967 3531 2970 3542 2974 3554 2978 3566 2981 3577 2987 3588 2992 3600 2999 3610 3005 3620 3011 3630 3019 3640 3026 3648 3035 3657 3043 3667 3052 3674 3062 3681 3071 3688 3082 3695 3092 3701 3103 3707 3115 3713 3125 3718 3138 3723 3150 3726 3163 3731 3176 3735 3189 3738 3203 3739 3217 3742 3230 3744 3246 3745 3259 3747 3273 3747 3289 3746 3303 3746 3319 3745 3334 3743 3349 3743 3349
% 
\special{pn 13}%
\special{pa 3480 2560}%
\special{pa 3511 2566}%
\special{pa 3542 2574}%
\special{pa 3572 2584}%
\special{pa 3601 2599}%
\special{pa 3627 2617}%
\special{pa 3651 2638}%
\special{pa 3673 2661}%
\special{pa 3692 2687}%
\special{pa 3707 2715}%
\special{pa 3721 2744}%
\special{pa 3730 2774}%
\special{pa 3738 2805}%
\special{pa 3743 2837}%
\special{pa 3746 2869}%
\special{pa 3746 2901}%
\special{pa 3745 2933}%
\special{pa 3743 2949}%
\special{sp}%
% LINE 2 2 3 0
% 2 2600 3400 2600 2800
% 
\special{pn 8}%
\special{pa 2600 3000}%
\special{pa 2600 2400}%
\special{dt 0.045}%
\special{pa 2600 2400}%
\special{pa 2600 2401}%
\special{dt 0.045}%
% ELLIPSE 2 0 3 1
% 4 2600 2000 3000 2800 3000 2800 3000 2800
% 
\special{pn 8}%
\special{ar 2600 1600 400 800  0.0000000 6.2831853}%
% VECTOR 2 0 3 2
% 4 2300 2300 2900 1700 2600 2600 2600 1400
% 
\special{pn 8}%
\special{pa 2300 1900}%
\special{pa 2900 1300}%
\special{fp}%
\special{sh 1}%
\special{pa 2900 1300}%
\special{pa 2839 1333}%
\special{pa 2862 1338}%
\special{pa 2867 1361}%
\special{pa 2900 1300}%
\special{fp}%
\special{pa 2600 2200}%
\special{pa 2600 1000}%
\special{fp}%
\special{sh 1}%
\special{pa 2600 1000}%
\special{pa 2580 1067}%
\special{pa 2600 1053}%
\special{pa 2620 1067}%
\special{pa 2600 1000}%
\special{fp}%
% SPLINE 1 0 3 3
% 41 2900 1600 2888 1619 2876 1637 2864 1653 2852 1667 2841 1681 2830 1693 2819 1703 2809 1712 2799 1720 2789 1726 2780 1731 2770 1735 2761 1736 2753 1737 2745 1736 2736 1734 2729 1731 2722 1725 2715 1719 2707 1711 2701 1702 2696 1691 2690 1679 2685 1666 2680 1650 2674 1634 2670 1616 2666 1598 2663 1578 2659 1555 2656 1532 2655 1508 2652 1482 2650 1455 2648 1427 2648 1397 2647 1366 2646 1333 2646 1300 2646 1265
% 
\special{pn 13}%
\special{pa 2900 1200}%
\special{pa 2883 1227}%
\special{pa 2864 1253}%
\special{pa 2844 1278}%
\special{pa 2822 1301}%
\special{pa 2797 1321}%
\special{pa 2769 1335}%
\special{pa 2737 1334}%
\special{pa 2711 1316}%
\special{pa 2695 1288}%
\special{pa 2683 1259}%
\special{pa 2673 1229}%
\special{pa 2666 1197}%
\special{pa 2661 1166}%
\special{pa 2656 1134}%
\special{pa 2654 1102}%
\special{pa 2651 1070}%
\special{pa 2649 1038}%
\special{pa 2648 1006}%
\special{pa 2647 974}%
\special{pa 2646 942}%
\special{pa 2646 910}%
\special{pa 2646 878}%
\special{pa 2646 865}%
\special{sp}%
% SPLINE 1 0 3 4
% 41 2300 2410 2312 2391 2324 2373 2336 2357 2348 2343 2359 2329 2370 2317 2381 2307 2391 2298 2401 2290 2411 2284 2420 2279 2430 2275 2439 2274 2447 2273 2455 2274 2464 2276 2471 2279 2478 2285 2485 2291 2493 2299 2499 2308 2504 2319 2510 2331 2515 2344 2520 2360 2526 2376 2530 2394 2534 2412 2537 2432 2541 2455 2544 2478 2545 2502 2548 2528 2550 2555 2552 2583 2552 2613 2553 2644 2554 2677 2554 2710 2554 2745
% 
\special{pn 13}%
\special{pa 2300 2010}%
\special{pa 2317 1983}%
\special{pa 2336 1957}%
\special{pa 2356 1932}%
\special{pa 2378 1909}%
\special{pa 2403 1889}%
\special{pa 2431 1875}%
\special{pa 2463 1876}%
\special{pa 2489 1894}%
\special{pa 2505 1922}%
\special{pa 2517 1951}%
\special{pa 2527 1981}%
\special{pa 2534 2013}%
\special{pa 2539 2044}%
\special{pa 2544 2076}%
\special{pa 2546 2108}%
\special{pa 2549 2140}%
\special{pa 2551 2172}%
\special{pa 2552 2204}%
\special{pa 2553 2236}%
\special{pa 2554 2268}%
\special{pa 2554 2300}%
\special{pa 2554 2332}%
\special{pa 2554 2345}%
\special{sp}%
% LINE 2 2 3 0
% 2 4110 1930 1710 730
% 
\special{pn 8}%
\special{pa 4110 1530}%
\special{pa 1710 330}%
\special{dt 0.045}%
\special{pa 1710 330}%
\special{pa 1711 330}%
\special{dt 0.045}%
% LINE 2 2 3 0
% 2 3940 3500 1540 2300
% 
\special{pn 8}%
\special{pa 3940 3100}%
\special{pa 1540 1900}%
\special{dt 0.045}%
\special{pa 1540 1900}%
\special{pa 1541 1900}%
\special{dt 0.045}%
% STR 2 0 3 0
% 3 4330 2090 4330 2190 2 0
% a1
\put(43.3000,-17.9000){\makebox(0,0)[lb]{$A_1$}}%
% STR 2 0 3 0
% 3 3990 3920 3990 4020 2 0
% k-k=0
\put(39.9000,-36.2000){\makebox(0,0)[lb]{\large $\kappa_0-\kappa_\infty=0$}}%
\end{picture}%
\end{center}
\caption{A Confluence of Nodal Curves in  the case $\tilde{E_6}$ ($P_{IV}$).}
\label{fig:p4}
\end{figure}

Let us consider the neighborhood of $(\kappa_0, \kappa_{\infty}, t) = (0, 0, t) \in \cM_R \times \cB_R$ and 
the hyperplanes as in (\ref{eq:e6-h}).  
Then,  over the subvariety  $H_0 \cap H_{\infty} = \{ ( 0, 0, t) \}$,  the family 
$(\cS - \cD)_{0, 0}$ contains both of families of nodal curves $\cC_{0} \cup \cC_{\infty}$ ($A_2$-configuration), (see Figure \ref{fig:p4} ).  
We call this phenomenon {\em the confluence of nodal curves} of 
Okamoto--Painlev\'e pairs.

Besides hyperplanes $H_0$, $ H_{\infty}$, we also have the hyperplane 
$$
H_{\kappa_0 = \kappa_{\infty}} = \{  \kappa_0 = \kappa_{\infty} \}.   
$$
Then one can easily see that over hyperplane 
$H_{\kappa_0 = \kappa_{\infty}}$ there exists a family 
of $(-2)$-curves defined by 
$$
\cC_{\kappa_0 = \kappa_{\infty}} \cap \tilde{U_0}_{\kappa_0 = \kappa_{\infty}} = \{ x_0 y_0 -  \kappa_0 = 0 \}.
$$
Note that if $\kappa_0$ goes to $0$, then the defining 
equation of the family becomes $x_0 y_0 = 0$.  Therefore on $(\cS - \cD)_{0, 0}$ we have 
a homological relation:
$$
\cC_{\kappa_0 = \kappa_{\infty}} = \cC_{0} \cup \cC_{\infty}.  
$$
(See Figure \ref{fig:p4}).  
On $\cC_{\kappa_0 = \kappa_{\infty}}$, 
the Painlev\'e system (\ref{eq:p4}) can be reduced to
\begin{eqnarray}
\displaystyle{\frac{d x_0 }{d t}} &  = & -x_0^2 -2 t x_0 + 2 \kappa_0, \label{eq:rix}\\
\displaystyle{\frac{d y_0 }{d t}} & = & -2y_0^2 +2 t y_0 + \kappa_0, \label{eq:riy}.
\end{eqnarray}
Note that if $\kappa_0 \not= 0$ 
the equations (\ref{eq:rix}) and (\ref{eq:riy}) can be transformed to each other 
by the coordinate change $x_0 = \kappa_0/y_0$.  

The hyperplanes are reflection hyperplanes in $\cM_{R}$ with respect to 
the reflections of the affine Weyl group $W(\tilde{A_2})$, which acts on both 
$\cM_{R}$ or $\cS$ as B\"acklund transformations (cf. \cite{U-W1} and \cite{NTY}).  
For example, by B\"acklund transformations, the Riccati equations (\ref{eq:p4-0}), (\ref{eq:p4-inf}) and 
(\ref{eq:rix}) are birational equivalent to each other. See Theorem 3.3 in \cite{U-W1}.

 \subsection{Rational solutions}

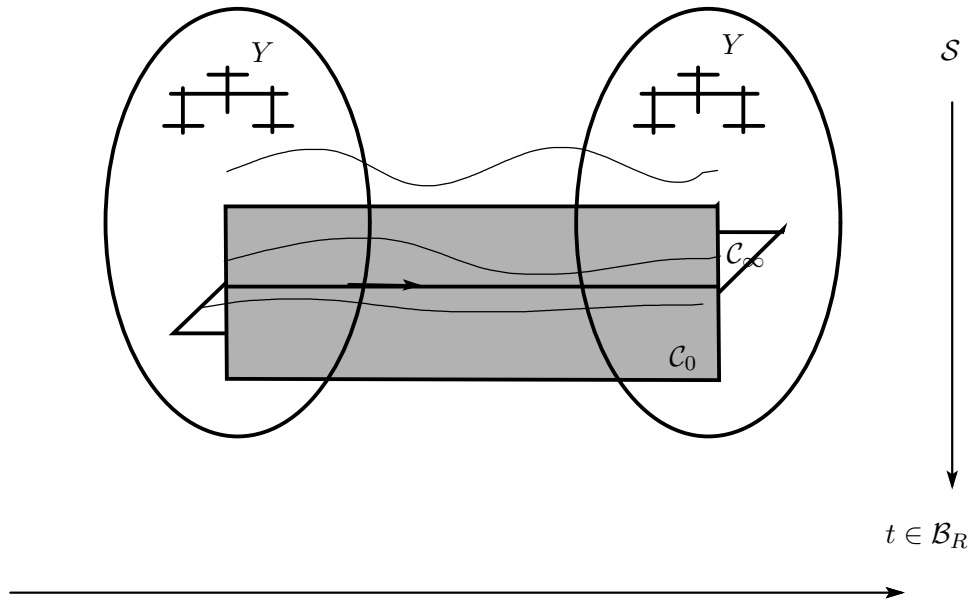
\begin{figure}
\begin{center}
%WinTpicVersion2.15
\unitlength 0.1in
\begin{picture}(49.35,30.64)(2.70,-37.44)
% LINE 0 0 3 0
% 2 1110 1526 1718 1526
% 
\special{pn 20}%
\special{pa 1110 1126}%
\special{pa 1718 1126}%
\special{fp}%
% LINE 0 0 3 0
% 2 1408 1388 1408 1623
% 
\special{pn 20}%
\special{pa 1408 988}%
\special{pa 1408 1223}%
\special{fp}%
% LINE 0 0 3 0
% 2 1174 1496 1174 1731
% 
\special{pn 20}%
\special{pa 1174 1096}%
\special{pa 1174 1331}%
\special{fp}%
% LINE 0 0 3 0
% 2 1307 1425 1510 1425
% 
\special{pn 20}%
\special{pa 1307 1025}%
\special{pa 1510 1025}%
\special{fp}%
% LINE 0 0 3 0
% 2 1075 1694 1278 1694
% 
\special{pn 20}%
\special{pa 1075 1294}%
\special{pa 1278 1294}%
\special{fp}%
% LINE 0 0 3 0
% 2 1545 1694 1748 1694
% 
\special{pn 20}%
\special{pa 1545 1294}%
\special{pa 1748 1294}%
\special{fp}%
% LINE 0 0 3 0
% 2 1643 1492 1643 1728
% 
\special{pn 20}%
\special{pa 1643 1092}%
\special{pa 1643 1328}%
\special{fp}%
% POLYGON 0 0 3 0
% 11 1677 2249 1124 2781 3763 2781 4316 2235 4316 2235 4309 2256 4309 2256 4310 2250 4310 2250 4310 2250 1677 2249
% 
\special{pn 20}%
\special{pa 1677 1849}%
\special{pa 1124 2381}%
\special{pa 3763 2381}%
\special{pa 4316 1835}%
\special{pa 4316 1835}%
\special{pa 4309 1856}%
\special{pa 4309 1856}%
\special{pa 4310 1850}%
\special{pa 4310 1850}%
\special{pa 4310 1850}%
\special{pa 1677 1849}%
\special{fp}%
% VECTOR 1 0 3 0
% 4 270 4139 4932 4139 4925 4139 4932 4139
% 
\special{pn 13}%
\special{pa 270 3739}%
\special{pa 4932 3739}%
\special{fp}%
\special{sh 1}%
\special{pa 4932 3739}%
\special{pa 4865 3719}%
\special{pa 4879 3739}%
\special{pa 4865 3759}%
\special{pa 4932 3739}%
\special{fp}%
\special{pa 4925 3739}%
\special{pa 4932 3739}%
\special{fp}%
\special{sh 1}%
\special{pa 4932 3739}%
\special{pa 4865 3719}%
\special{pa 4879 3739}%
\special{pa 4865 3759}%
\special{pa 4932 3739}%
\special{fp}%
% POLYGON 0 0 1 0
% 8 1397 2116 1404 3019 3980 3026 3973 2109 3973 2109 3973 2109 3966 2116 1397 2116
% 
\special{pn 20}%
\special{sh 0.300}%
\special{pa 1397 1716}%
\special{pa 1404 2619}%
\special{pa 3980 2626}%
\special{pa 3973 1709}%
\special{pa 3973 1709}%
\special{pa 3973 1709}%
\special{pa 3966 1716}%
\special{pa 1397 1716}%
\special{fp}%
% LINE 0 0 3 0
% 2 1404 2536 3973 2536
% 
\special{pn 20}%
\special{pa 1404 2136}%
\special{pa 3973 2136}%
\special{fp}%
% STR 2 0 3 0
% 3 3720 2880 3720 2950 2 0
% $\cC_{0}$
\put(37.2000,-25.5000){\makebox(0,0)[lb]{$\cC_{0}$}}%
% STR 2 0 3 0
% 3 4020 2340 4020 2410 2 0
% $\cC_{\infty}$
\put(40.2000,-20.1000){\makebox(0,0)[lb]{$\cC_{\infty}$}}%
% VECTOR 1 0 3 0
% 2 2041 2522 2405 2529
% 
\special{pn 13}%
\special{pa 2041 2122}%
\special{pa 2405 2129}%
\special{fp}%
\special{sh 1}%
\special{pa 2405 2129}%
\special{pa 2339 2108}%
\special{pa 2352 2128}%
\special{pa 2338 2148}%
\special{pa 2405 2129}%
\special{fp}%
% STR 2 0 3 0
% 3 4850 3810 4850 3880 2 0
% $t \in \cB_R$
\put(48.5000,-34.8000){\makebox(0,0)[lb]{$t \in \cB_R$}}%
% ELLIPSE 0 0 3 0
% 4 1460 2200 767 3320 767 3320 676 3453
% 
\special{pn 20}%
\special{ar 1460 1800 693 1120  2.3619674 6.2831853}%
\special{ar 1460 1800 693 1120  0.0000000 2.3561945}%
% LINE 0 0 3 0
% 2 3574 1526 4182 1526
% 
\special{pn 20}%
\special{pa 3574 1126}%
\special{pa 4182 1126}%
\special{fp}%
% LINE 0 0 3 0
% 2 3871 1388 3871 1623
% 
\special{pn 20}%
\special{pa 3871 988}%
\special{pa 3871 1223}%
\special{fp}%
% LINE 0 0 3 0
% 2 3638 1496 3638 1731
% 
\special{pn 20}%
\special{pa 3638 1096}%
\special{pa 3638 1331}%
\special{fp}%
% LINE 0 0 3 0
% 2 3771 1425 3974 1425
% 
\special{pn 20}%
\special{pa 3771 1025}%
\special{pa 3974 1025}%
\special{fp}%
% LINE 0 0 3 0
% 2 3539 1694 3742 1694
% 
\special{pn 20}%
\special{pa 3539 1294}%
\special{pa 3742 1294}%
\special{fp}%
% LINE 0 0 3 0
% 2 4009 1694 4212 1694
% 
\special{pn 20}%
\special{pa 4009 1294}%
\special{pa 4212 1294}%
\special{fp}%
% LINE 0 0 3 0
% 2 4107 1492 4107 1728
% 
\special{pn 20}%
\special{pa 4107 1092}%
\special{pa 4107 1328}%
\special{fp}%
% ELLIPSE 0 0 3 0
% 4 3924 2200 3231 3320 3231 3320 3140 3453
% 
\special{pn 20}%
\special{ar 3924 1800 693 1120  2.3619674 6.2831853}%
\special{ar 3924 1800 693 1120  0.0000000 2.3561945}%
% STR 2 0 3 0
% 3 1530 1290 1530 1360 2 0
% $Y$
\put(15.3000,-9.6000){\makebox(0,0)[lb]{$Y$}}%
% STR 2 0 3 0
% 3 3994 1248 3994 1318 2 0
% $Y$
\put(39.9400,-9.1800){\makebox(0,0)[lb]{$Y$}}%
% STR 2 0 3 0
% 3 5140 1280 5140 1350 2 0
% $\cS$
\put(51.4000,-9.5000){\makebox(0,0)[lb]{$\cS$}}%
% VECTOR 1 0 3 0
% 2 5200 1570 5200 3579
% 
\special{pn 13}%
\special{pa 5200 1170}%
\special{pa 5200 3179}%
\special{fp}%
\special{sh 1}%
\special{pa 5200 3179}%
\special{pa 5220 3112}%
\special{pa 5200 3126}%
\special{pa 5180 3112}%
\special{pa 5200 3179}%
\special{fp}%
% SPLINE 2 0 3 0
% 5 1390 2403 2349 2291 2811 2452 3980 2375 3987 2375
% 
\special{pn 8}%
\special{pa 1390 2003}%
\special{pa 1422 1994}%
\special{pa 1455 1984}%
\special{pa 1487 1975}%
\special{pa 1520 1966}%
\special{pa 1552 1957}%
\special{pa 1584 1948}%
\special{pa 1617 1939}%
\special{pa 1649 1931}%
\special{pa 1681 1923}%
\special{pa 1714 1915}%
\special{pa 1746 1908}%
\special{pa 1778 1901}%
\special{pa 1810 1894}%
\special{pa 1842 1888}%
\special{pa 1874 1883}%
\special{pa 1906 1878}%
\special{pa 1938 1874}%
\special{pa 1970 1870}%
\special{pa 2001 1867}%
\special{pa 2033 1865}%
\special{pa 2064 1863}%
\special{pa 2096 1863}%
\special{pa 2127 1863}%
\special{pa 2158 1864}%
\special{pa 2190 1865}%
\special{pa 2221 1868}%
\special{pa 2252 1872}%
\special{pa 2282 1877}%
\special{pa 2313 1883}%
\special{pa 2344 1890}%
\special{pa 2374 1898}%
\special{pa 2405 1907}%
\special{pa 2435 1917}%
\special{pa 2465 1927}%
\special{pa 2495 1939}%
\special{pa 2525 1950}%
\special{pa 2555 1962}%
\special{pa 2585 1974}%
\special{pa 2615 1986}%
\special{pa 2645 1998}%
\special{pa 2676 2009}%
\special{pa 2706 2020}%
\special{pa 2736 2030}%
\special{pa 2766 2040}%
\special{pa 2797 2048}%
\special{pa 2828 2056}%
\special{pa 2858 2062}%
\special{pa 2889 2067}%
\special{pa 2920 2072}%
\special{pa 2951 2075}%
\special{pa 2983 2077}%
\special{pa 3014 2079}%
\special{pa 3045 2079}%
\special{pa 3077 2079}%
\special{pa 3109 2078}%
\special{pa 3140 2077}%
\special{pa 3172 2075}%
\special{pa 3204 2072}%
\special{pa 3236 2069}%
\special{pa 3268 2065}%
\special{pa 3300 2061}%
\special{pa 3332 2057}%
\special{pa 3364 2052}%
\special{pa 3397 2047}%
\special{pa 3429 2042}%
\special{pa 3461 2036}%
\special{pa 3494 2031}%
\special{pa 3526 2026}%
\special{pa 3558 2020}%
\special{pa 3591 2015}%
\special{pa 3623 2010}%
\special{pa 3655 2005}%
\special{pa 3688 2000}%
\special{pa 3720 1995}%
\special{pa 3752 1991}%
\special{pa 3785 1987}%
\special{pa 3817 1984}%
\special{pa 3849 1981}%
\special{pa 3881 1978}%
\special{pa 3914 1977}%
\special{pa 3946 1975}%
\special{pa 3978 1975}%
\special{pa 3987 1975}%
\special{sp}%
% SPLINE 2 0 3 0
% 6 1257 2648 1712 2592 2601 2683 3511 2648 3896 2627 3896 2627
% 
\special{pn 8}%
\special{pa 1257 2248}%
\special{pa 1289 2243}%
\special{pa 1321 2237}%
\special{pa 1352 2232}%
\special{pa 1384 2227}%
\special{pa 1416 2222}%
\special{pa 1448 2217}%
\special{pa 1479 2212}%
\special{pa 1511 2208}%
\special{pa 1543 2204}%
\special{pa 1575 2201}%
\special{pa 1606 2198}%
\special{pa 1638 2196}%
\special{pa 1670 2194}%
\special{pa 1702 2192}%
\special{pa 1733 2192}%
\special{pa 1765 2192}%
\special{pa 1797 2192}%
\special{pa 1829 2193}%
\special{pa 1861 2195}%
\special{pa 1892 2197}%
\special{pa 1924 2199}%
\special{pa 1956 2202}%
\special{pa 1988 2205}%
\special{pa 2020 2209}%
\special{pa 2051 2213}%
\special{pa 2083 2217}%
\special{pa 2115 2221}%
\special{pa 2147 2225}%
\special{pa 2179 2230}%
\special{pa 2210 2235}%
\special{pa 2242 2239}%
\special{pa 2274 2244}%
\special{pa 2306 2248}%
\special{pa 2338 2253}%
\special{pa 2370 2258}%
\special{pa 2402 2262}%
\special{pa 2433 2266}%
\special{pa 2465 2270}%
\special{pa 2497 2273}%
\special{pa 2529 2277}%
\special{pa 2561 2280}%
\special{pa 2593 2282}%
\special{pa 2625 2285}%
\special{pa 2657 2286}%
\special{pa 2689 2288}%
\special{pa 2721 2289}%
\special{pa 2753 2289}%
\special{pa 2785 2290}%
\special{pa 2817 2290}%
\special{pa 2849 2289}%
\special{pa 2880 2289}%
\special{pa 2912 2288}%
\special{pa 2944 2287}%
\special{pa 2976 2286}%
\special{pa 3008 2284}%
\special{pa 3040 2282}%
\special{pa 3072 2280}%
\special{pa 3104 2278}%
\special{pa 3136 2276}%
\special{pa 3168 2274}%
\special{pa 3200 2272}%
\special{pa 3232 2269}%
\special{pa 3264 2267}%
\special{pa 3296 2264}%
\special{pa 3328 2262}%
\special{pa 3360 2259}%
\special{pa 3392 2257}%
\special{pa 3424 2254}%
\special{pa 3456 2252}%
\special{pa 3488 2250}%
\special{pa 3520 2247}%
\special{pa 3552 2245}%
\special{pa 3584 2243}%
\special{pa 3616 2241}%
\special{pa 3648 2239}%
\special{pa 3680 2238}%
\special{pa 3712 2236}%
\special{pa 3744 2234}%
\special{pa 3776 2233}%
\special{pa 3808 2231}%
\special{pa 3840 2230}%
\special{pa 3872 2228}%
\special{pa 3896 2227}%
\special{sp}%
% SPLINE 2 0 3 0
% 6 1404 1934 1992 1815 2433 2004 3280 1787 3966 1934 3973 1927
% 
\special{pn 8}%
\special{pa 1404 1534}%
\special{pa 1436 1520}%
\special{pa 1468 1507}%
\special{pa 1500 1494}%
\special{pa 1532 1481}%
\special{pa 1563 1469}%
\special{pa 1595 1457}%
\special{pa 1627 1446}%
\special{pa 1659 1436}%
\special{pa 1690 1427}%
\special{pa 1722 1419}%
\special{pa 1753 1412}%
\special{pa 1784 1407}%
\special{pa 1816 1403}%
\special{pa 1847 1401}%
\special{pa 1878 1400}%
\special{pa 1908 1401}%
\special{pa 1939 1404}%
\special{pa 1969 1410}%
\special{pa 2000 1417}%
\special{pa 2030 1427}%
\special{pa 2060 1438}%
\special{pa 2089 1451}%
\special{pa 2119 1466}%
\special{pa 2149 1481}%
\special{pa 2178 1497}%
\special{pa 2207 1512}%
\special{pa 2236 1528}%
\special{pa 2266 1543}%
\special{pa 2295 1558}%
\special{pa 2324 1571}%
\special{pa 2353 1583}%
\special{pa 2382 1592}%
\special{pa 2411 1600}%
\special{pa 2441 1605}%
\special{pa 2470 1607}%
\special{pa 2499 1607}%
\special{pa 2528 1605}%
\special{pa 2558 1600}%
\special{pa 2587 1594}%
\special{pa 2617 1586}%
\special{pa 2647 1576}%
\special{pa 2677 1565}%
\special{pa 2707 1553}%
\special{pa 2737 1540}%
\special{pa 2767 1526}%
\special{pa 2797 1512}%
\special{pa 2828 1498}%
\special{pa 2859 1484}%
\special{pa 2890 1469}%
\special{pa 2921 1456}%
\special{pa 2952 1442}%
\special{pa 2984 1430}%
\special{pa 3015 1418}%
\special{pa 3047 1408}%
\special{pa 3079 1399}%
\special{pa 3112 1392}%
\special{pa 3144 1386}%
\special{pa 3177 1383}%
\special{pa 3211 1382}%
\special{pa 3244 1383}%
\special{pa 3278 1387}%
\special{pa 3312 1393}%
\special{pa 3346 1402}%
\special{pa 3381 1414}%
\special{pa 3415 1427}%
\special{pa 3450 1441}%
\special{pa 3484 1456}%
\special{pa 3518 1472}%
\special{pa 3552 1489}%
\special{pa 3586 1505}%
\special{pa 3620 1521}%
\special{pa 3653 1535}%
\special{pa 3685 1549}%
\special{pa 3717 1561}%
\special{pa 3748 1571}%
\special{pa 3779 1579}%
\special{pa 3809 1583}%
\special{pa 3838 1585}%
\special{pa 3866 1583}%
\special{pa 3892 1578}%
\special{pa 3918 1568}%
\special{pa 3943 1553}%
\special{pa 3966 1534}%
\special{pa 3973 1527}%
\special{sp}%
\end{picture}%
\end{center}
\caption{Rational solution coming from $\cC_0 \cap \cC_{\infty}$ for  $\tilde{E_6}$ ($P_{IV}$).}
\label{fig:rat-p4}
\end{figure}

We shall remark briefly on rational solutions of Painlev\'e equations.  
In the above example, when $(\kappa_0, \kappa_{\infty}) = (0, 0)$, the functions 
\begin{equation}
(x_0,  y_0)  \equiv (0, 0) 
\end{equation}
give a solution of the system (\ref{eq:p4}), hence gives a rational solution for the Painlev\'e equation $P_{IV}$.  
From the view point of the geometry of Okamoto--Painlev\'e pairs, it is clear that the intersection of two different families of  nodal curves $\cC_{0}, \cC_{\infty}$  
gives a solution of the 
 Painlev\'e equation.  
In fact, Painlev\'e vector field $\tilde{v}$ in (\ref{eq:vf}) is tangent to each family 
of rational curves by Proposition \ref{prop:riccati}, hence tangent to their intersection. 
(See Figure \ref{fig:rat-p4}).  It is not surprising that not all rational solutions of 
Painlev\'e equations can  be obtained in this way.  For example, as we explained after 
Theorem \ref{thm:u-w}, the equation  $S_{II}(0) $ in (\ref{eq:p-2}) has the rational 
solution $(x_0, y_0) = (0, \frac{t}{2})$, but no Riccati solution.   
It should be  an interesting problem to understand the rational or algebraic solutions 
from the view point of the geometry of Okamoto--Painlev\'e pairs.  

Here we only remark that there are many works for the classification problems of 
rational and algebraic solutions.  
$($See e.g., \cite{DM}, \cite{Maz}, \cite{Mu1}, \cite{NO}, \cite{O3}, \cite{U-W1}, \cite{U-W2}$)$.

\newpage

%\vspace{2cm}

\appendix

\section{\bf Local cohomology group $H^1_D(\Theta_S(-\log D))$}\label{sec:lcoh}

Let $(S,Y)$ be a rational Okamoto--Painlev\'e pair of non-fibered type and of additive type which corresponds to Painlev\'e equations (i.e. of type $\tilde D_i  (4 \leq i \leq 8)$ or $\tilde E_i (6 \leq i \leq 8)$), and set $D=Y_{red}$.

Applying  the classification of nodal curves on $S-D$, we will 
investigate the local cohomology group $H^1_D(\Theta_S(-\log D))$.   
Note that the local cohomology group 
can be regarded  as the space of time variables for differential equations
associated to  $(S,D)$ (cf. \S 3. \cite{STT}). 

\vskp

We state our conjecture for the local cohomology:
\begin{Conjecture}[Conjecture 3.1. \cite{STT}, \cite{T}]\label{conj:local}
 Let $(S, Y)$ be  a  rational Okamoto-Painlev\'e pair $(S, Y)$ as above.  Then  we have 
\begin{equation}
H^1_D(\Theta_S(-\log D) ) \simeq \C.
\end{equation}
\end{Conjecture}

For the positivity of the dimension of the cohomology group,  we have the following result:
\begin{Theorem}[Theorem 2.1. \cite{T}]\label{thm:time}
% Let $(S, Y)$ be  a generalized rational Okamoto-Painlev\'e pair $(S, Y)$ with the conditions above. Then  we have 
\begin{equation}
\dim H^0 (D, \Theta_S(- \log D) \otimes N_{D} ) = 1.
\end{equation}
Here  we put $N_D = {\cal O}_S(D)/{\cal O}_S$.  

In particular, a natural inclusion 
$$
H^0 (D, \Theta_S(- \log D) \otimes N_{D} ) \hookrightarrow H^1_D(\Theta_S(-\log D) ), 
$$ 
implies
\begin{equation}\label{eq:time}
\dim H^1_D(\Theta_S(-\log D) ) \geq  1.
\end{equation}
\end{Theorem}

On the other hand, in this section,  we shall prove

\begin{Theorem}\label{thm:appendix}
Let
$$
\begin{array}{ccc}
 \cS &  \hookleftarrow  & \cD  \\
\pi \downarrow & \swarrow & \varphi \\
 \cM_R \times \cB_R  
\end{array} 
$$
be the 
 semi-universal deformation of rational Okamoto--Painlev\'e pairs $(S,D)$ whose type is 
 one of $\tilde{E_8}, \tilde{E_7},\tilde{D_8}, \tilde{D_6},\tilde{E_6},\tilde{D_5}$ and $\tilde{D_4}$ (i.e., 
 except for $R = \tilde{D_7}$).
Then there is a Zariski open set $U \subset \cM_R \times \cB_R$ such that for any $(\balpha,t) \in U$,
$$
\dim H^1_{\cD_{(\balpha,t)}}(\Theta_{\cS_{(\balpha,t)}}( - \log \cD_{(\balpha,t)})=1. 
$$
\end{Theorem}

\begin{Remark}\label{rem:d8-e8}
{\rm
For $(S, Y)$ of type $\tilde D_8$ or $\tilde E_8$,  Theorem  \ref{thm:time} proves Conjecture \ref{conj:local}. 
In fact,  we always have the 
inclusion $ H^1_D(\Theta_S( - \log D)) \hookrightarrow  H^1(S, \Theta_S( - \log D))$ and $\dim H^1(S, \Theta_S( - \log D)) =10-9=1$ for these cases.
}
\end{Remark}

From Remark \ref{rem:d8-e8}, 
in order to show Theorem \ref{thm:appendix}, we will estimate the dimension of the local cohomology 
group for a special rational Okamoto--Painlev\'e pairs of other type $R$.

\vskp

We first calculate some cohomology groups.

\begin{Lemma}\label{lem:h2tlC}
Let $(S,Y)$ be a  rational Okamoto--Painlev\'e pair, and $C$  a normal crossing divisor of $S$.
Moreover, let $C=\sum_{i=1}^s C_i$ be an irreducible decomposition of $C$, and we assume that $\{C_i\}_{i=1}^s$ is linearly independent in $H^2(S,\C) \simeq \Pic (S) \otimes \C$.
Then we have
$$
H^2(S,\Theta_S(-\log C))=\{ 0 \}.
$$
\end{Lemma}

{\it Proof.}
We have only to replace $D$ of [Lemma 2.2 and Corollary 2.1,\cite{STT}] with $C$.
\qed

\begin{Lemma}\label{lem:h1tlDC}
Let $(S,Y)$ be a generalized rational Okamoto--Painlev\'e pair such that $D=Y_{red} = \sum_{i=1}^r Y_i$ 
is a normal crossing divisor with at least two irreducible components, 
say $r \geq 2$, and let $C=\sum_{i=1}^s C_i$ be a normal crossing divisor of $S$.
We assume that
\begin{enumerate}
\item\label{as:c0} $C \subset S-D$,
\item $C_i \simeq \BP^1$,
\item $\{Y_i , C_j | 1 \le i \le r, 1 \le j \le s  \}$ is linearly independent.
\end{enumerate}
Then we have
$$
\dim H^1(S,\Theta_S(-\log (D+C)))=10-(r+s).
$$
\end{Lemma}

{\it Proof.}
Note that assumption \ref{as:c0} implies $D+C$ is normal crossing and $K_S \cdot C_i = -Y \cdot C_i=0$. We have $H^2(S,\Theta_S(-\log (D+C)))=0$ by applying Lemma \ref{lem:h2tlC} to $D+C$. Therefore by using the same argument as Proposition 2.2 in \cite{STT}, we have the assertion.
\qed

\begin{Remark}{\em
We have the following exact sequence of sheaves:
$$
0 \ra \Theta_S(-\log (D+C)) \ra \Theta_S(-\log (D+C-C_i)) \ra  N_{C_i /S} \ra 0
$$
where $N_{C_i /S}=\cO_S(C_i) /\cO_S$ denotes the normal bundle of the divisor $C_i \subset S$. Note that since $N_{C_i /S}=\cO_{C_i}(-2)$, we have $H^0(N_{C_i /S})=\{ 0 \}$. Then the morphism
$$
H^0(\Theta_S(-\log (D+C))) \ra H^0(\Theta_S(-\log (D+C-C_i)))
$$
is injective. Moreover we have $\dim H^0(\Theta_S(-\log (D+C-C_i))) - \dim H^0(\Theta_S(-\log (D+C))) =1$ by Lemma \ref{lem:h1tlDC}. 
This implies that there exist a deformation $(S',D')$ of $(S,D)$ such that only the curve $C_i$ vanish and other nodal curves remain.
}
\end{Remark}

Now we obtain the following

\begin{Proposition}\label{prop:h1DtlD}
Let $(S,Y)$ be a  rational Okamoto--Painlev\'e pair ``of non-fibered type '' such that $D=Y_{red}$ is a normal crossing divisor with at least two irreducible components, say $r \geq 2$. We suppose the existence of a divisor $C=\sum_{i=1}^{9-r} C_i$ of $S$ satisfying the conditions in Lemma \ref{lem:h1tlDC}.
Then we have
$$
\dim H^1_D(\Theta_S(-\log D)) \le 1.
$$
\end{Proposition}

{\it Proof.}
Let us consider the following exact sequence of local cohomology groups (cf. [Corollary 1.9,\cite{Gr}])
$$
H^0(S-D,\Theta_S(-\log (D+C))) \ra H^1_D(\Theta_S(-\log (D+C))) \ra H^1(S,\Theta_S(-\log (D+C))).
$$

We have an inclusion $H^0(S-D,\Theta_S(-\log (D+C))) \hookrightarrow H^0(S-D,\Theta_S(-\log D))=H^0(S-D,\Theta_S)$. Since $(S,Y)$ is of non-fibered type, from (2) of Proposition 2.1 in \cite{STT}, we have $H^0(S-D,\Theta_S)=\{ 0 \}$. Therefore we have 
$$
H^0(S-D,\Theta_S(-\log (D+C)))=\{ 0 \}.
$$
By applying Lemma \ref{lem:h1tlDC}, we see 
$$
H^1(S,\Theta_S(-\log (D+C))) \simeq \C.
$$
Moreover, since $C \subset S-D$, we have $H^1_D(\Theta_S(-\log D)) \simeq H^1_D(\Theta_S(-\log (D+C)))$, which proves the assertion.

\qed

\begin{Lemma}\label{lem:existC}
For the types $\tilde D_4,\tilde D_5,\tilde D_6,\tilde E_7$ and $\tilde E_8$, 
there exists a rational Okamoto--Painlev\'e pair $(S,Y)$ 
of non-fibered type satisfying the assumption of Proposition \ref{prop:h1DtlD}.
\end{Lemma}

{\it Proof.}  For each case, we only have to show the existence of nodal curves $ C_j \subset S - D$  
$j = 1, \cdots, 9 - r$  on  a rational Okamoto--Painlev\'e pair $(S, Y)$ of non-fibered type.  The existence 
of $(-2)$-curves follows from  Theorem \ref{thm:main-2}.

\begin{Remark}\label{rem:h1DtlD} {\rm
For any  rational Okamoto--Painlev\'e pair $(S,Y)$ of $\tilde D_7$, 
there is no $(-2)$-curve  $C$ on  $S - D$ 
satisfying the condition in Lemma \ref{lem:h1tlDC} 
(cf. Table \ref{tab:subE8fibered}, \ref{tab:config-2nonfibered}).
}
\end{Remark}

Lemma \ref{lem:existC} and Theorem \ref{thm:time} lead us the following corollary, which 
also implies Theorem \ref{thm:appendix}.

\begin{Corollary}
For the types $\tilde D_4,\tilde D_5,\tilde D_6, \tilde{D_8}, \tilde E_7$ and $\tilde E_8$, 
there exists a rational Okamoto--Painlev\'e pair $(S,Y)$ of non-fibered type such that
$$
\dim H^1_D(\Theta_S(-\log D)) = 1.
$$
\end{Corollary}

\vspace{1cm}
\begin{center}
{\bf Acknowledgements}
\end{center}
\vspace{0.5cm}
We would like to thank  Tetsu Masuda, Masatoshi Noumi, Kyoichi Takano,  
Yasuhiko Yamada and Kota Yoshioka for valuable discussions 
during the preparation of this paper.

\vspace{2cm}

\end{document}